\documentclass[11pt]{amsart}

\usepackage{graphicx}
\usepackage{float}
\usepackage{fancyhdr}
\usepackage{latexsym,amsmath,amsfonts,amscd, amsthm}
\usepackage{changebar}
\usepackage{color}
\usepackage{bm}
\usepackage{tikz}
\usepackage{multirow}
\usepackage{epsfig,subfigure}
\usepackage{enumerate}
\usepackage{multirow}
\usetikzlibrary{arrows,backgrounds,snakes}
\usepackage{xcolor}
\usepackage[margin=1in]{geometry}
\usepackage{ulem}

\makeatletter
\newcommand{\figcaption}{\def\@captype{figure}\caption}
\newcommand{\tabcaption}{\def\@captype{table}\caption}
\makeatother

\theoremstyle{definition}

\theoremstyle{remark}

\numberwithin{equation}{section}

\begin{document}

\title[Non-uniqueness for compressible Euler equations]{Numerical study of non-uniqueness for 2D
compressible isentropic Euler equations}

\author{Alberto Bressan, Yi Jiang and Hailiang Liu}
\address{Mathematics Department, The Pennsylvania State University, University Park, PA 16802, USA}
\email{axb62@psu.edu}
\address{Department of Mathematics and Statistics, Southern Illinois University Edwardsville, Edwardsville, IL 62026, USA.}
\email{yjianaa@siue.edu}
\address{Mathematics Department, Iowa State University, Ames, IA 50011, USA} 
\email{hliu@iastate.edu}


\maketitle

\begin{abstract}
In this paper, we numerically study a class of solutions with spiraling singularities in vorticity for two-dimensional, inviscid, compressible Euler systems, where the initial data have an algebraic singularity in vorticity at the origin. These are different from the multi-dimensional Riemann problems widely studied in the literature. Our computations provide numerical evidence of the existence of initial value problems with multiple solutions, thus revealing a fundamental obstruction toward the well-posedness of the governing equations. The compressible Euler equations are solved using the positivity-preserving discontinuous Galerkin method.
\end{abstract}


\section{Introduction} For strictly hyperbolic conservation laws in one space dimension, the
existence and uniqueness of entropy weak BV solutions is  well established \cite{BB05, Br00, BLY99,Da16, Gl65, HR}. Whether a similar theory can be achieved for multi-dimensional problems has 
remained an outstanding open question. On the positive side,
a wealth of results has been obtained for specific classes of problems, such as 
shock reflections \cite{CF}, or multidimensional Riemann problems:
see for example \cite{CF, Zheng} and references therein.
In many of these cases, a unique solution could be explicitly found. 
On the other hand, recent work by De Lellis, Sz\'ekelyhidi, and collaborators
\cite{LS09, LS10}  has shown the existence of a huge family of weak solutions to the Euler equations, all with the same initial data. Their construction, based on convex integration and a Baire category argument, produces an infinite family of solutions of turbulent nature, none of which can be explicitly described. 
As shown in \cite{CLK15}, in this setting the usual entropy 
admissibility conditions, imposed on weak solutions to conservation laws,
fail to select a unique solution.

At this stage, it seems unlikely that some new, physically meaningnful criteria
can be found, leading to the well posedness of the multidimensional equations.
On the contrary, simple examples of initial data, apparently leading to two distinct solutions,
has been recently studied in \cite{BS20}, for the incompressible 
two-dimensional Euler flow.
In a vorticity formulation, these equations can be written as
\begin{align}\label{eq:2DinEuler}
\begin{cases}
& \omega_t +\nabla^\bot \psi \cdot  \nabla \omega~=~0, \\
& \Delta \psi~=~\omega,
\end{cases}
\end{align}
where $\psi$ is the stream function, $u,v$ are two components of the velocity,
$\omega={\rm curl} (u, v)$ is the vorticity,  and $\nabla^\bot=(-\partial_y, \partial_x)$.  
As initial data, one takes a vorticity concentrated on two wedges, 
symmetric w.r.t.~the origin.  Inside these wedges, 
the vorticity is unbounded, with a singularity at the origin.
By approximating the same initial data in two different ways, numerical simulations
performed by Wen Shen \cite{ShenCode} show that two very different limits are obtained.

A natural question, which we investigate in the present paper, is whether
similar examples of non-uniqueness can occur also for compressible flow.
We focus on the two-dimensional isentropic compressible Euler equations:
\begin{align}\label{eq:2DEuler}
\begin{cases}
&\rho _t +(\rho u)_x +(\rho v)_y ~= ~0,\\
&(\rho u)_t + (\rho u^2 + p)_x + (\rho uv)_y ~=~0,\\
&(\rho v)_t + (\rho uv)_x + (\rho v^2 + p)_y ~=~0.
\end{cases}
\end{align}
Here $\rho$ is the density, $(u, v)$ is the fluid velocity, and $p(\rho)=A \rho ^{\gamma}$ is the pressure,  with $A>0, 
\gamma>1$. In the following we shall take the adiabatic   
constant $\gamma =1.4$ and $A=1$, unless otherwise stated.

The goal of the present work is to explore the multiplicity of solutions to compressible equations (\ref{eq:2DEuler}) by a numerical approach. 
Motivated by the  numerical construction for the incompressible Euler equation  (\ref{eq:2DinEuler}) in \cite{BS20}, we start with a similar form of singular vorticity profile, and construct three families of initial data  (given in polar coordinates) 
which approximate the same vorticity in the limiting case. Using these initial data, together with an initial density $\rho_0= r^\beta$  which is a power of the radial coordinate, we solve the system (\ref{eq:2DEuler}) and check whether they lead to distinct solutions, at any positive time. By carefully tuning the three parameters in initial data (exponent in the density, exponent in the vorticity, and angular support of the vorticity),  in our extensive numerical experiments we find several cases where
non-uniqueness of solutions can be observed through the vorticity profile.  
More specifically, with the other two of the three parameters fixed at certain values, we find that the continuous dependency of initial data is violated when 1) the exponent term in initial density is small enough, or 2) the exponent term in initial vorticity is within a certain range, or 3) the initial vorticity function is supported on a small enough angle. Hopefully,  these numerical results can be further validated by rigorous a posteriori error estimates, leading eventually to a computer-assisted proof of this striking phenomena.  

\subsection{Further related work} 
While the uniqueness of weak solutions to Euler equations with general velocity/vorticity profiles remains an open question, there is  a  body of literature that has been devoted to numerical study of possible non-uniqueness of the incompressible Euler and compressible Euler equations.  Different approaches have been suggested in \cite{Pu89, FLLZ06, Le18, LM20} for incompressible Euler equations and \cite{FKMT17, FLMW20, FLSW20}  for the compressible Euler equations.   

For the 2D incompressible Euler equation, promising candidates for scenarios of non-uniqueness are flows involving vortex sheets. The first non-unique vortex sheet evolution comes from the intriguing result of Pullin \cite{Pu89}.  Pullin considered multiple self-similar vortex sheet solutions from a single, initially flat single-signed vortex sheet with a specially chosen $x$-coordinate. His simulation based on self-similar configurations is suggestive of non-uniqueness  for the initial value problem. How to obtain such an example directly was left open. For the construction of some solutions with vorticity forming algebraic spirals near the origin, we refer to \cite{El13}.  In \cite{FLLZ06}, numerical evidence of non-uniqueness in the evolution of vortex sheets is given for the 2D incompressible equations with initial data containing smooth, single-signed vortex sheets. The same solution behavior was also observed in \cite{LM20} but with a different numerical method -- a spectral viscosity method to approximate the two-dimensional Euler equations with rough initial data is proposed and shown to converge to a weak solution for a large class of initial data, including when the initial vorticity is in the so-called Delort class i.e., it is a sum of a signed measure and an integrable function.

For the 2D compressible Euler equations, Elling \cite{El06} presented a numerical counterexample to the well posedness of entropy solutions in the presence of shock waves in similarity coordinates.  In \cite{FKMT17, FLMW20, FLSW20} the emphasis is on how to achieve numerical convergence of all interesting solutions. Numerical  experiments in \cite{FKMT17}  strongly  suggest  that  there  is  no  convergence  of  approximations  generated by  standard  numerical  schemes  as  the  mesh  is  refined, hence they considered the notion of entropy measure-value solutions introduced by DiPerna (1985). In \cite{FLMW20} statistical solutions are further considered by adding multi-point spatial correlations upon the measure-valued solutions. From their numerical experiments, they concluded that one observes convergence of all interesting statistical observables in that framework. In a similar spirit, the authors in \cite{FLSW20} proposed a method to compute the Young measures associated to sequences of numerical solutions based on the concept of K-convergence. 

However, all these works are based on generating solutions from rather singular initial data.  In contrast, our  initial data is much more regular than those in these papers.

 \subsection{Remarks on numerical results} 
 There are no rigorous convergence results to entropy solutions for any numerical schemes  approximating  multi-dimensional  systems  of  conservation  laws. We regard the numerical results obtained as an indication that $W^{1, p}$ initial velocity can lead to non-uniqueness, but a more extensive numerical study is definitely warranted.  The main purpose of this article is to suggest examples and provide convincing experiments. One may ask whether our computational solution may be converging to the entropy solution. This is less clear theoretically.  Due to the use of invariant-region-preserving methods, which is consistent with the entropy or energy in the present setting, we believe that the numerical approximation does yield an entropy solution of the 2D Euler equations. Finally, regarding the numerical evidence presented herein, it is possible that truncation error may be playing a role in our observations. Nontheless, our results showing non-uniqueness on refined grids are consistent with those on coarse grids. We are aware that the experiments performed might not be fully resolved, but this is not a problem since we do not deal with well-posedness, and numerically non-uniqueness is stronger than instability. 
 
The remainder of the paper is organized as follows. In the next section, we review the non-uniqueness results for the incompressible case presented in \cite{BS20}. In Section 3, we discuss the non-uniqueness for compressible equations, where the initial data designed for numerical tests are introduced and numerical results as well as implementation details are presented. In particular, we describe how different parameters in the initial data affect the uniqueness of solutions.  
In Section 4, we provide further numerical results to improve our understanding of how the compressibility can 
affect the structure of solutions. Finally, concluding remarks are given in Section 5. 

\section{Review of non-uniqueness results for the incompressible flow}
For the incompressible Euler flow (\ref{eq:2DinEuler}), solutions with
spiraling singularities were numerically constructed in \cite{BS20} .
The initial vorticity (given in polar coordinates) takes the form 
\begin{equation}\label{rb}
\omega_0(r,\theta)~=~ r^{-\alpha} \phi(\theta), \qquad\qquad (x, y)=(r \cos \theta, r \sin \theta).
\end{equation}
Here  $0<\alpha<2$  while   $\phi\in {\mathcal C}^\infty({\mathbb R})$ is a non-negative, smooth, periodic function
which satisfies
\[
\phi(\theta) ~=~\phi(\pi+\theta),\qquad\qquad
\phi(\theta)~=~0\quad\hbox{if}~~\theta\in \left[\frac{\pi}{ 4}\,,\, \pi\right].
\]
Notice that the initial vorticity $\omega_0$ is supported on two wedges,
symmetric w.r.t.~the origin,
and becomes arbitrarily large  as $|x|\to 0$.

The function $\omega_0$ can now be
approximated  by two families of  bounded initial data, by taking 
\begin{equation}\label{cex}
\omega'_{0,\epsilon} (r,\theta) =\left\{
\begin{array}{cl}
\omega _0(r,\theta)\quad &\text{if }|r|>\epsilon,\\[3mm]
\epsilon ^{-\alpha} \quad &\text{if }|r|\leq \epsilon,
\end{array}\right.\qquad
\quad \quad\omega''_{0,\epsilon} (r,\theta) =\left\{
\begin{array}{cl}
\omega _0(r,\theta)\quad &\text{if }|r|>\epsilon,\\[3mm]
0 \quad &\text{if }|r|\leq \epsilon.
\end{array}\right.\end{equation}
As $\epsilon \to 0$, both families converge to  $\omega_0$ in ${\bf L}^p_{\rm loc}
({\mathbb R})$,  for a suitable $p$ depending on the choice of $\alpha$. 

Since $\omega'_{0,\epsilon} ,\omega''_{0,\epsilon} \in {\bf L}^\infty(\mathbb{R})$
for every $\epsilon>0$, 
by Yudovich's theorem \cite{Yu63}  each of these initial data yields a unique solution. 

However, the numerical simulations in \cite{BS20} indicate that, as $\epsilon \to 0$, two distinct limit solutions are obtained. 
In the first solution, both wedges wind up together into a single spiral. 
On the other hand, in the second solution, each wedge curls up on itself and two distinct spirals are observed. This indicates that the ill-posedness of the 
two-dimensional incompressible Euler equation (\ref{eq:2DinEuler}) 
in $W^{1,p}_{loc}$ is ``incurable", since there is no way to choose a unique solution continuously depending on the initial data.

Some partial steps toward a rigorous validation of these numerical results were taken in \cite{BM20, BS20}. More precisely, in \cite{BS20} some a posteriori
error estimates were proved, for numerical approximations on a domain where the solution remains smooth.
In addition, in \cite{BM20} the authors constructed two types of analytical solutions: in a neighborhood of infinity, and in a neighborhood of the spiral's center where
the vorticity is unbounded.

\section{Non-uniqueness for compressible equations}
In this section we present the results of several numerical simulations, with carefully designed initial data,  checking whether the spiraling solutions found in the incompressible case are still produced.  
The underlying motivation is that, even for compressible flow, the vorticity
is passively transported along particle trajectories.  Therefore, if the vorticity is initially supported on two wedges, then at any time $t>0$,
we expect that the vorticity will still be supported on a set which is topologically 
equivalent to  two wedges. 

To be more specific, we work on a square domain $\Omega =[-a,a]\times [-a,a]$ with periodic boundary conditions. 
Using polar coordinates as in (\ref{rb}), 
we consider an initial density of the form
\begin{equation}\label{rho0}
\rho_0(r,\theta)\,=\,r^\beta, 
\end{equation}
for some $\beta \geq 0$.
The components of the initial velocity are obtained from 
\begin{align}
\label{velocity}
u = \psi _y, \quad v = -\psi _x,
\end{align}
where the stream function $\psi$ is recovered from the vorticity
using the Poisson equation
\begin{align}
\label{vorticity}
\psi _{xx}+\psi _{yy}=\omega, 
\end{align}
subject to an appropriate boundary condition.   Since what matters is the behavior of the solution near the origin, which is not much affected by the boundary conditions, we simply adopt the zero boundary conditions for $\psi$.

As initial vorticity, we take the same kind of profile as in (\ref{rb}), namely
\begin{equation}\label{barom}
\bar \omega (r,\theta) \,= \,r^{-\alpha} \phi(\theta).
\end{equation}
Here $0<\alpha<2$, while $\phi$ is $\pi$- periodic: $\phi(\theta+\pi)=\phi(\theta)$,
and 
\begin{align*}
\phi(\theta)=
\begin{cases}
\theta _0 -|\theta|, \quad &\text{if } |\theta| <  \theta _0,\\
0, \quad &\text{if } \theta \in[-\pi/2, -\theta_0]\cup[\theta_0, \pi/2],
\end{cases}
\end{align*}
for some given angle $0<\theta_0<\pi/2$.

\begin{figure}[!htbp]
\centering
\includegraphics[width=.3\textwidth]{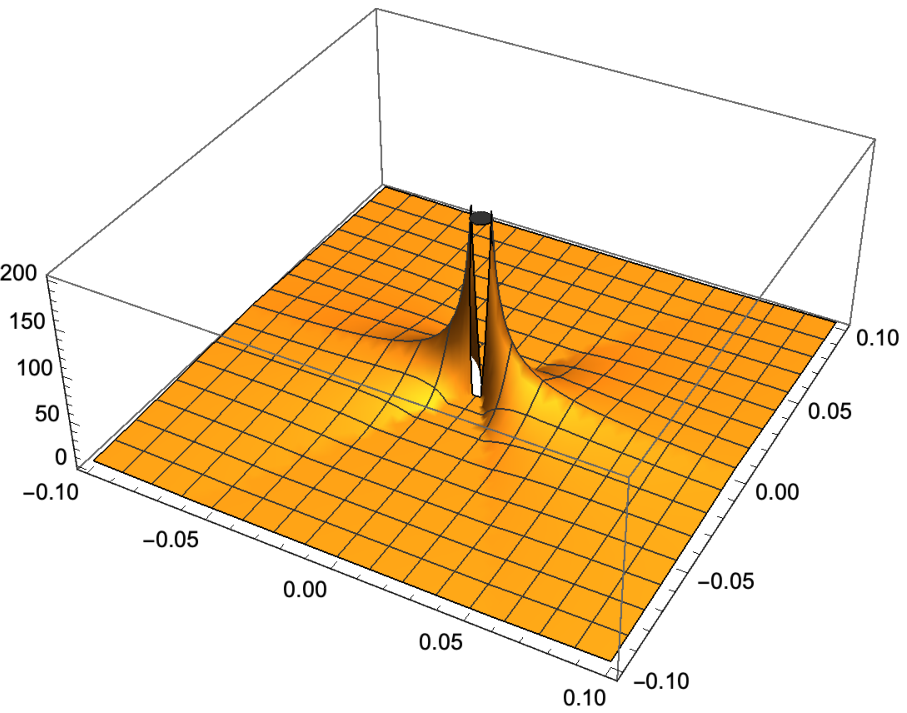}
\includegraphics[width=.3\textwidth]{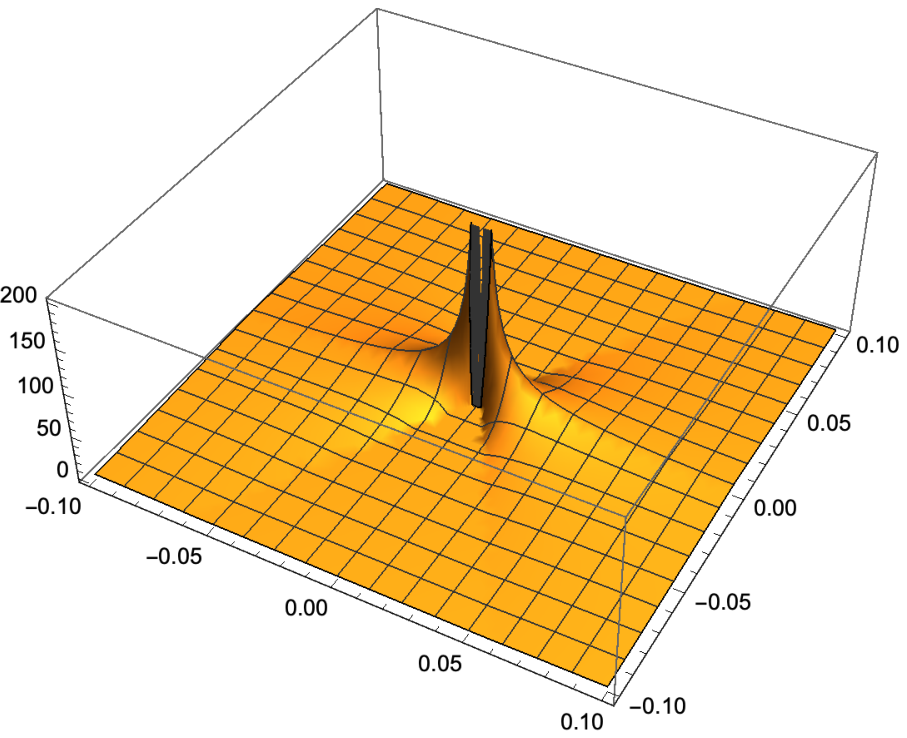}
\includegraphics[width=.3\textwidth]{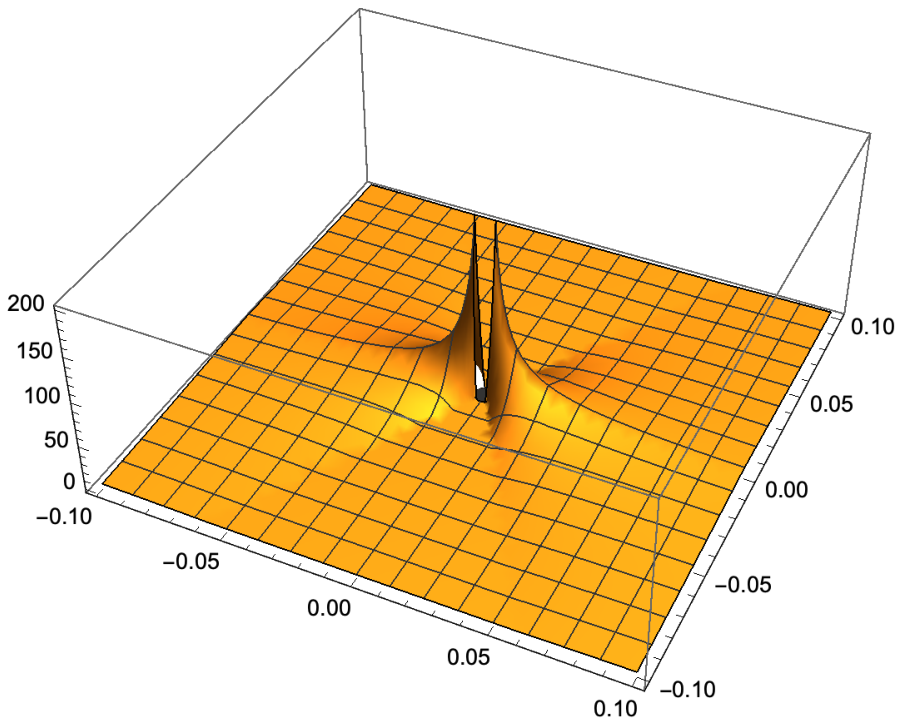}
\caption{\small
Three families of initial vorticity function in a small neighborhood of the origin: $(x,y)\in [-0.05,0.05]\times[-0.05,0.05]$, where $\alpha = 0.95$, $\theta _0 = \frac{\pi}{3}$, $\epsilon = 0.004$. From left to right: Case 0, Case 1 and Case 2.}
\label{InitialVorticity2}
\end{figure}

We approximate the initial vorticity $\bar \omega$ by three families of vorticity functions:
\begin{align*}
&\underline{\text{Case 0}}: \bar \omega ^{\epsilon,0} (r,\theta) =
\begin{cases}
\bar \omega(r,\theta), \quad &\text{if }|r|>\epsilon,\\
\epsilon ^{-\alpha}, \quad &\text{if }|r|\leq \epsilon,
\end{cases}
\\
&\underline{\text{Case 1}}: \bar \omega ^{\epsilon,1}(r,\theta)=
\begin{cases}
\bar \omega(r,\theta), \quad &\text{if }|r|>\epsilon,\\
\epsilon ^{-\alpha}\phi(\theta), \quad &\text{if }|r|\leq \epsilon,
\end{cases}
\\
&\underline{\text{Case 2}}: \bar \omega ^{\epsilon,2}(r,\theta)=
\begin{cases}
\bar \omega(r,\theta), \quad &\text{if }|r|>\epsilon,\\
0, \quad &\text{if }|r|\leq \epsilon.
\end{cases}
\end{align*}
As $\epsilon\to 0$, all these three functions converge  to $\bar \omega$ in 
${\bf L}^1_{\rm loc}$. See their plots in Figure \ref{InitialVorticity2}.

In our simulations, we will mainly focus on the comparison between the results from Case 0 and Case 2. We shall see in Section \ref{sec:alpha} that Case 1 provides an ``intermediate" solution, which, as time increases, converges either to the``one-spiral" or to the ``two-spirals solution". This suggests that, after one or two spiraling vortices have formed, these two solutions are both locally stable, and attract all nearby solutions.

\subsection{Implementation details}
Numerical implementation goes as follows: we consider a uniform discretization in space $\Delta x=\Delta y=\frac{2a}{N}$, where the number of computational cells is $N\times N$. The procedure
for numerically solving the system is the following.
\begin{itemize}
\item[(i)]Given a vorticity profile $\omega$, we solve the Poisson equation (\ref{vorticity}) using the five-point Laplacian scheme to obtain the approximation of the stream function $\psi$ at grid points;
\item[(ii)]then we use the second order central finite difference to get the approximation of two velocities $u$ and $v$ at grid points respectively;
\item[(iii)]in each square computational cell, we use three out of four values at corners (grid points) to reconstruct the $P^1$ polynomial approximation of the two velocity functions respectively;
\item[(iv)]using these two velocity approximations together with the density function $\rho_0$ as the initial data, we solve the Euler equations using $P^1$-DG method in space and the third order strong-stability-preserving Runge-Kutta method \cite{ShuOsher88} in time, where the invariant-region-preserving limiter introduced in \cite{JL18} is applied. Note that, for the isentropic system we consider here, only the positivity of the density needs to be preserved, and the limiter reduces to the usual positivity-preserving-limiter \cite{ZhangShuDG}. 
The numerical flux used in the DG method is the Lax-Friedrich flux, which is an invariant-region-preserving flux as shown in \cite{JL18}.
\end{itemize}

We follow the steps below to construct the vorticity profile at the final time $T$.
\begin{itemize}
\item[(i)] Use the $P^1$-DG solutions (that is, $P^1$ polynomial approximation of the density function, and two momentum functions) to construct the two velocity functions $u$ and $v$ at the center of each computational cell;
\item[(ii)]use the second order central finite difference to approximate $u_y$ and $v_x$, and then find the vorticity through (\ref{velocity}) and (\ref{vorticity}).
\end{itemize}

In the following experiments, we investigate what choice of parameters in the initial data (that is, $\beta$, $\alpha$, and $\theta_0$) could lead to non-uniqueness phenomenon. Without further specification, we take $a=0.2$, $\epsilon=0.004$, and  $N=200$. Note that by such choice, $\Delta x=0.002$ is smaller than $\epsilon$ so that the difference between different cases of the initial data is guaranteed to be captured. In particular, in the discussion of effects of $\beta$, when non-unique solutions are indicated, we have a mesh refinement study as well as an asymptotic study as $\epsilon$ goes to zero to further confirm our results. A related quantitative comparison between solutions as time evolves is also presented.

Since our initial data is assigned to vorticity, we focus on the behavior of vortivity solutions in the following experiments. The non-uniqueness phenomenon can also be observed in density solutions but slightly less intuitive. All the plots are made by using the surface plot function in MATLAB where the color is specified by a cut-off of the relative value of the vorticity vector (with respect to its maximum value). 
We look at the projection of the vorticity on the $x$-$y$ plane unless otherwise specified.

\subsection{Numerical results}
\subsubsection{Effects of $\beta$}\label{SubsecBeta}
We first test the initial data, Case 0 and Case 2, with different values of $\beta$, the power term in the initial density function. Based on our experience with the incompressible equations, we fix $\alpha =0.95$ and $\theta _0 = \frac{\pi}{8}$. \\

\begin{figure}[!htbp]
\centering
\includegraphics[width=.3\textwidth]{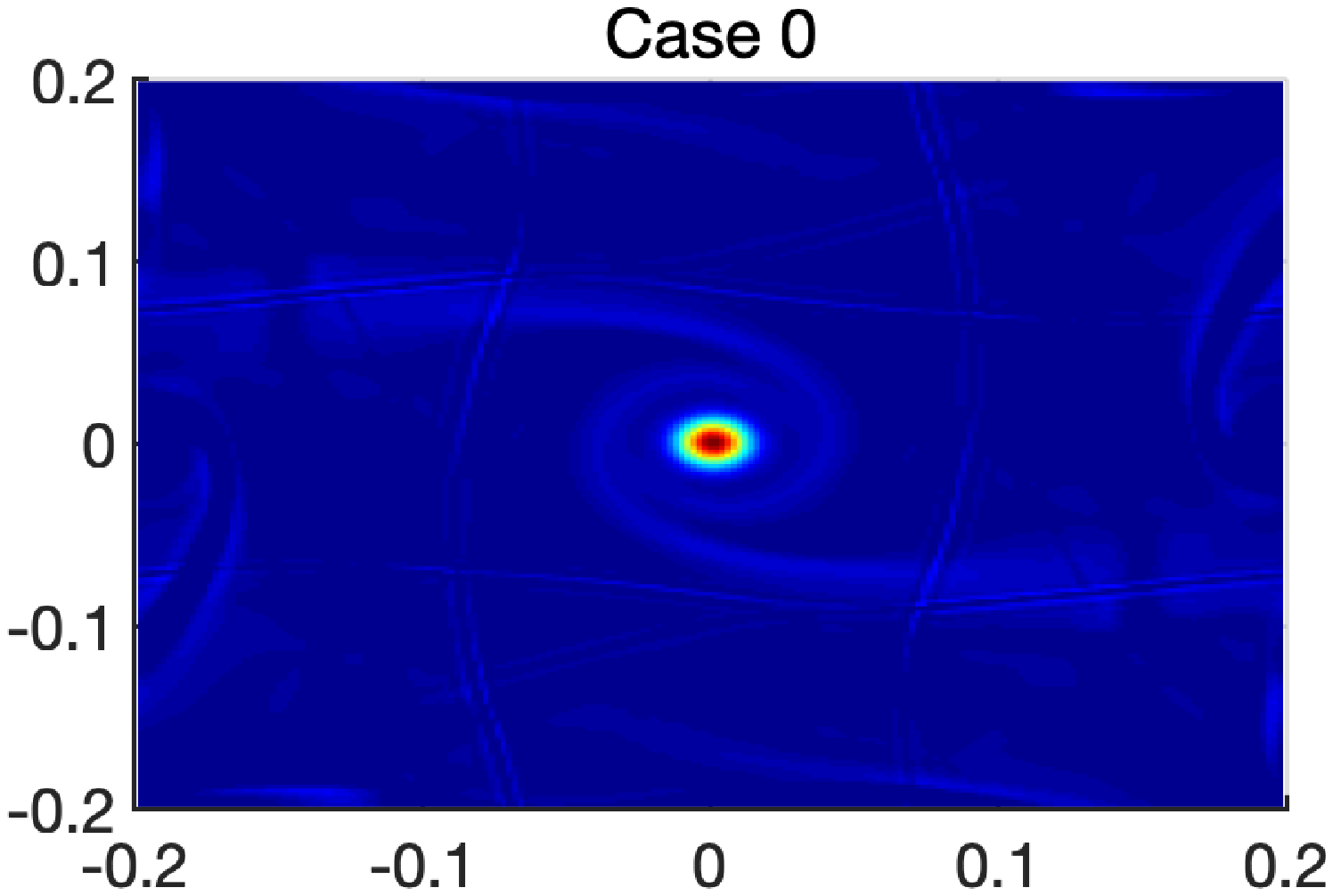}
\includegraphics[width=.3\textwidth]{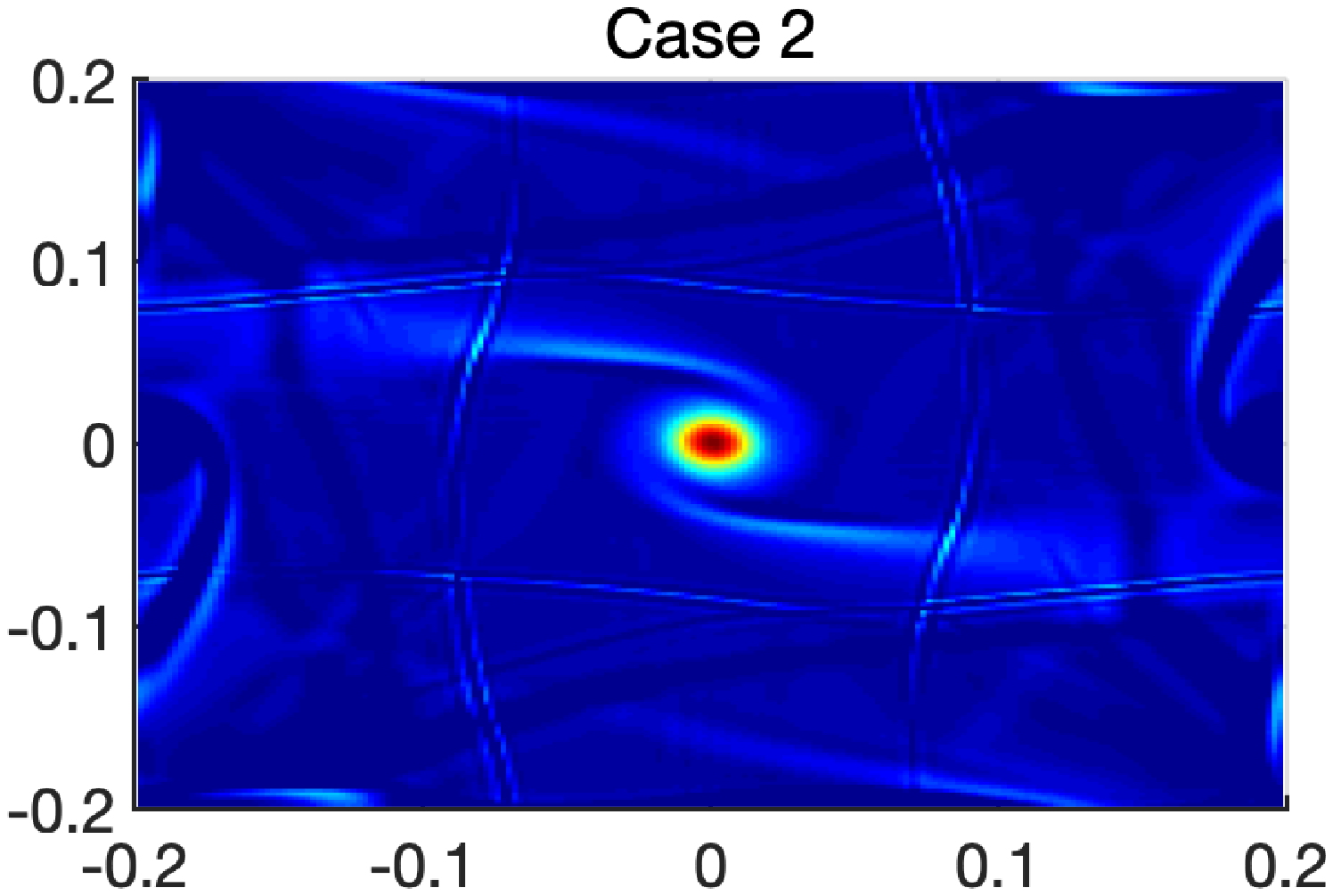}
\caption{Example 1: vorticity profiles at $T=1$. $\beta = 1$, $\alpha =0.95$ and $\theta _0 = \frac{\pi}{8}$.}
\label{EX1}
\end{figure}

\noindent \textit{Example 1.} $\beta = 1$. Figure \ref{EX1} shows that the two cases generate very similar vorticities, which are single spirals.\\

\begin{figure}[!htbp]
\centering
\includegraphics[width=.3\textwidth]{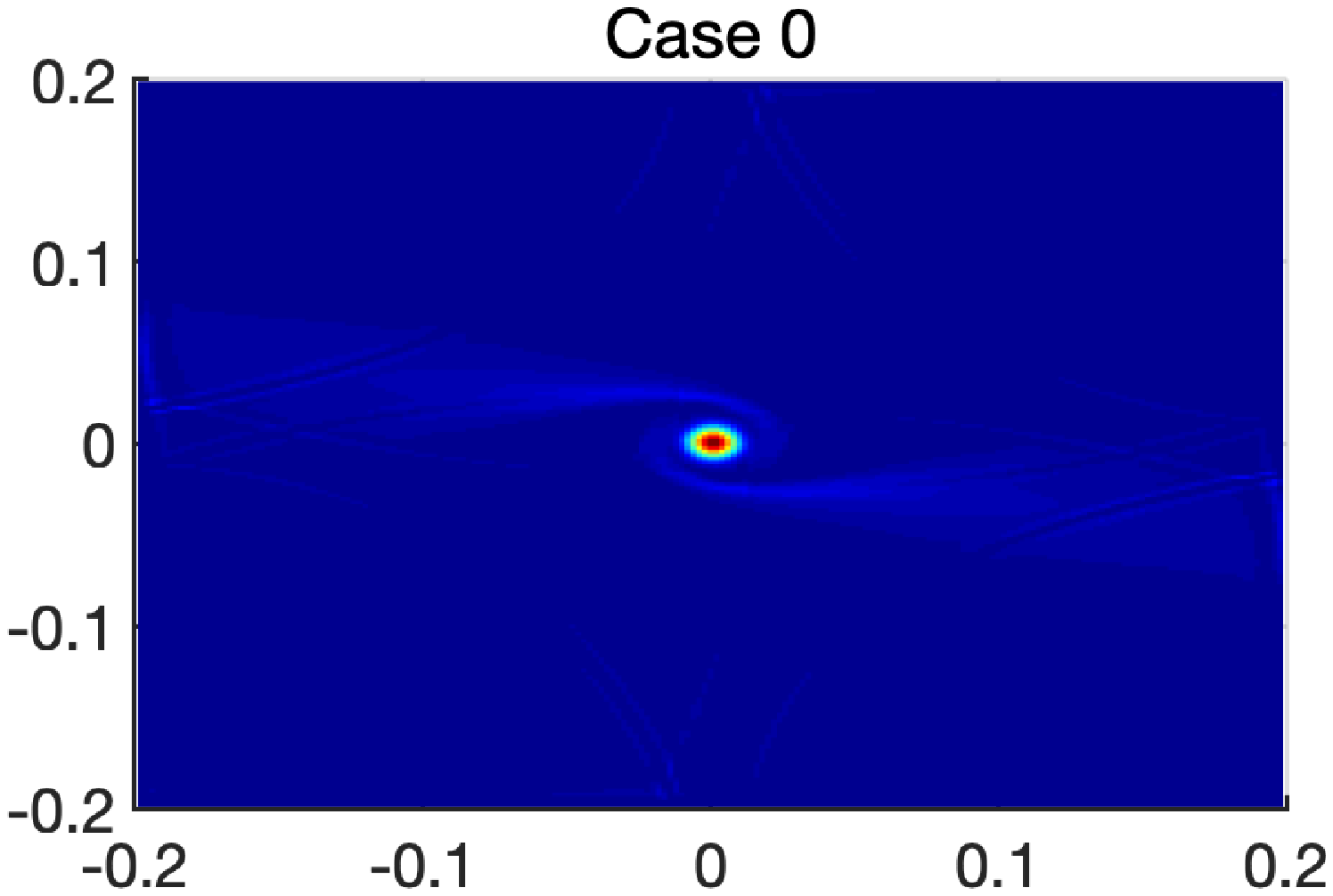}
\includegraphics[width=.3\textwidth]{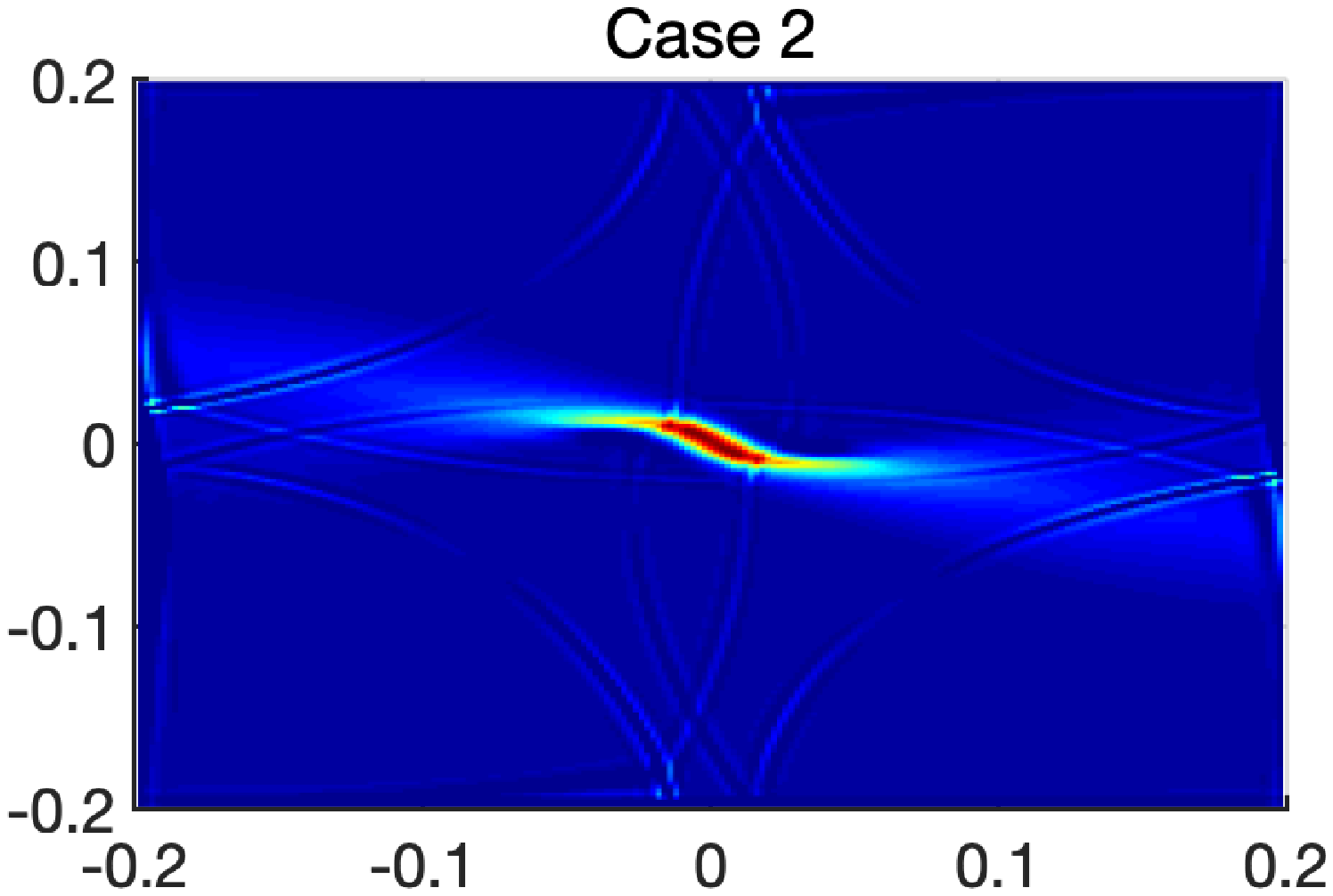}\\
\includegraphics[width=.3\textwidth]{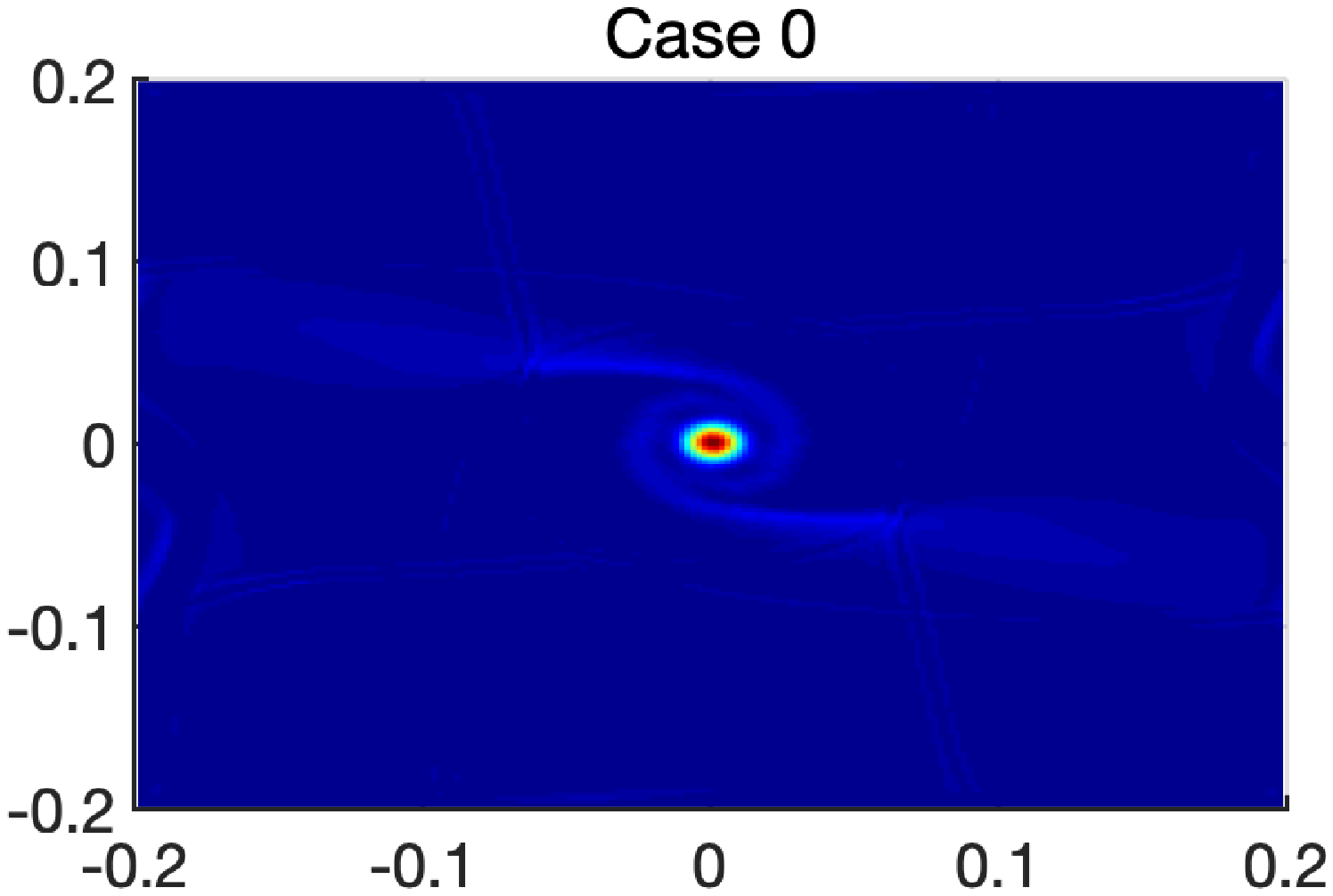}
\includegraphics[width=.3\textwidth]{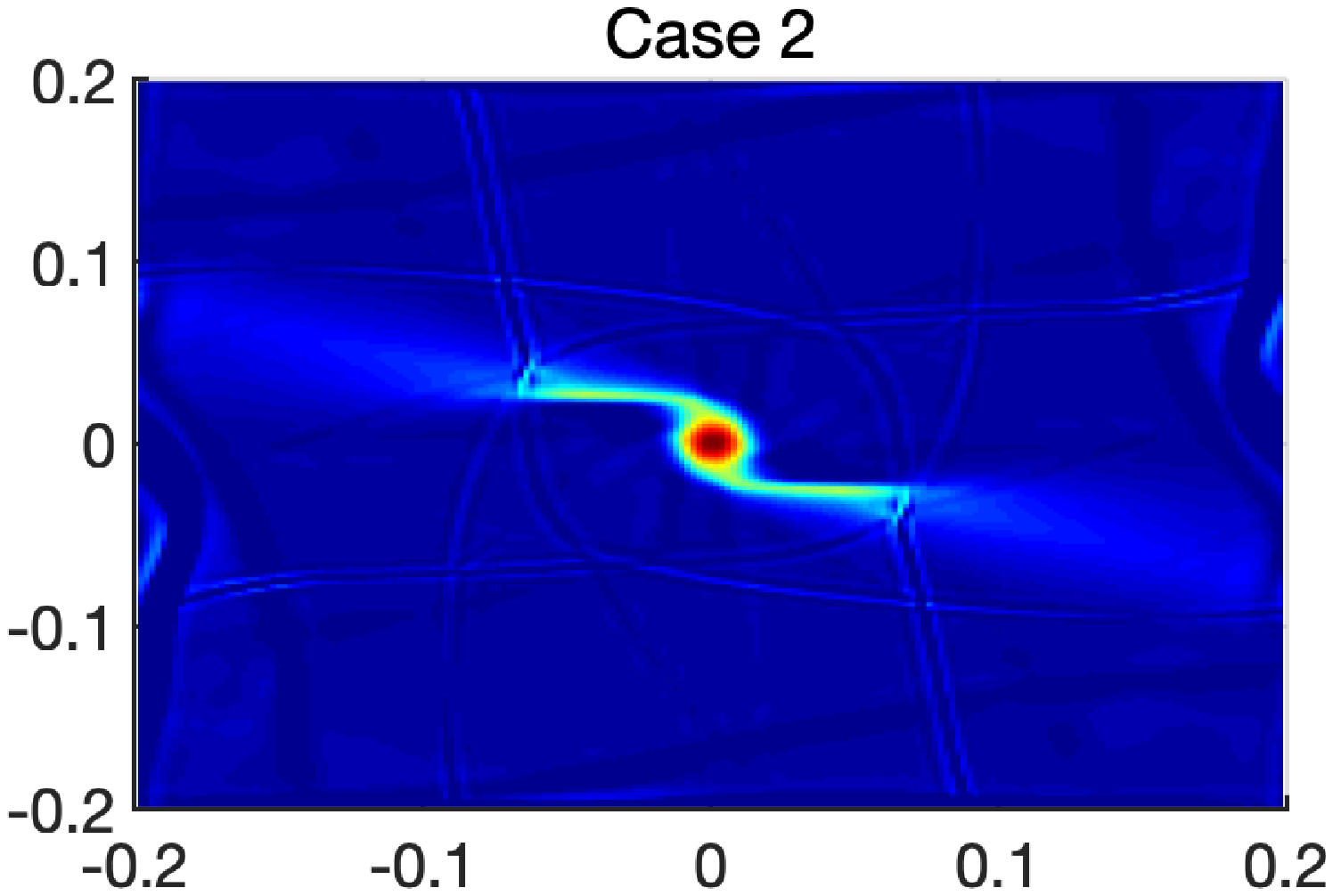}
\caption{Example 2: vorticity profiles at T=0.5 (top) and T=1 (bottom). $\beta = 0.5$, $\alpha =0.95$ and $\theta _0 = \frac{\pi}{8}$.}
\label{EX22}
\end{figure}

\noindent \textit{Example 2.} $\beta = 0.5$. We first run the simulation at $T=0.5$. We notice that the result from Case 2 has a different shape compared to that from Case 0. However, when we test it with $T=1$, the vorticity profile in Case 2 also becomes a single spiral. See Figure \ref{EX22}.\\

\begin{figure}[!htbp]
\centering
\includegraphics[width=.3\textwidth]{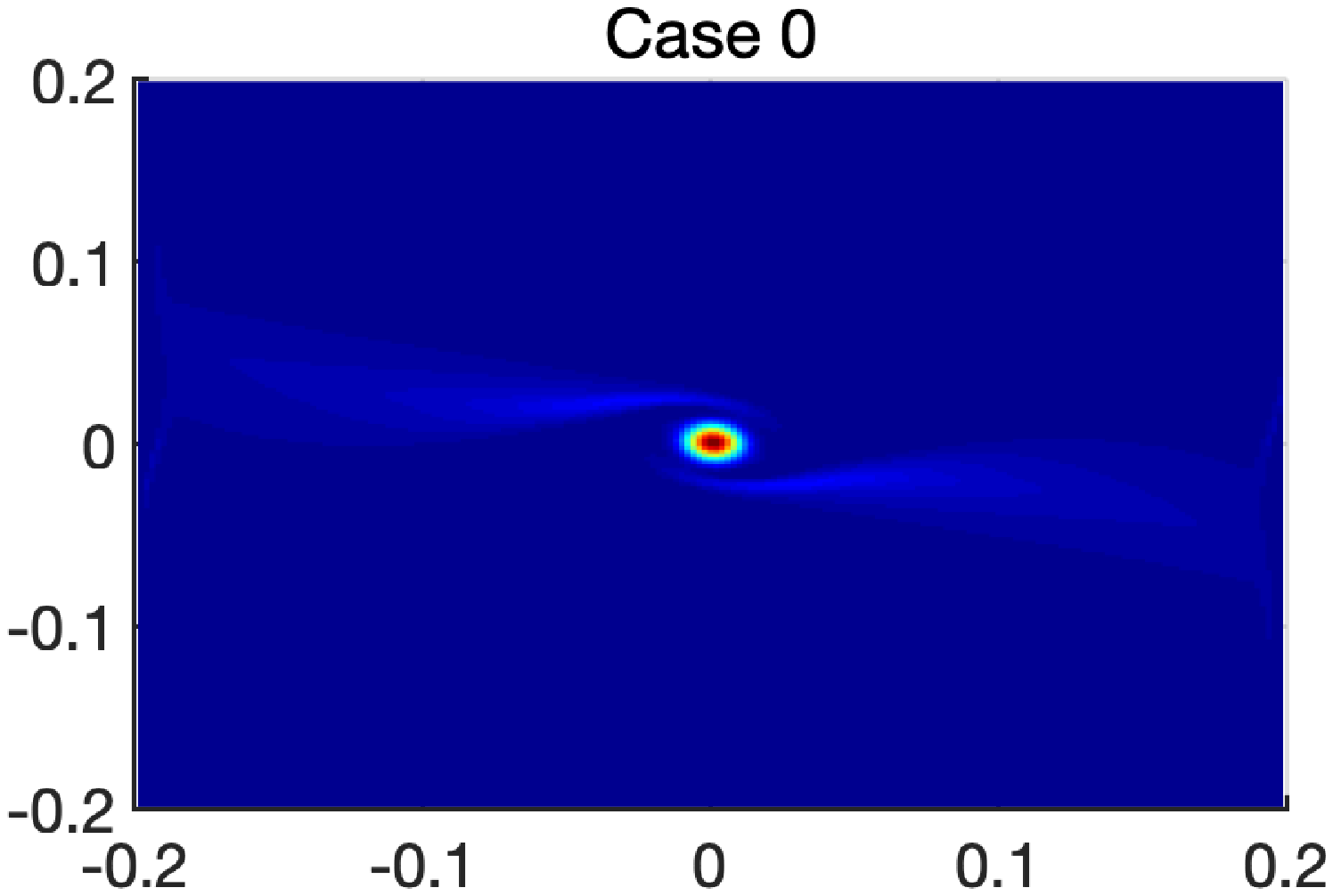}
\includegraphics[width=.3\textwidth]{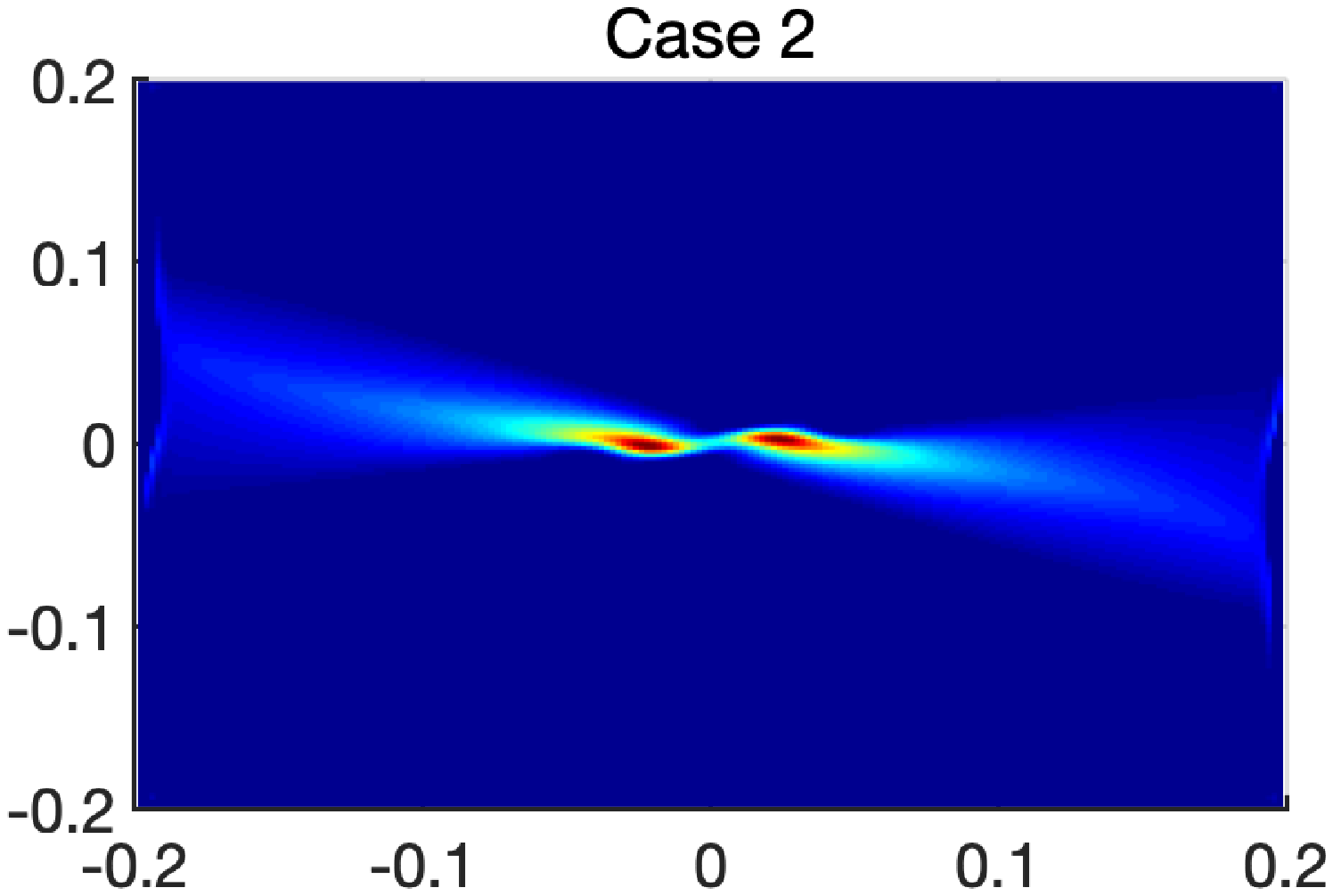}
\caption{Example 3: vorticity profiles at T=0.5. $\beta = 0$, $\alpha =0.95$ and $\theta _0 = \frac{\pi}{8}$.}
\label{EX31}
\end{figure}

\begin{figure}[!htbp]
\centering
\includegraphics[width=.3\textwidth]{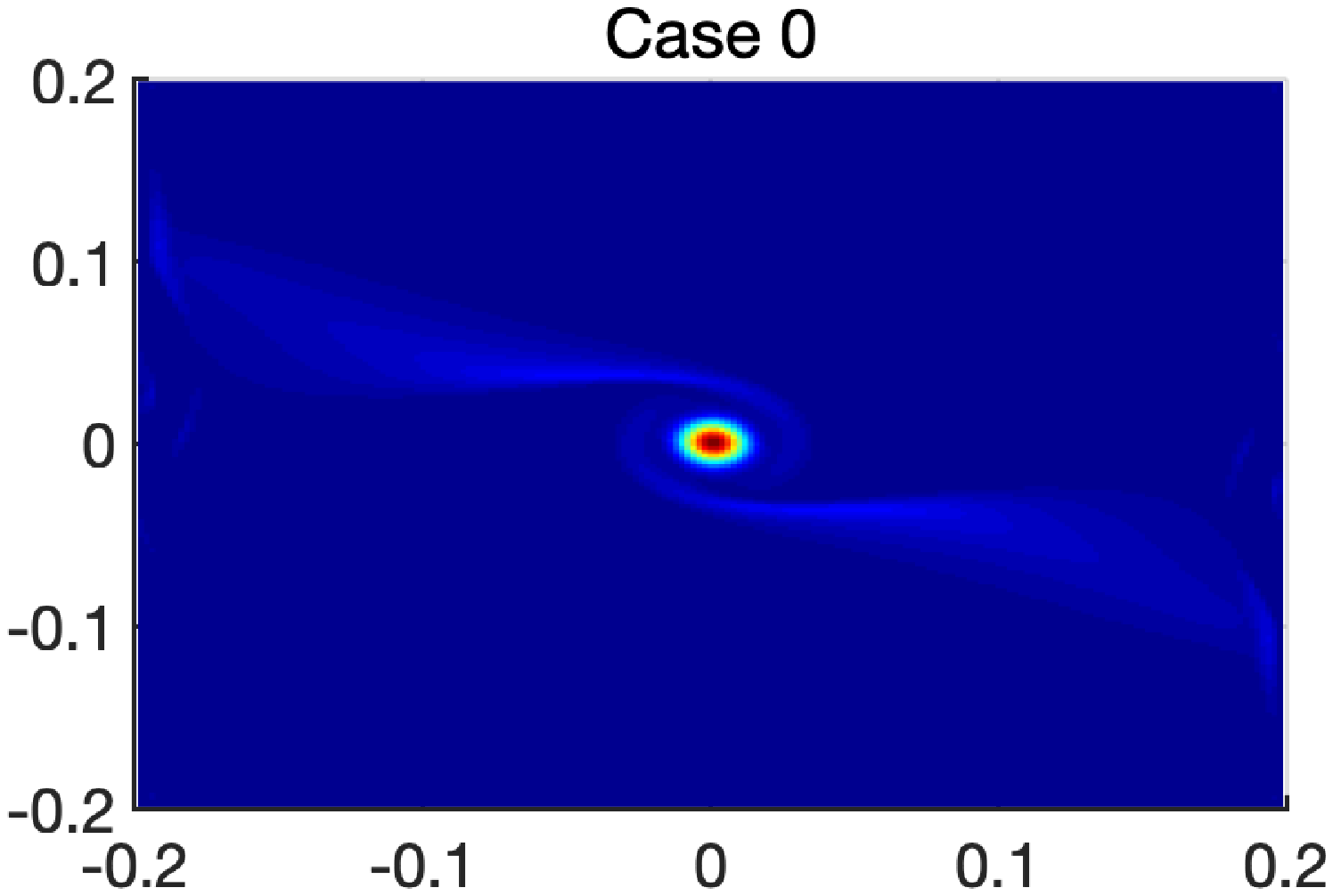}
\includegraphics[width=.3\textwidth]{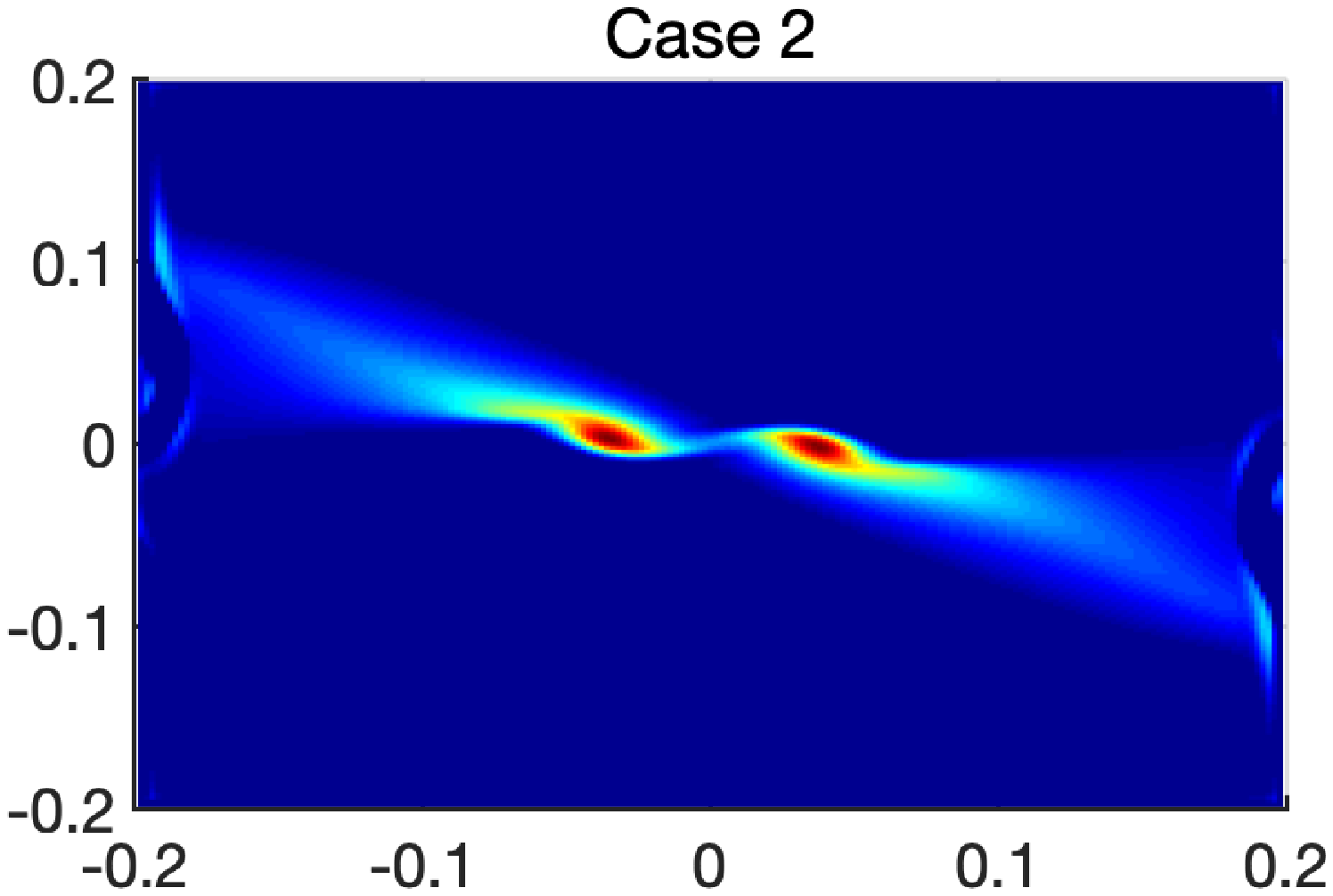}\\
\includegraphics[width=.3\textwidth]{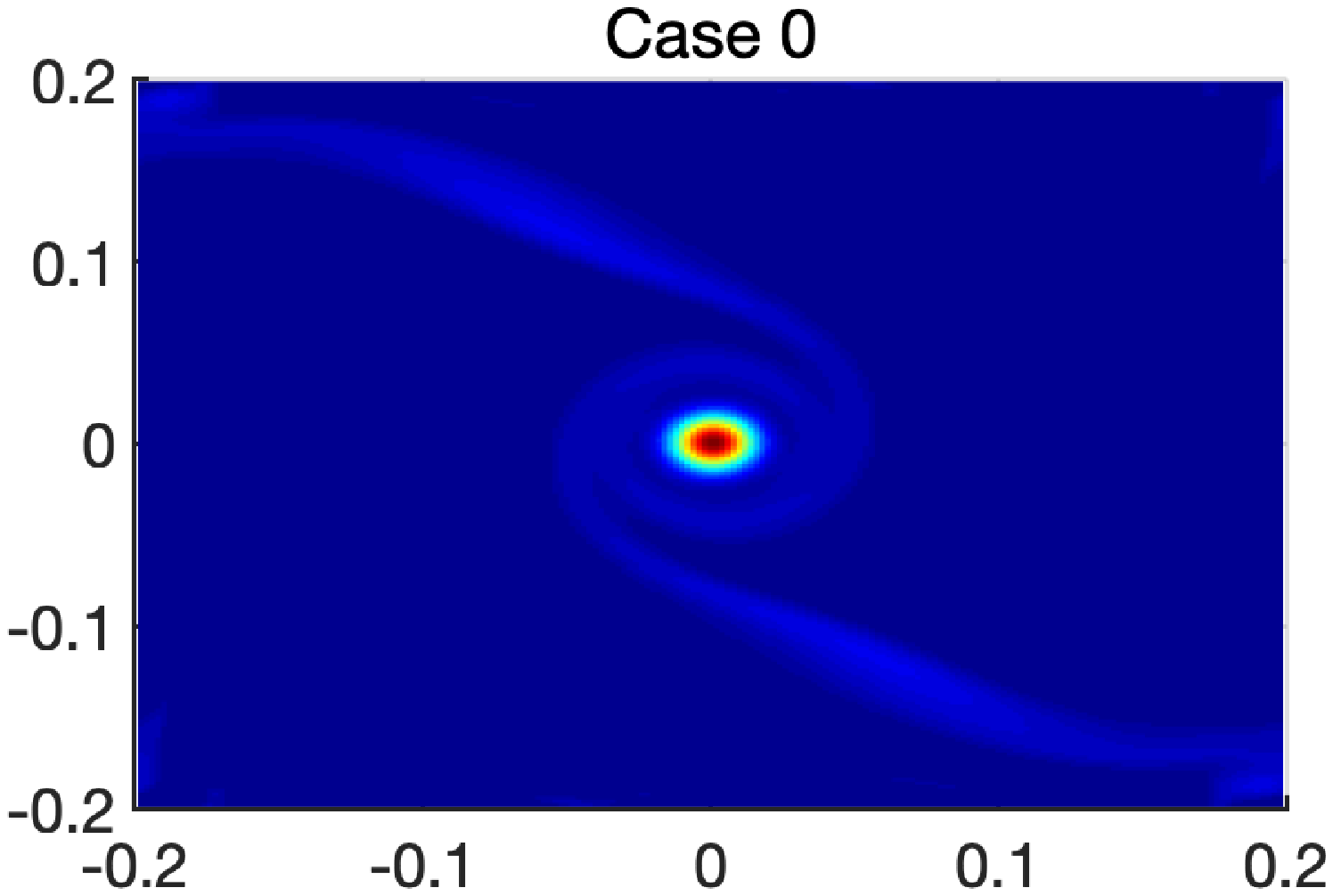}
\includegraphics[width=.3\textwidth]{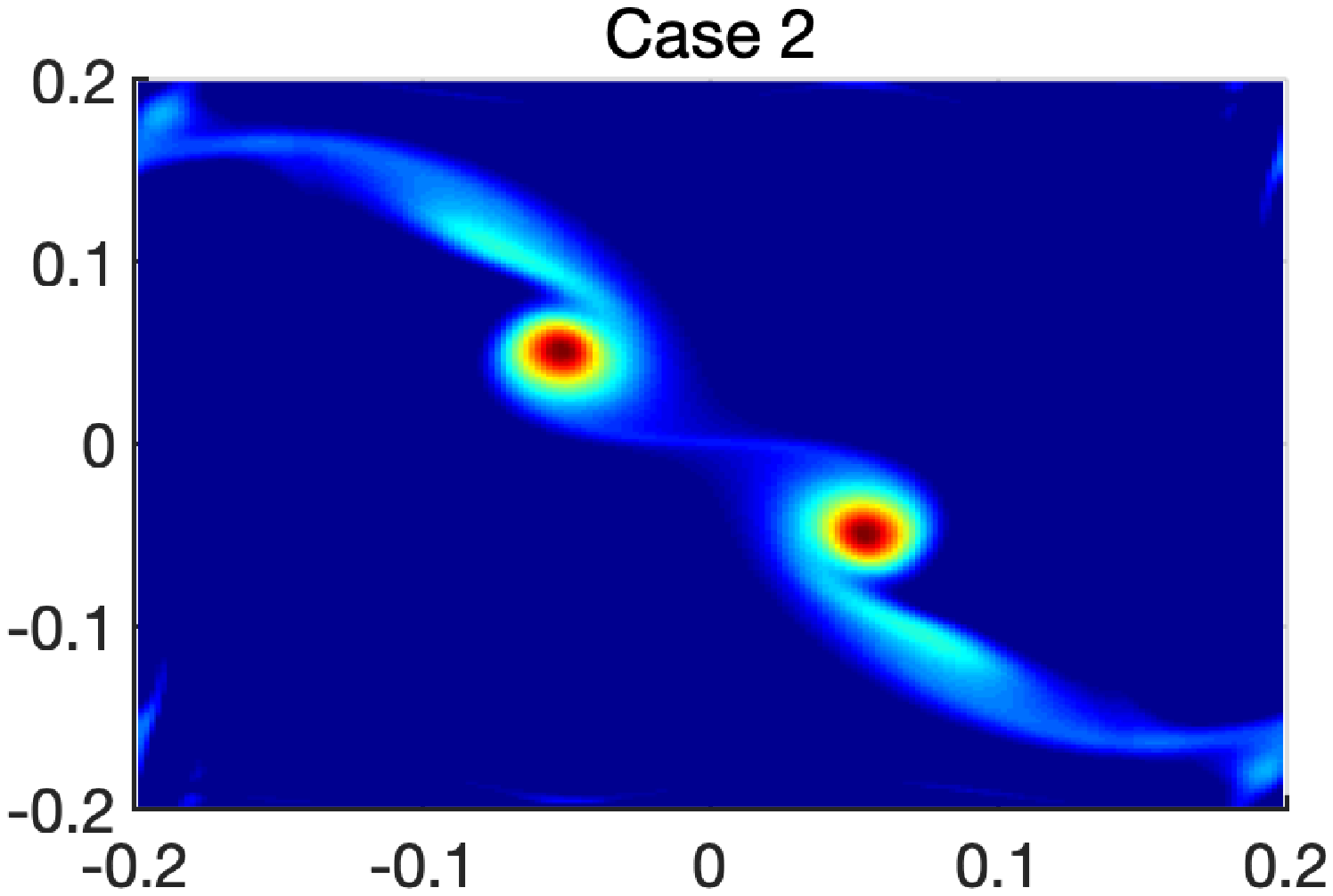}
\caption{Example 3: vorticity profiles at larger times. $\beta = 0$, $\alpha =0.95$ and $\theta _0 = \frac{\pi}{8}$. Top: $T=1$; bottom: $T=3$.}
\label{EX32}
\end{figure}

\begin{figure}[!htbp]
\centering
\includegraphics[width=.3\textwidth]{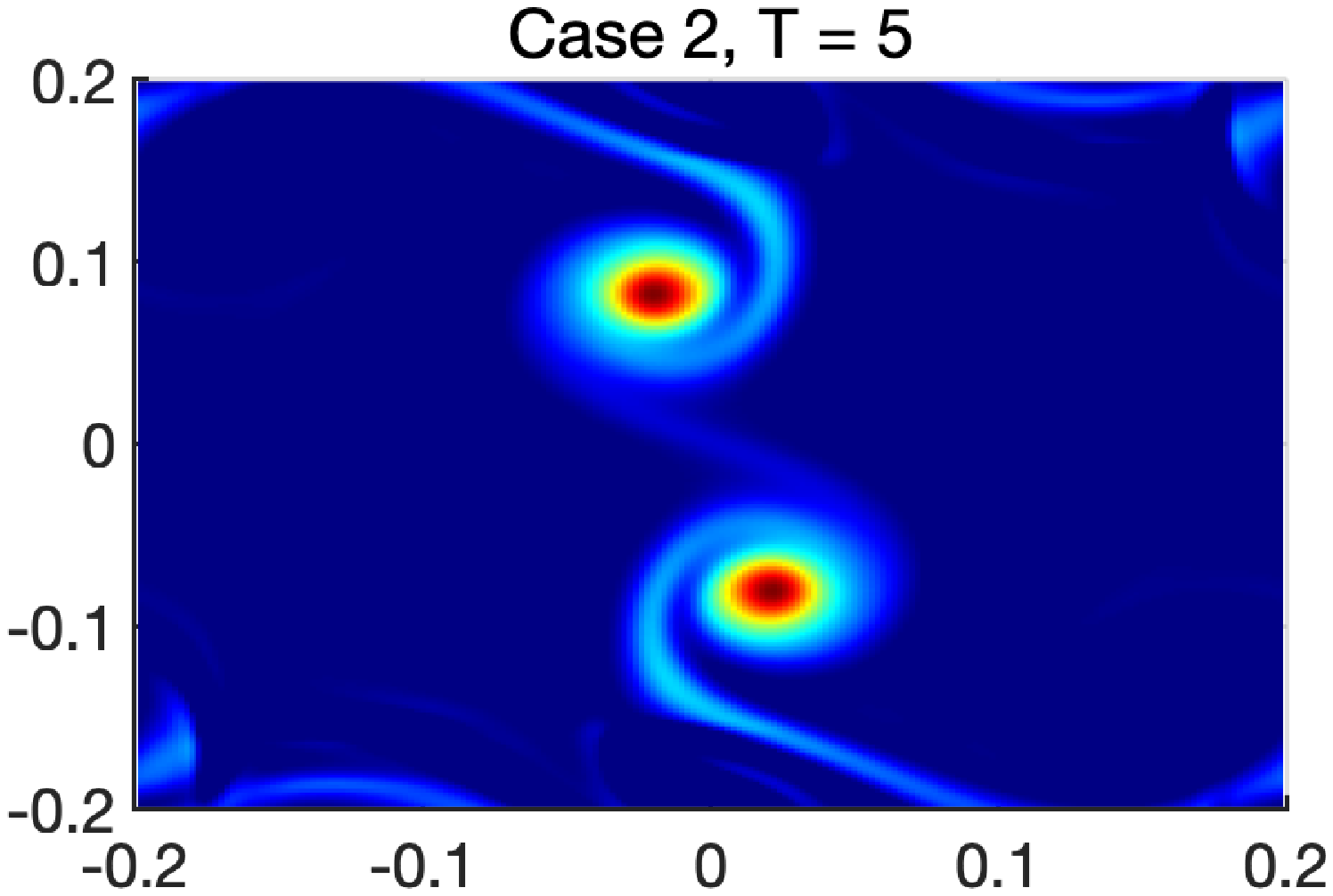}
\includegraphics[width=.3\textwidth]{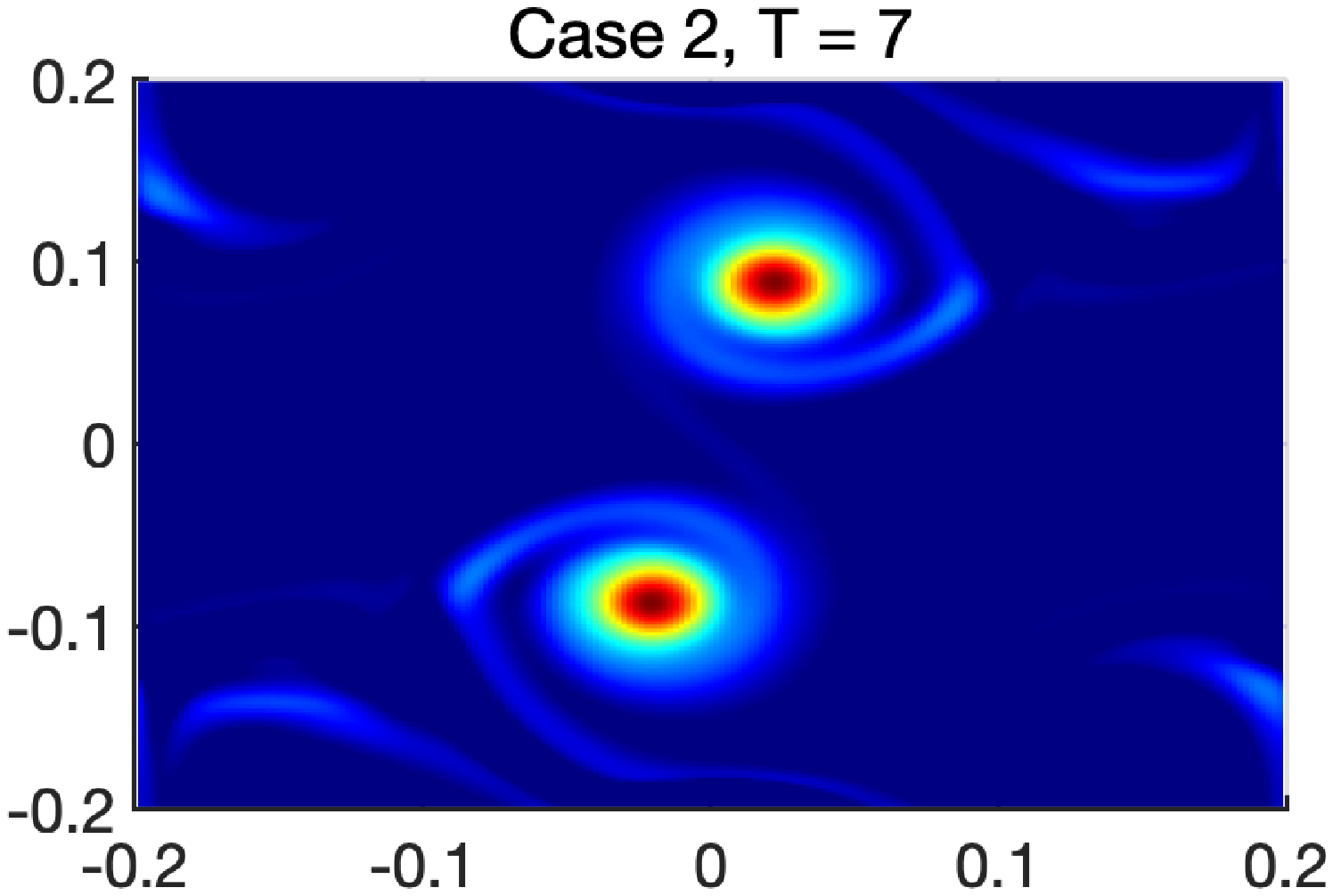}
\includegraphics[width=.3\textwidth]{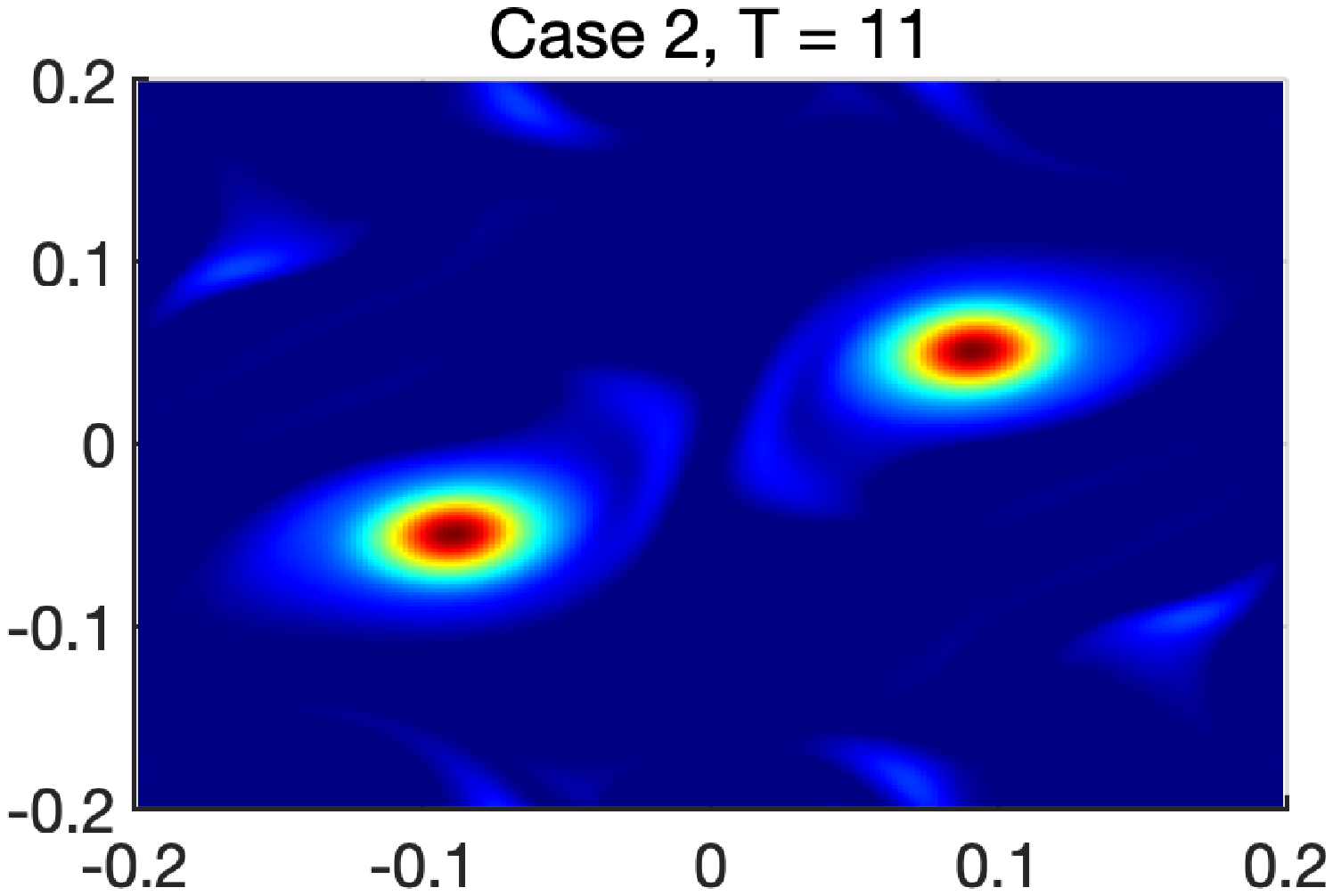}
\caption{Example 3: Case 2 vorticity profiles at larger times. $\beta = 0$, $\alpha =0.95$ and $\theta _0 = \frac{\pi}{8}$.}
\label{EX33}
\end{figure}

\noindent \textit{Example 3.} $\beta = 0$. This is the case where the initial density is constant: $\rho _0=1$. We first run the simulation at $T=0.5$. We notice that two cases generate two different shapes of vorticities, where the one from Case 0 is a single spiral while the one from Case 2 has two peaks. See Figure \ref{EX31}. We then test them with larger times $T=1$ and $T=3$. It shows that two spirals are generated in Case 2, which indicates the non-uniqueness of solutions. See Figure \ref{EX32}. {More Case 2 vorticity profiles at larger times ($T=3,5,7$) are shown in Figure \ref{EX33}}.\\

\begin{figure}[htbp]
\centering
\includegraphics[width=.3\textwidth]{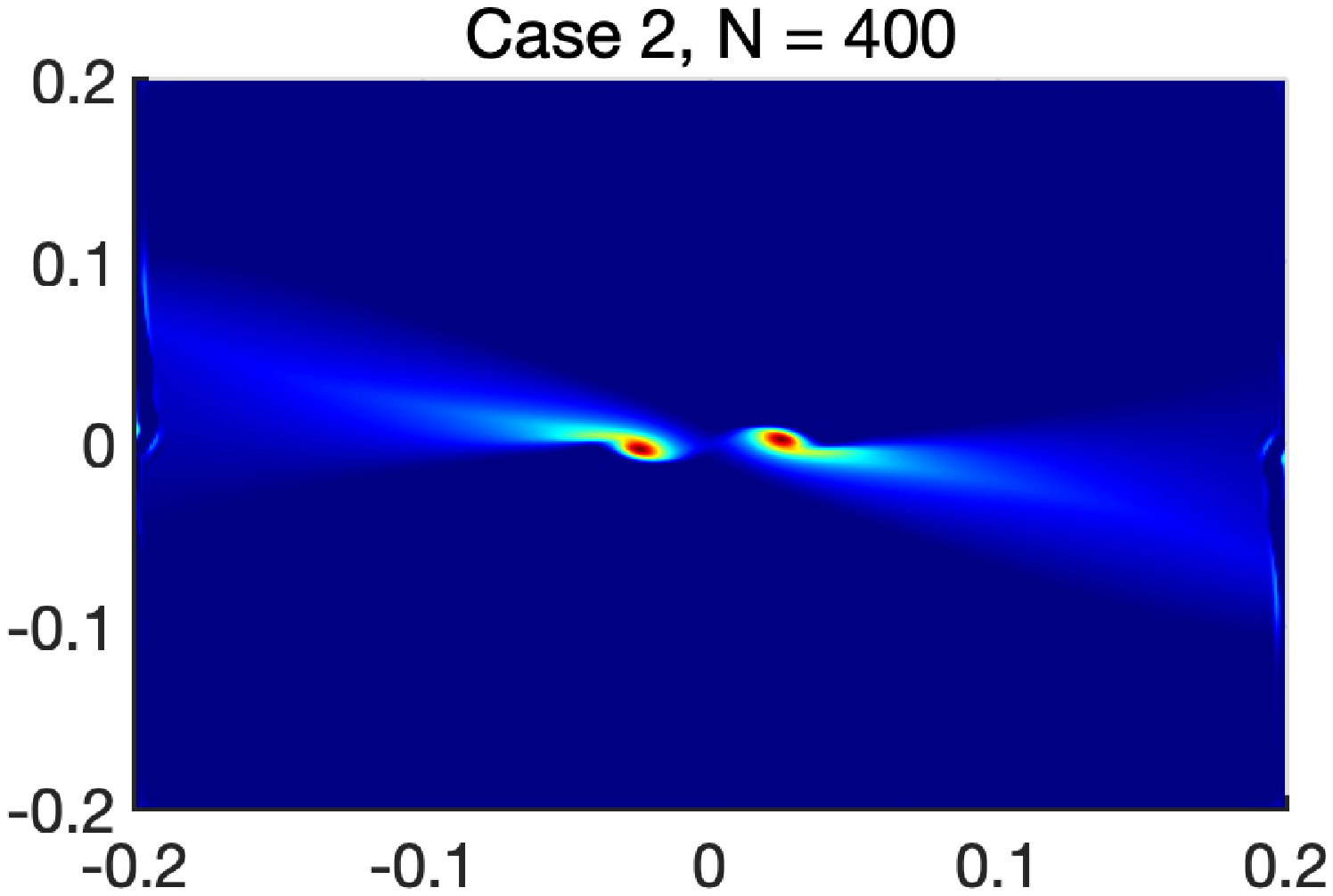}
\includegraphics[width=.3\textwidth]{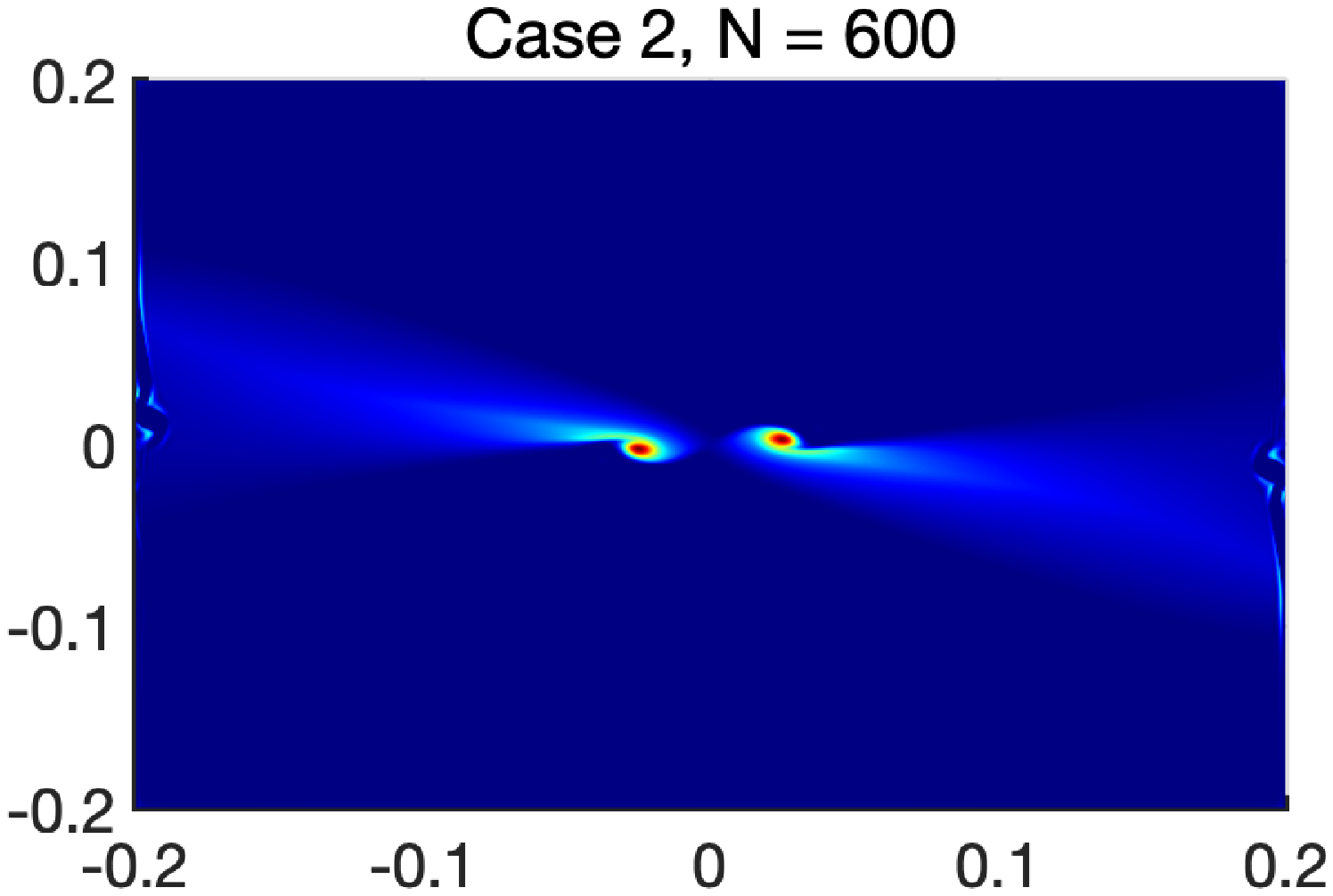}
\includegraphics[width=.3\textwidth]{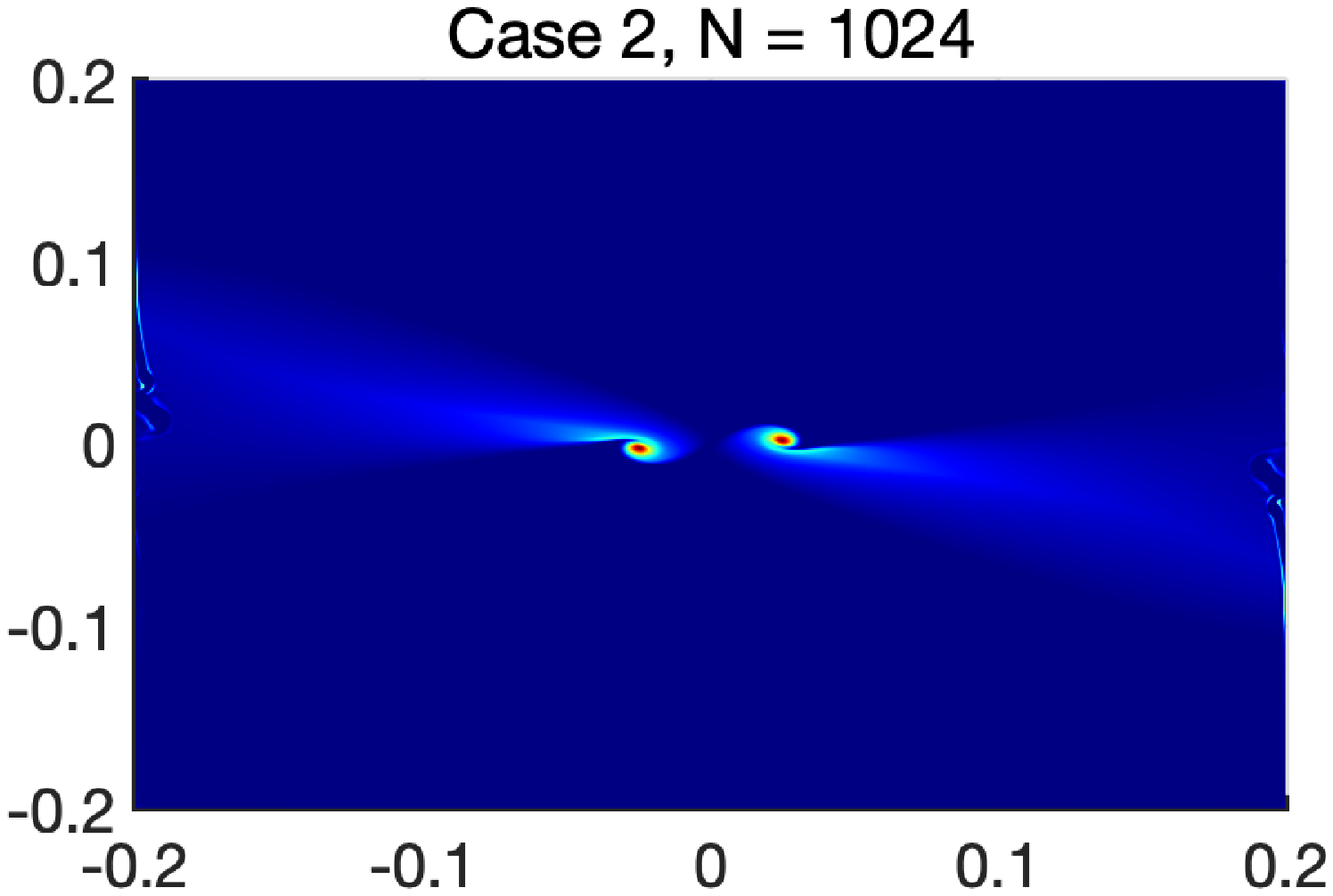}\\
\includegraphics[width=.3\textwidth]{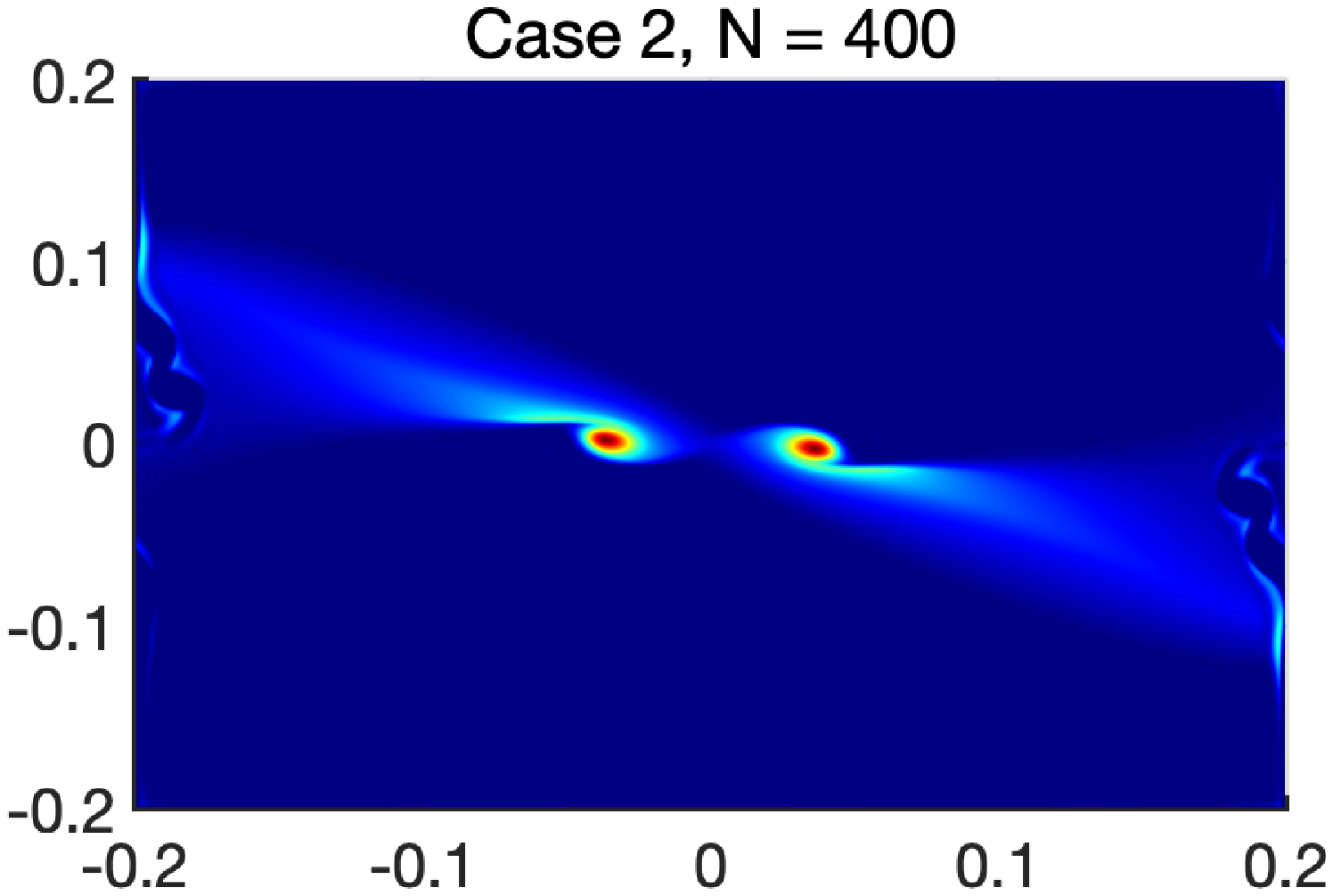}
\includegraphics[width=.3\textwidth]{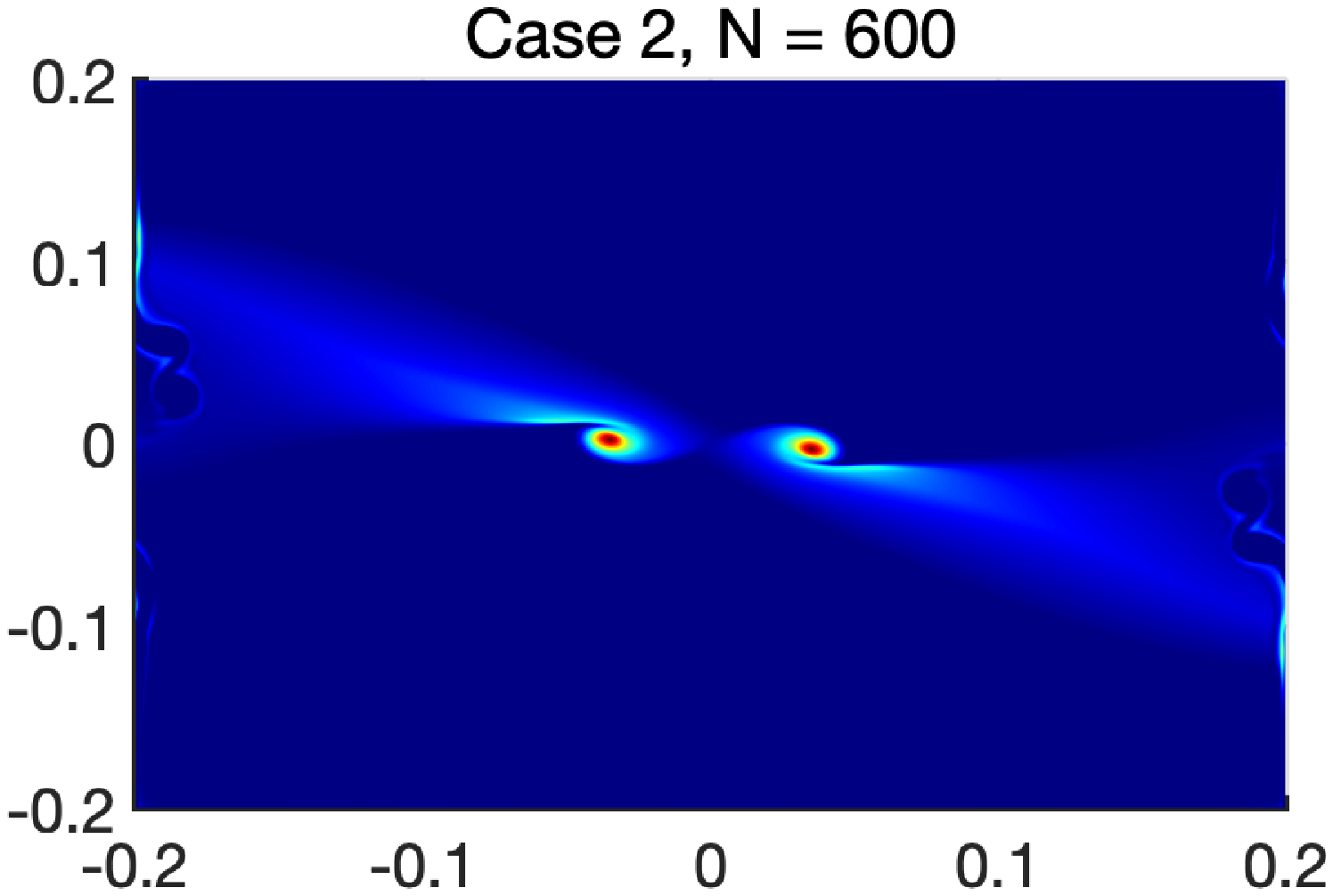}
\includegraphics[width=.3\textwidth]{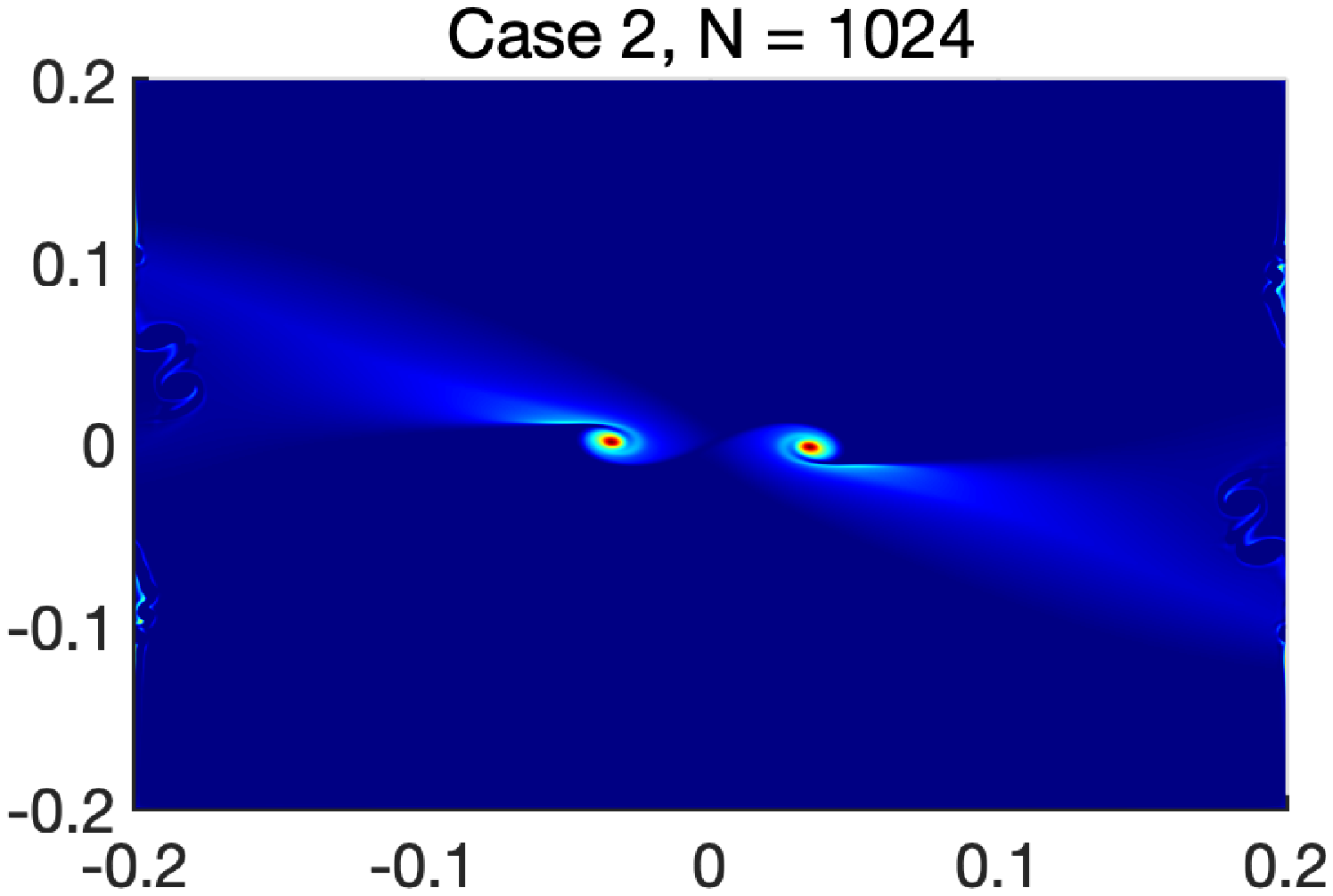}
\caption{Example 3: vorticities obtained from Case 2 initial data with different meshzie and at different times. $\beta =0$, $\alpha =0.95$, $\theta =\frac{\pi}{8}$. Top: $T=0.6$; bottom: $T=1$.}
\label{Vrefine}
\end{figure}

\begin{figure}[htbp]
\centering
\includegraphics[width=.45\textwidth]{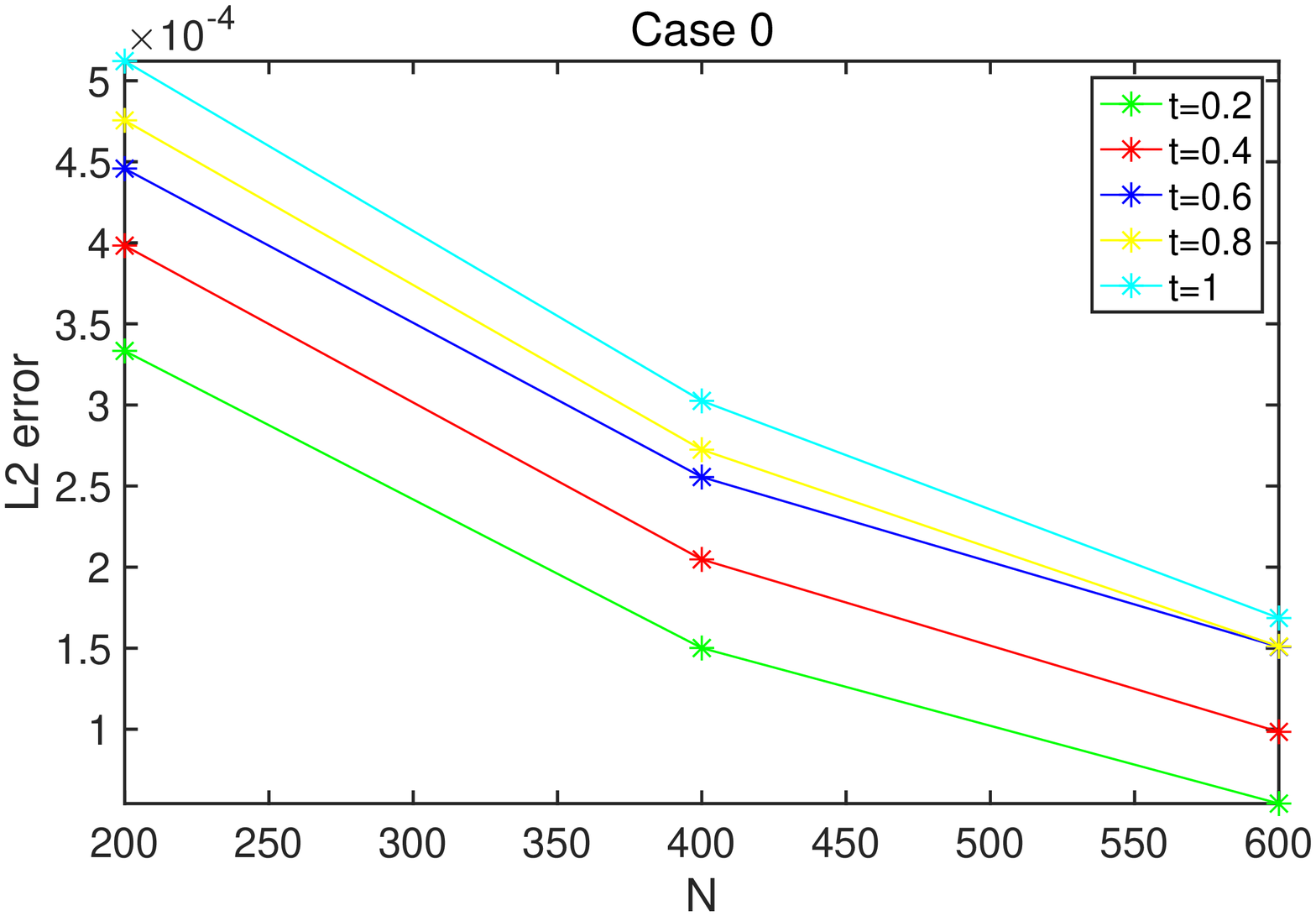}
\includegraphics[width=.45\textwidth]{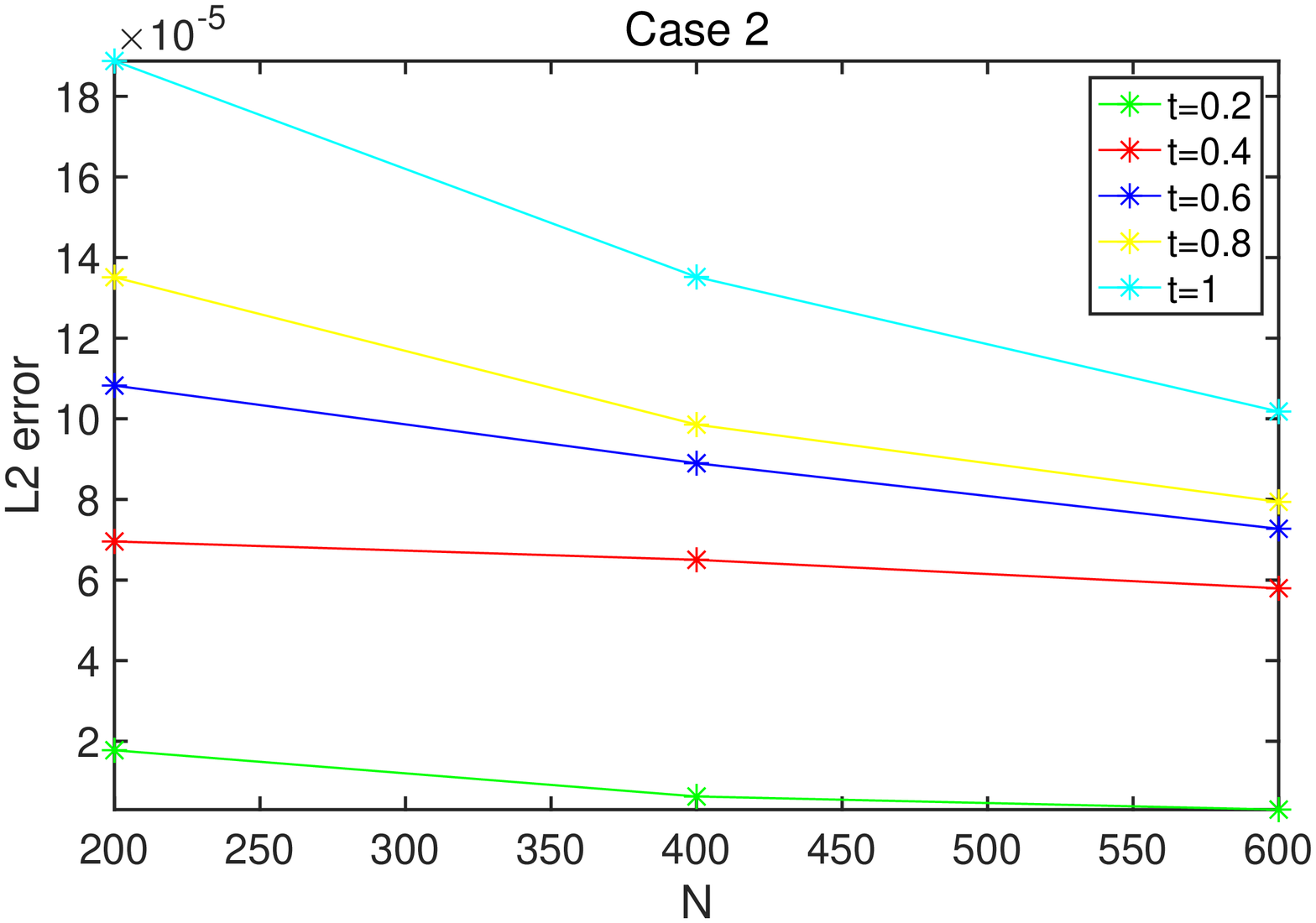}
\caption{Example 3: $L^2$ error of the density when $\beta =0$, $\alpha =0.95$, $\theta =\frac{\pi}{8}$.}
\label{L2plot}
\end{figure}

To further confirm our results, we look at the numerical vorticity profiles at a fixed time with different resolutions, all with $\epsilon = 0.004$. Since the different profiles are always obtained from Case 2 initial data, we only present the results from that case. Figure \ref{Vrefine} displays the results at two different times $T=0.6$ and $T=1$. These results indicate that the solutions converge and the non-uniqueness (indicated by two spirals rather than one) occurs consistently in comparison with the results of lower resolution shown before. We also compare the density solutions on meshes $N=200, 400, 600$ with a reference solution obtained on a refined mesh ($N=1024$). The $L^2$ errors (with respect to the reference solution) of the density in different cases are displayed in Figure \ref{L2plot}, which echoes the convergence of solutions.\\

\begin{figure}[htbp]
\centering
\includegraphics[width=.3\textwidth]{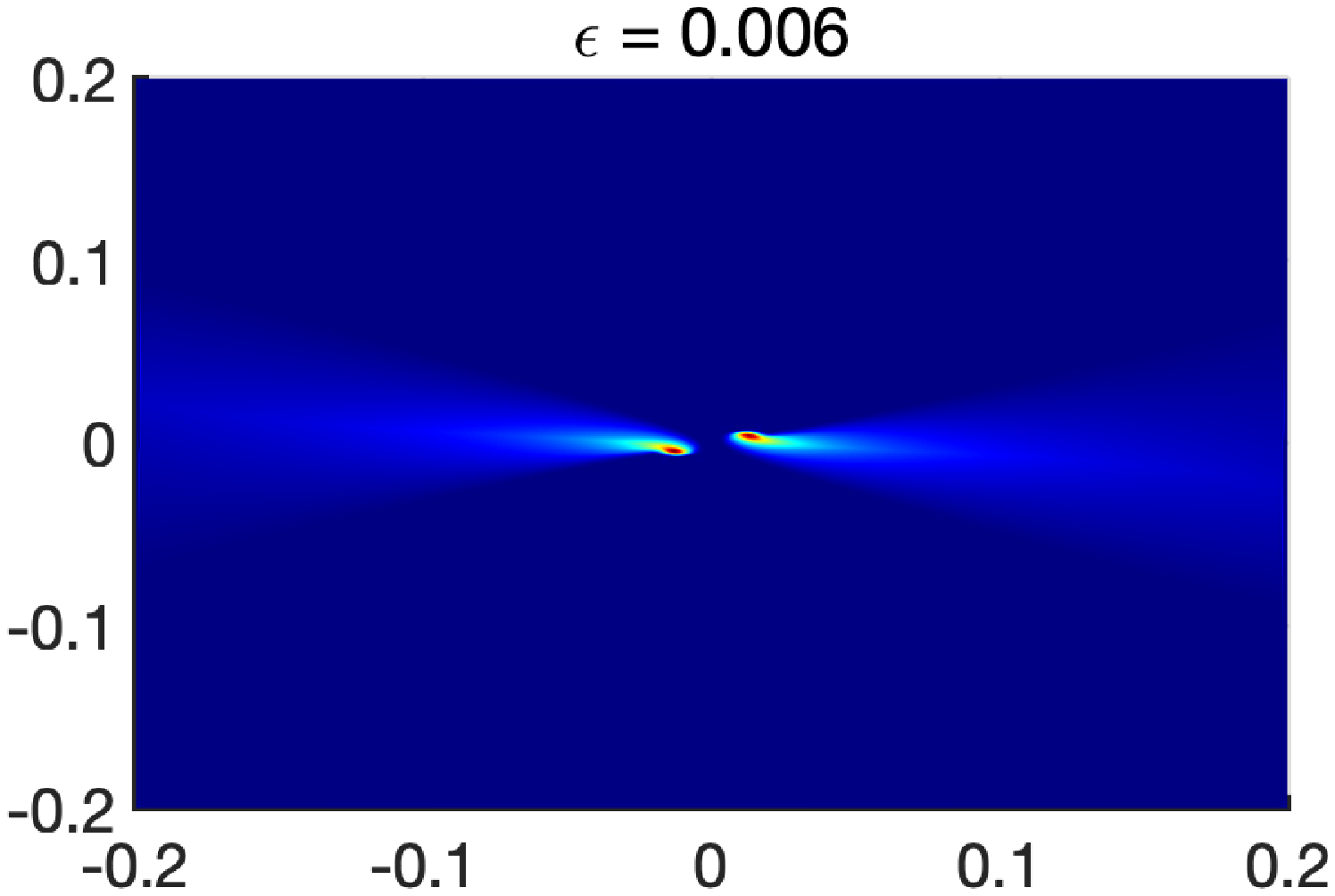}
\includegraphics[width=.3\textwidth]{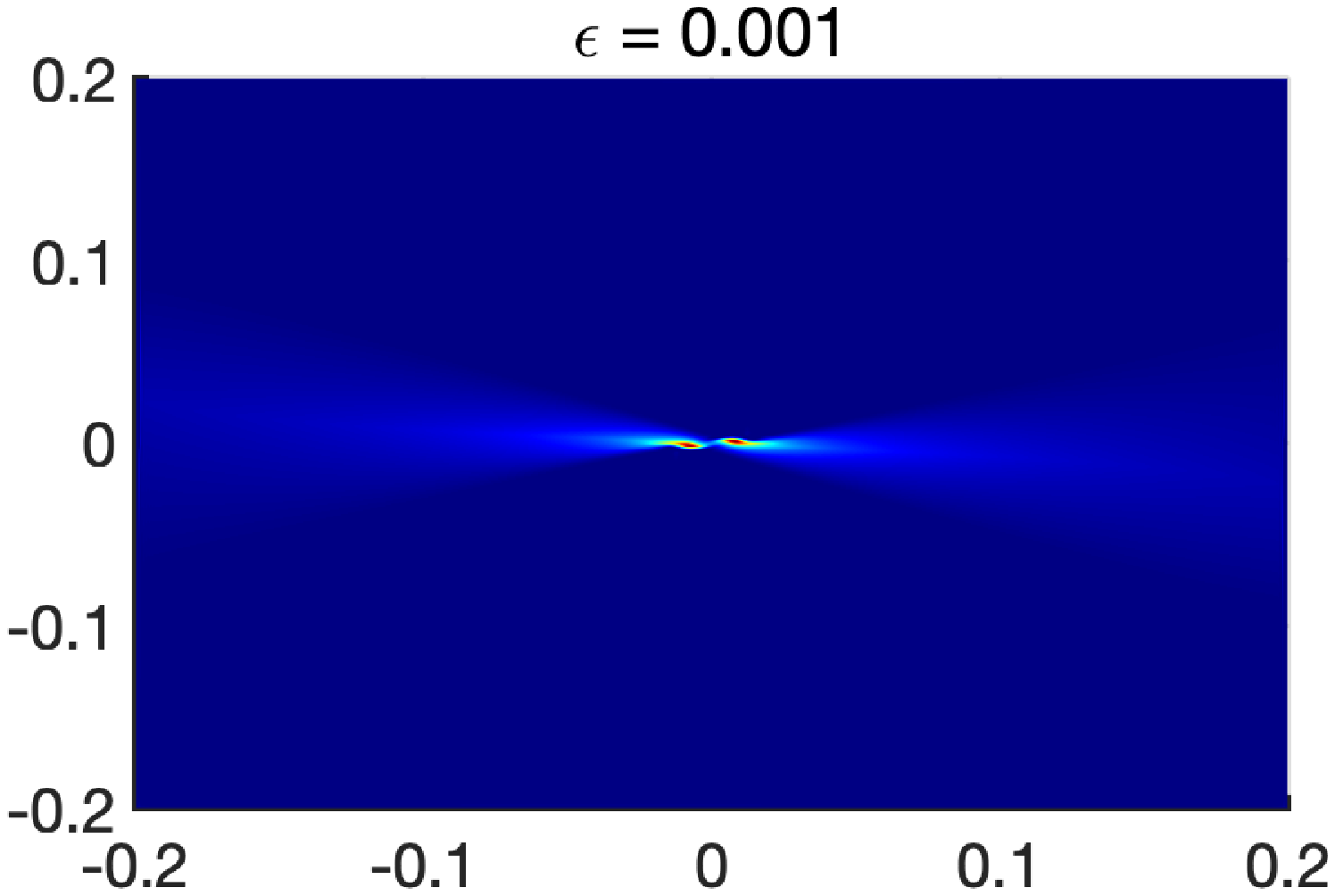}
\includegraphics[width=.3\textwidth]{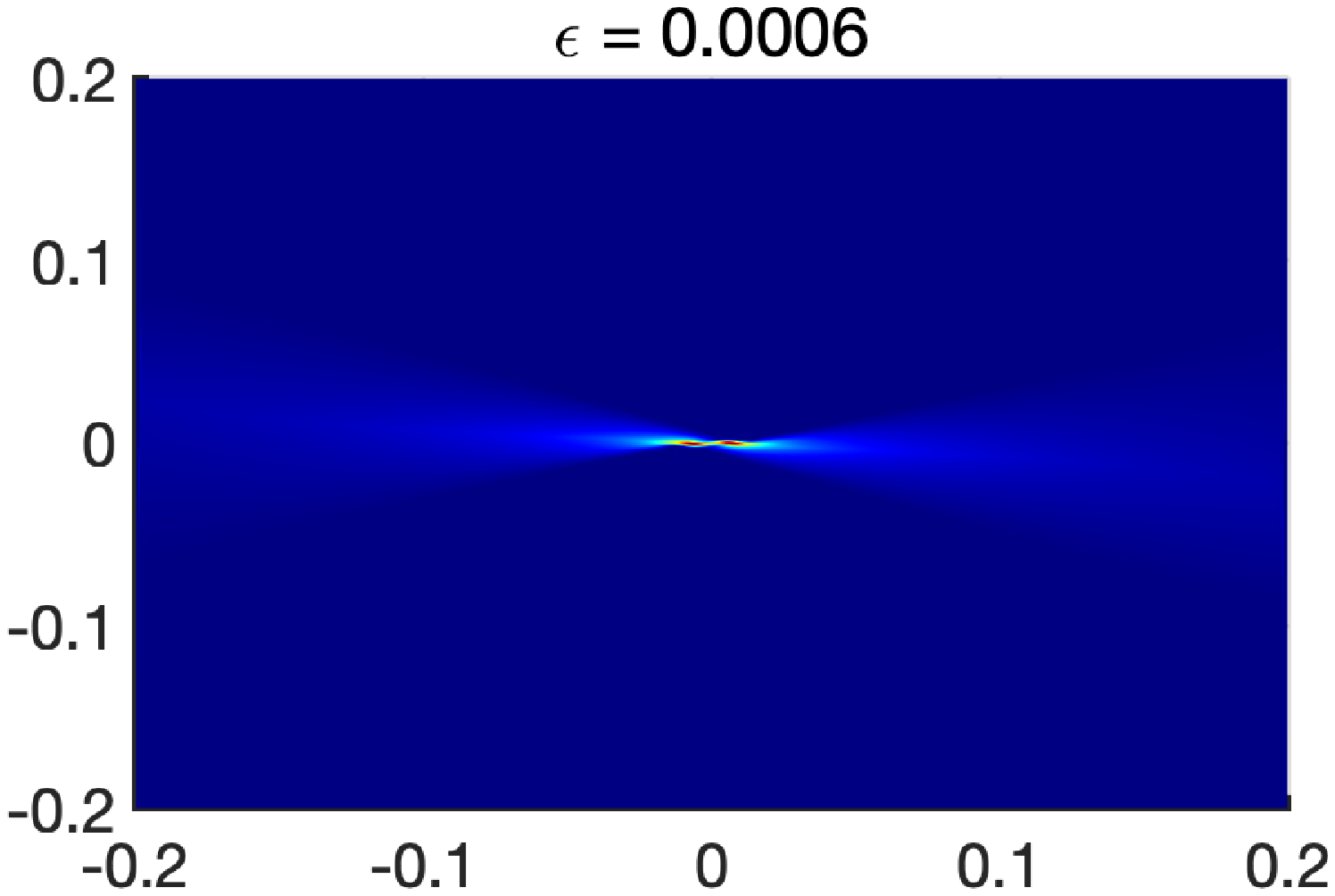}\\
\includegraphics[width=.3\textwidth]{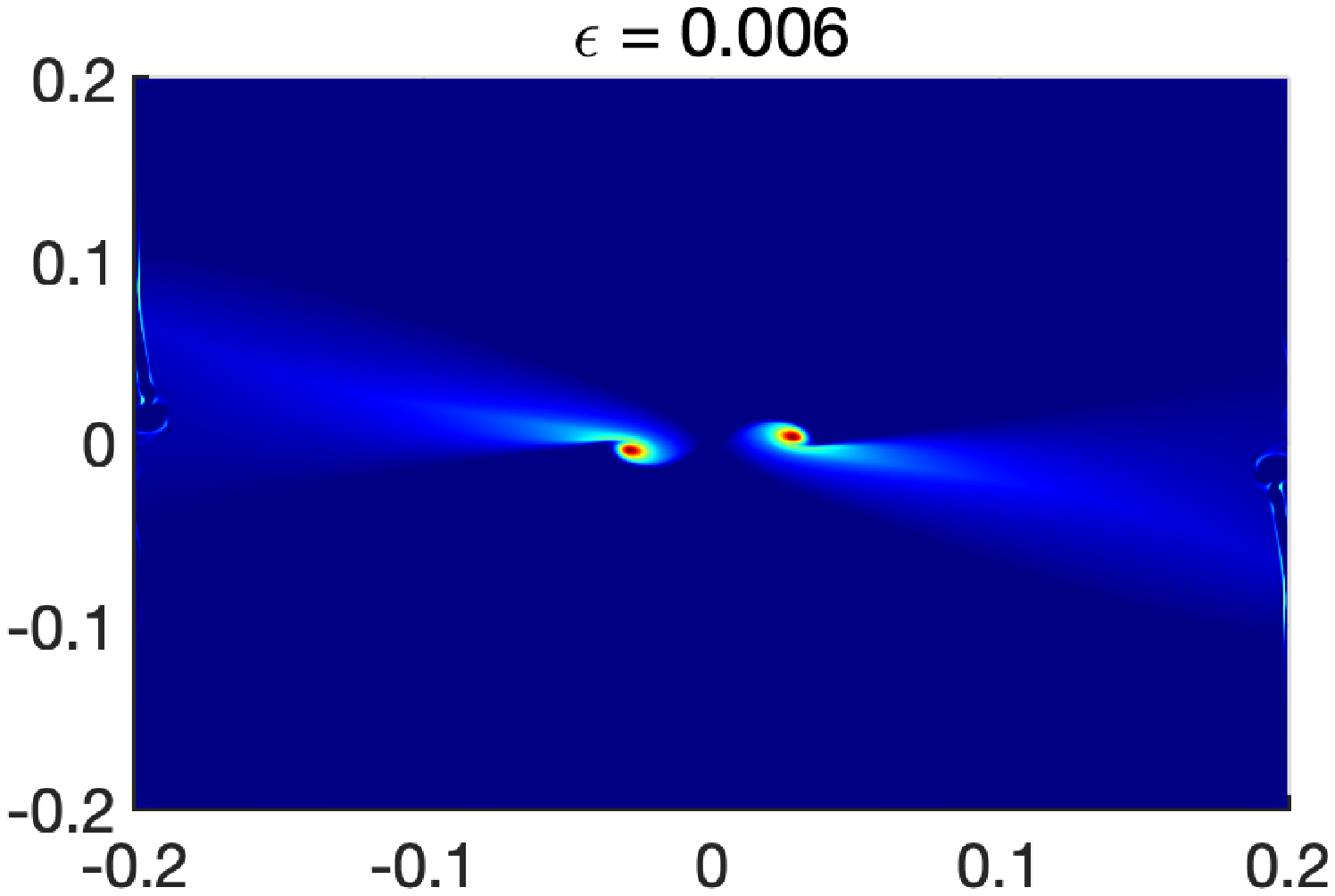}
\includegraphics[width=.3\textwidth]{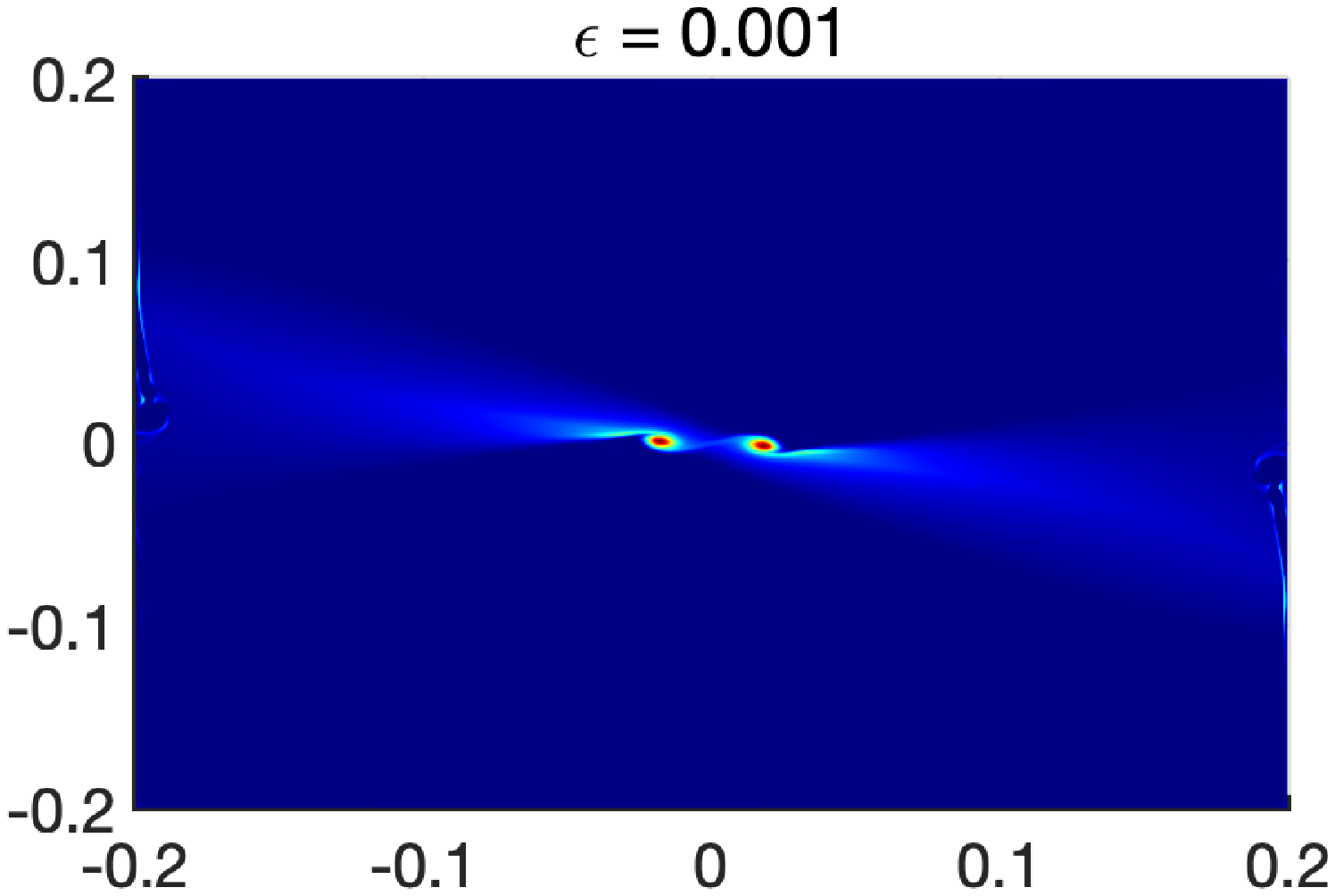}
\includegraphics[width=.3\textwidth]{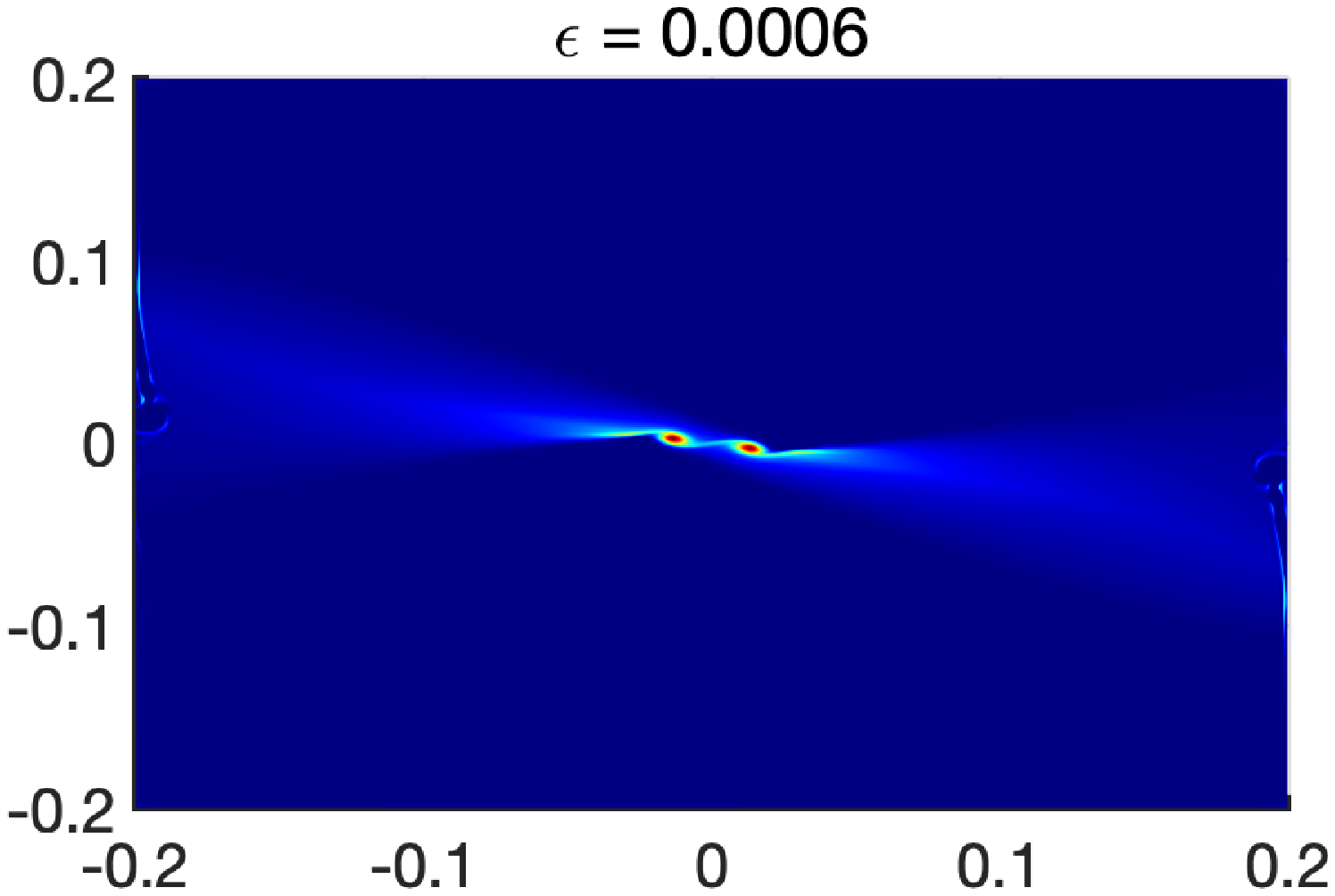}
\caption{Example 3: vorticities obtained from Case 2 initial data with different $\epsilon$ and at different times. $\beta =0$, $\alpha =0.95$, $\theta =\frac{\pi}{8}$. Top: $T=0.2$; bottom: $T=0.6$.}
\label{Ep735}
\end{figure}

Next, we investigate the asymptotic behavior of solutions as the perturbation parameter, $\epsilon$, goes to zero. We plot the numerical vorticity obtained from Case 2 initial data when $N=800$, that is $\Delta x=0.0005$, with three different $\epsilon$'s: 0.006, 0.001, 0.0006. Figure \ref{Ep735} shows the results at two different times $T=0.2$ and $T=0.6$. We can see that the two spirals are always formed and evolve more clearly as time processes. We point out here that since the perturbation parameter $\epsilon$ is used to divide the domain of the piecewise initial data, its value relative to the mesh size makes difference in generating the profile of two spirals. More specifically, the two spirals are shown more clearly when the difference between the mesh size and the perturbation parameter value is larger. This can be observed in both Figure \ref{Vrefine} and \ref{Ep735}. Meanwhile, the two spirals are formed more slowly when $\epsilon$ is closer to the mesh size, as shown in Figure \ref{Ep735}.\\

\begin{figure}[htbp]
\centering
\includegraphics[width=.45\textwidth]{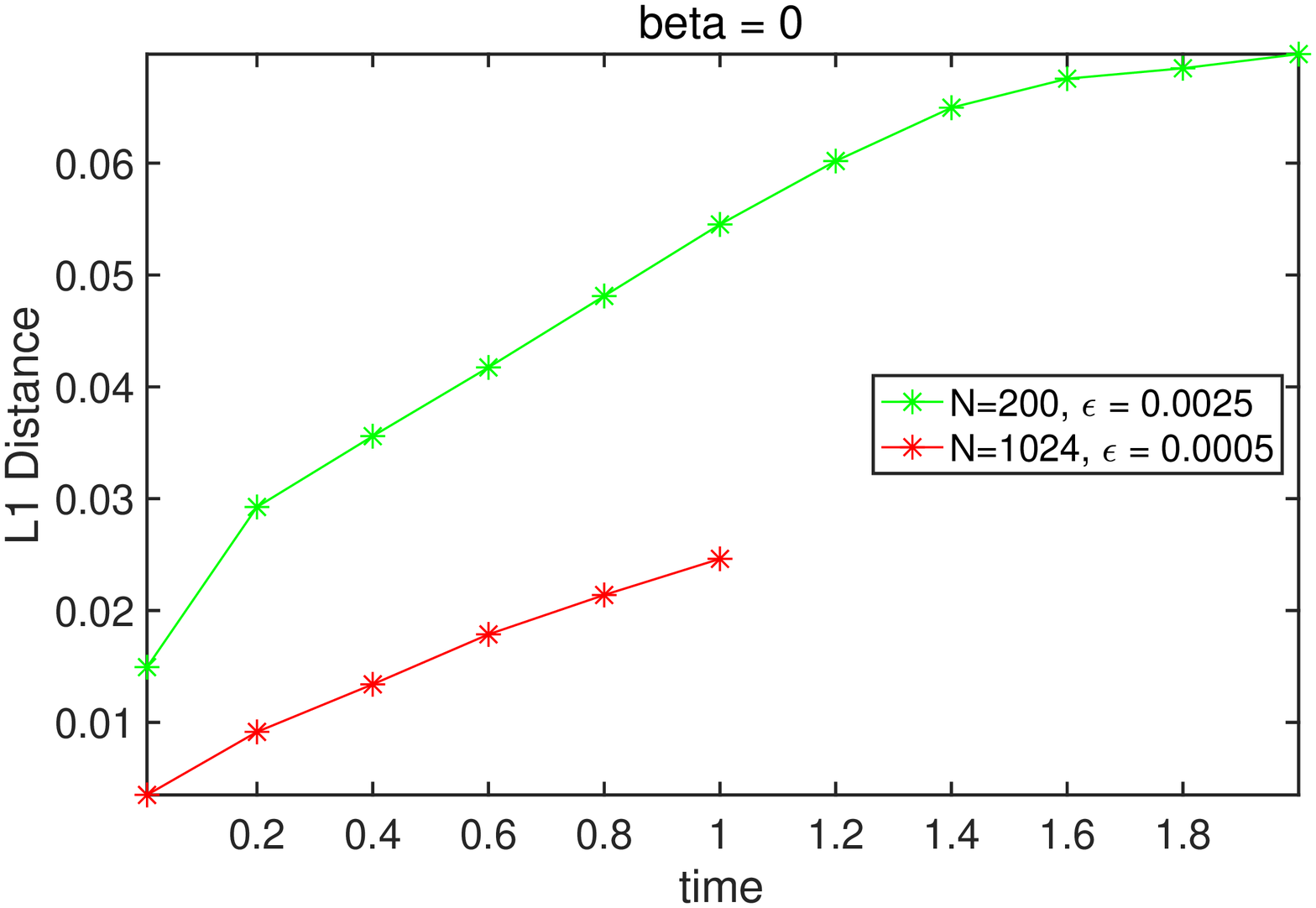}
\includegraphics[width=.45\textwidth]{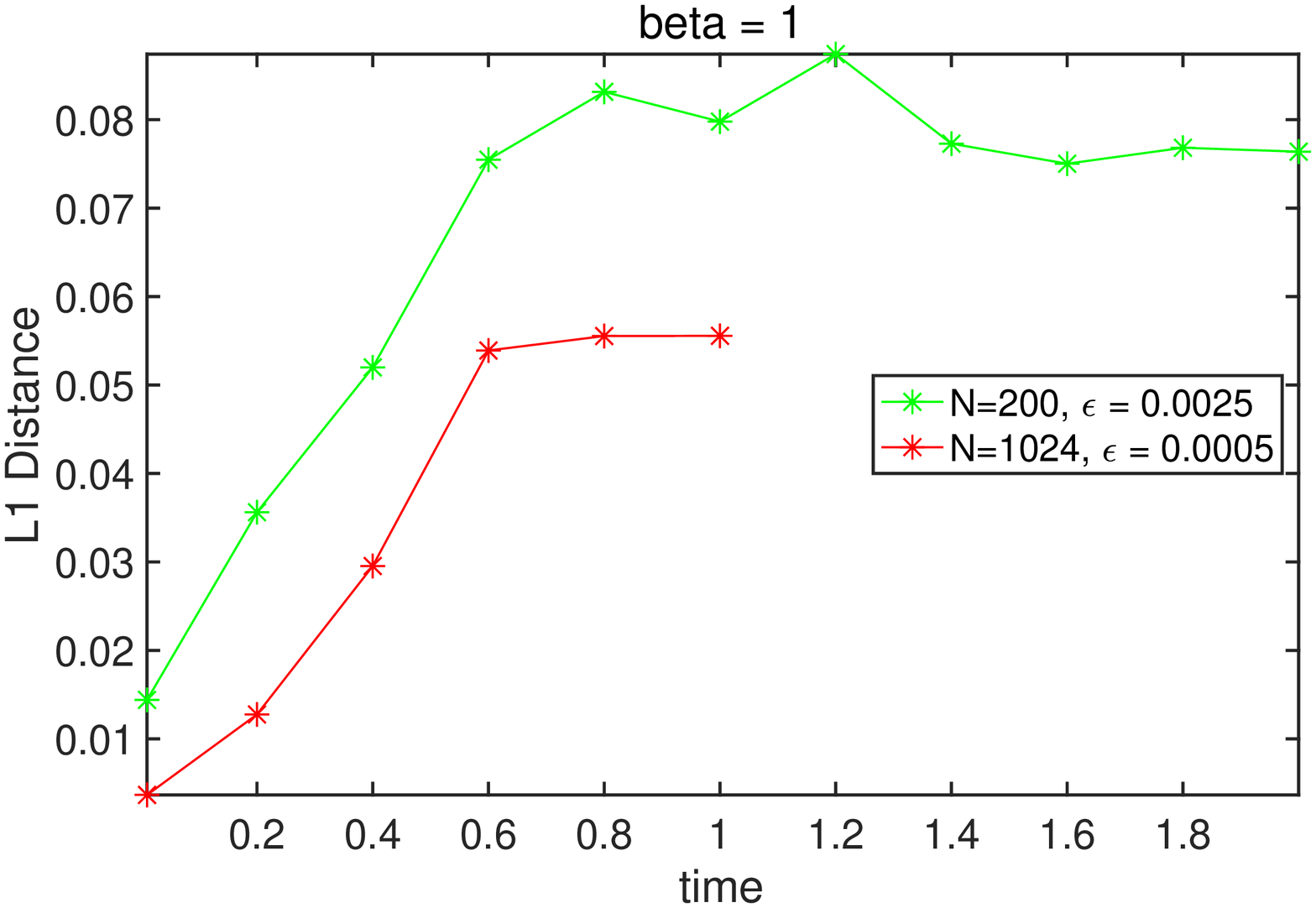}
\caption{Example 3: Behavior of $D^{\epsilon}(t)$ for different mesh resolutions. Left: expected non-unique case; right: expected unique case. The simulation results at larger time ($T>1$) for the refined mesh ($N=1024$) are unavailable due to the limited computation resources.}
\label{DisTime}
\end{figure}

We further quantify the difference between the solutions by defining the following metric
\begin{align*}
    D^{\epsilon}(t)=\| \omega ^{\epsilon, 0}(\cdot, \cdot, t)-\omega ^{\epsilon, 2}(\cdot, \cdot, t)\|_{L^1_{\Omega}}, 
\end{align*}
where $\omega ^{\epsilon, i}(x, y, t)$ is the vorticity solution obtained from Case $i$ initial data, $i=0, 2$. At $t=0$, as $\epsilon \rightarrow 0$, $\omega ^{\epsilon, 0}$ and $\omega ^{\epsilon, 2}$ both converge to $\bar{\omega}$ defined in Section 3, therefore we have $D^{\epsilon}(0)\rightarrow 0$.  For $t>0$, with $\epsilon$ fixed and small enough, we would expect that for non-unique solutions, this metric keeps increasing as time evolves; but for the unique case, it will remain small for positive $t$. We perform two groups of computations. One is on $200^2$ grids with $\epsilon =0.0025$ and the other is on $1024^2$ grids with $\epsilon \approx 0.0005$. In Figure \ref{DisTime}, we plot the defined metric as a function of time and compare the results when $\beta =0$ and $\beta =1$, respectively. On both meshes, same trends have been observed; when $\beta =0$, the metric increases in time, which is consistent with the expectation of the non-uniqueness; when $\beta =1$, the metric shows an initial increase and then stays almost unchanged, which is expected for the unique case while the initial increase is likely due to the numerical errors and nonlinear dependence of the time variable.
We remark here that a different comparison can be made by studying the behavior of the metric $D^{\epsilon}(t)$ at a fixed time as $\epsilon$ approaches zero. However, our experiments show that the metric appear to be sensitive to the specific form of the solutions as well as the effect the perturbation parameter has in the data. Such comparison may be helpful when the perturbation parameter takes values very close to zero, which, however, is not feasible in the present setting with the special structure of the initial data.\\

One may wonder if there are other values of $\beta$ that could lead to different solutions. 
With further tests on more different values of $\beta$, we conclude that for a fixed (carefully chosen) value of $\alpha$ and $\theta_0$ in the initial data, there exists a $\beta ^* >0$ such that

\begin{figure}[!htbp]
\centering
\includegraphics[width=.3\textwidth]{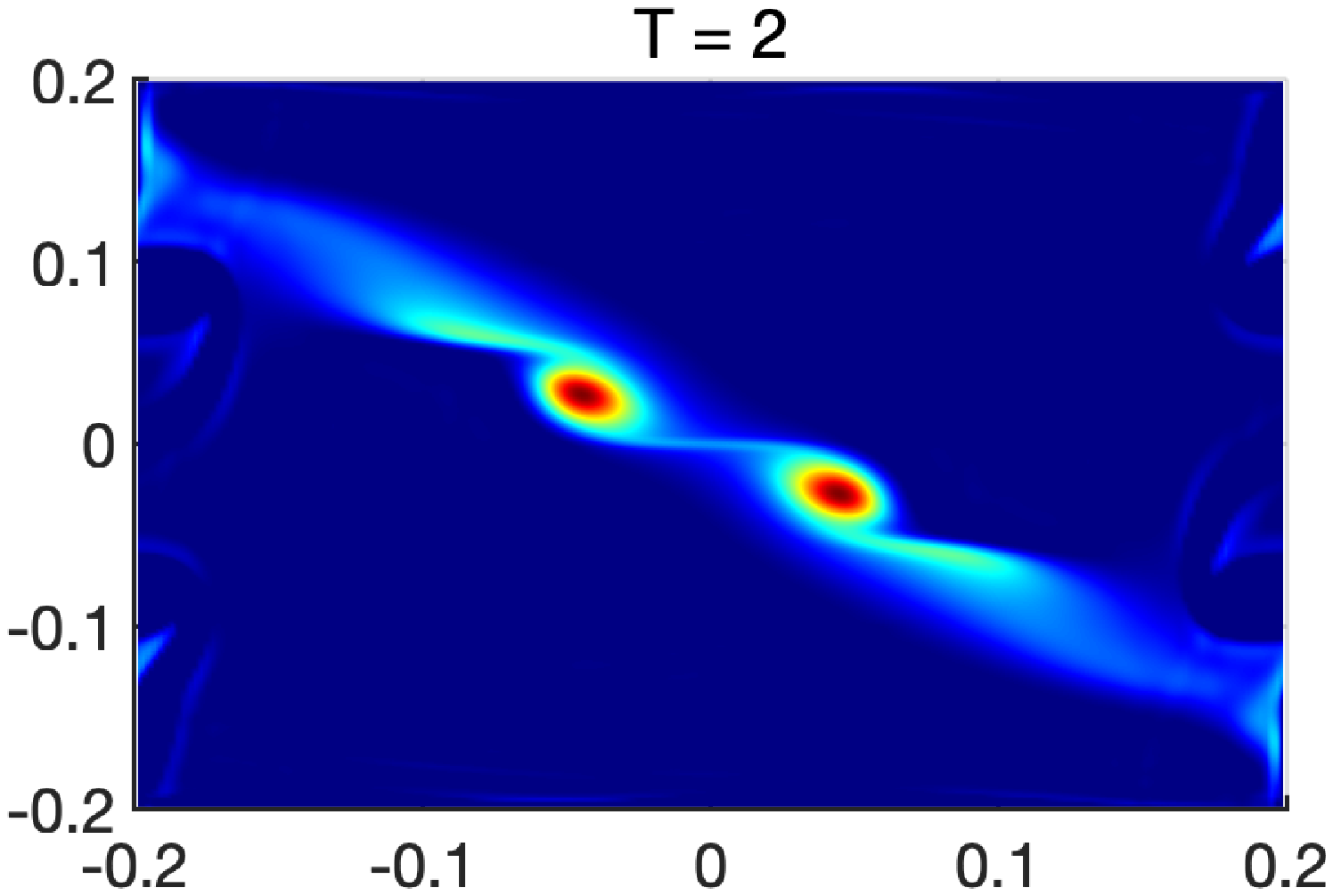}
\includegraphics[width=.3\textwidth]{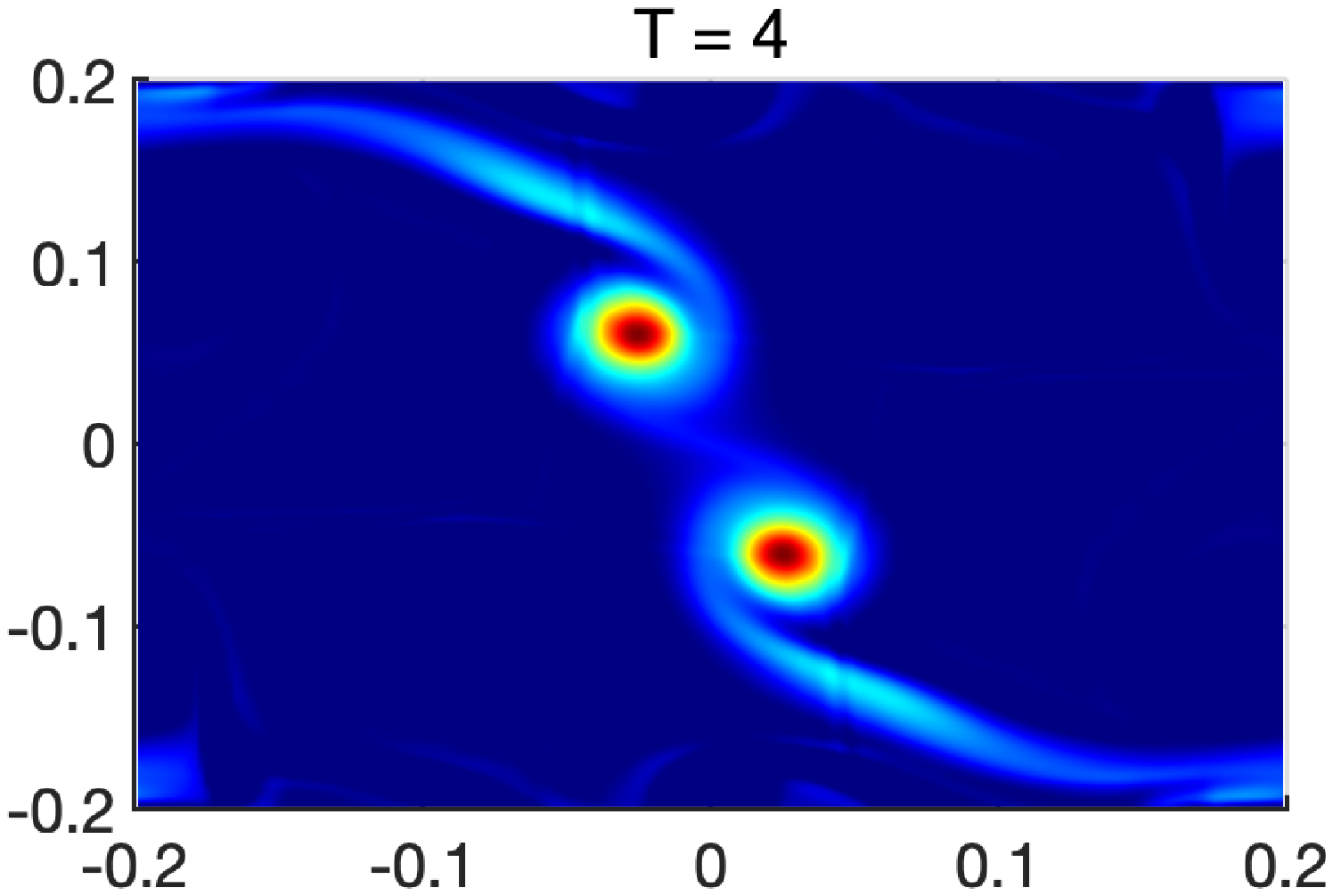}
\includegraphics[width=.3\textwidth]{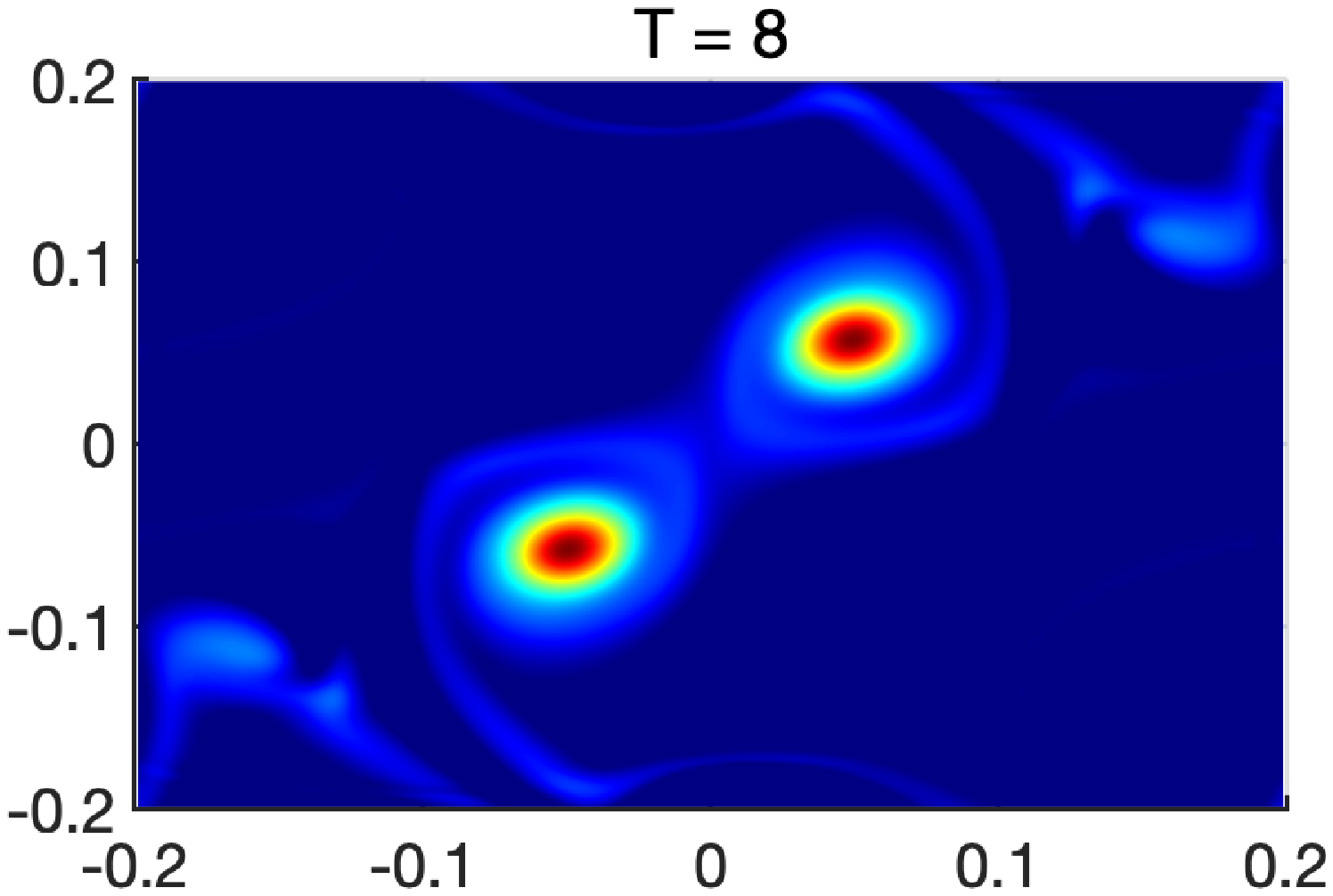}
\caption{Case 2 vorticity profiles at different times with $\beta = 0.1$. $\alpha =0.95$ and $\theta _0 = \frac{\pi}{8}$. }
\label{EX411}
\end{figure}

\begin{itemize}
\item[(1)]when $0\leq \beta < \beta ^*$, non-unique solutions can be observed for (\ref{eq:2DEuler}), which indicates that the continuous dependency of initial data is violated and it is ``incurable";
\item[(2)]when $\beta \geq \beta ^*$, the initial data for both Case 0 and Case 2 lead to vorticities as one single spiral eventually.
\end{itemize}
For instance, among three values of $\beta$: $0.1$, $0.2$ and $0.3$, it has been observed that only for $\beta =0.1$, the non-uniqueness occurs consistently as time processes. We present the results obtained from Case 2 for $\beta =0.1$ in Figure \ref{EX411}.\\

\subsubsection{Effects of $\alpha$}\label{sec:alpha}
In this section, we test how the choice of $\alpha$, the exponent in the initial vorticity function, affects the solutions to the system, with three different initial data respectively. Based on the results in the previous section, we fix $\beta =0.1$ and $\theta _0 = \frac{\pi}{8}$. In the following, we first present the results for different choices of $\alpha$ and then draw a conclusion.\\

\begin{figure}[!htbp]
\centering
\includegraphics[width=.3\textwidth]{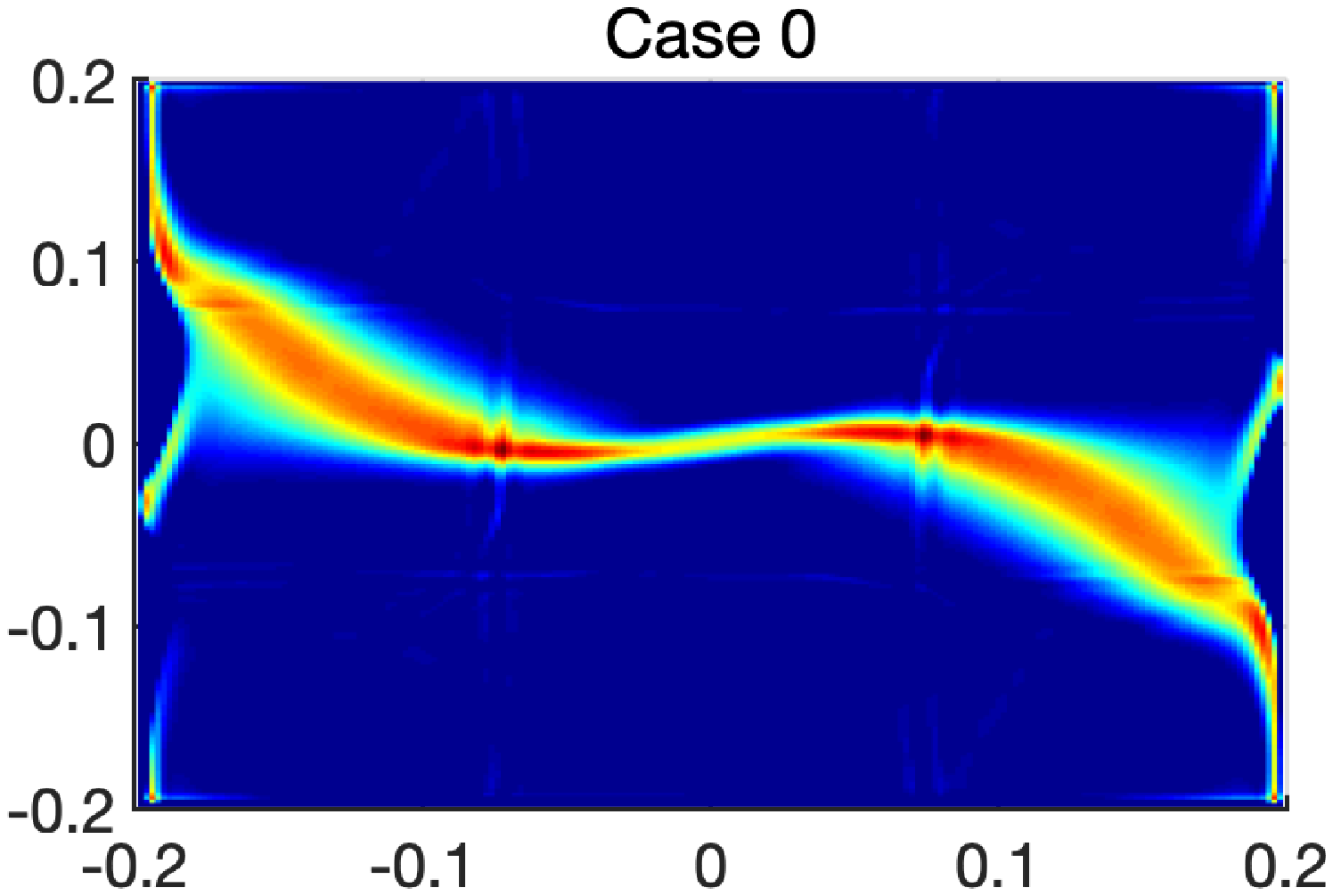}
\includegraphics[width=.3\textwidth]{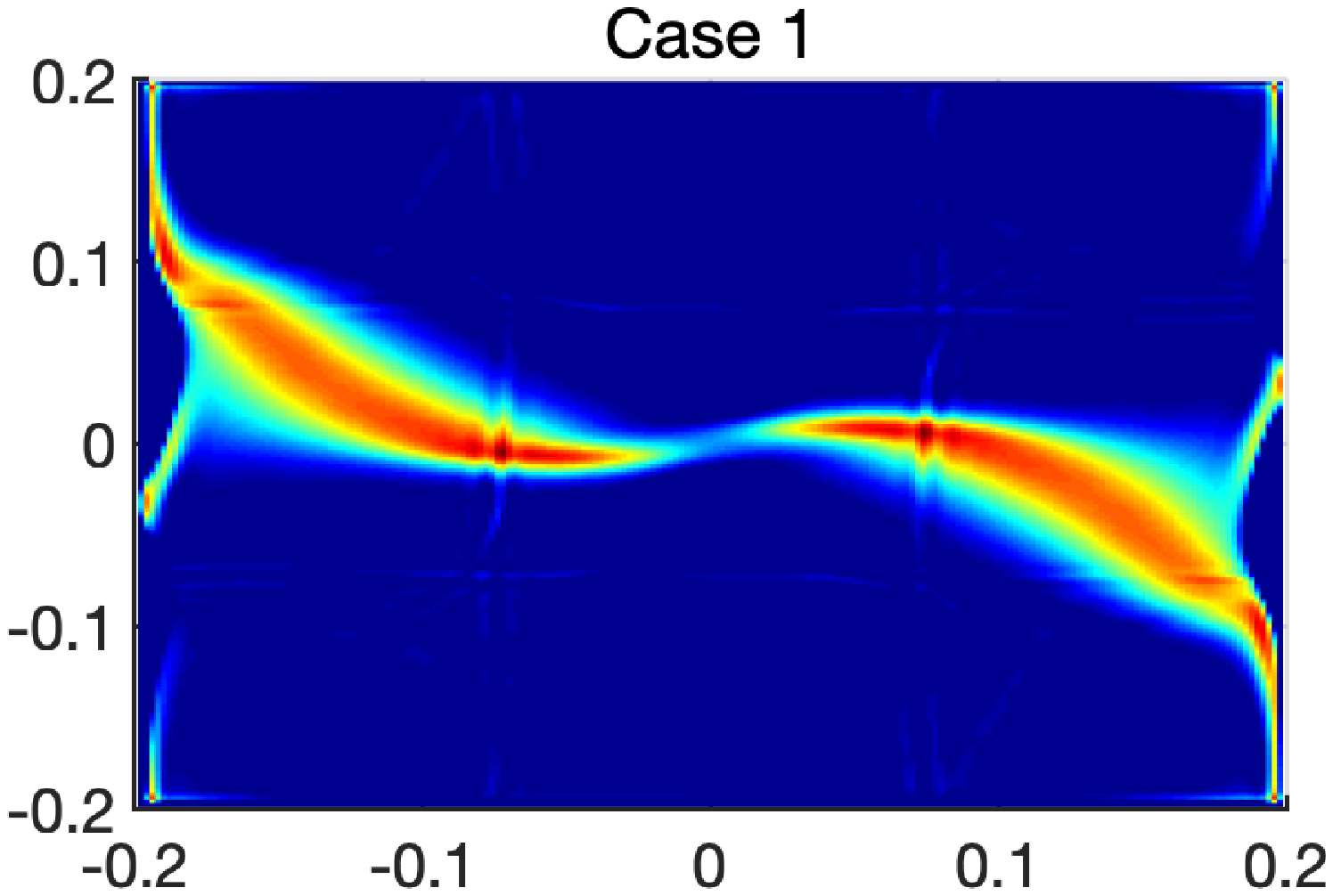}
\includegraphics[width=.3\textwidth]{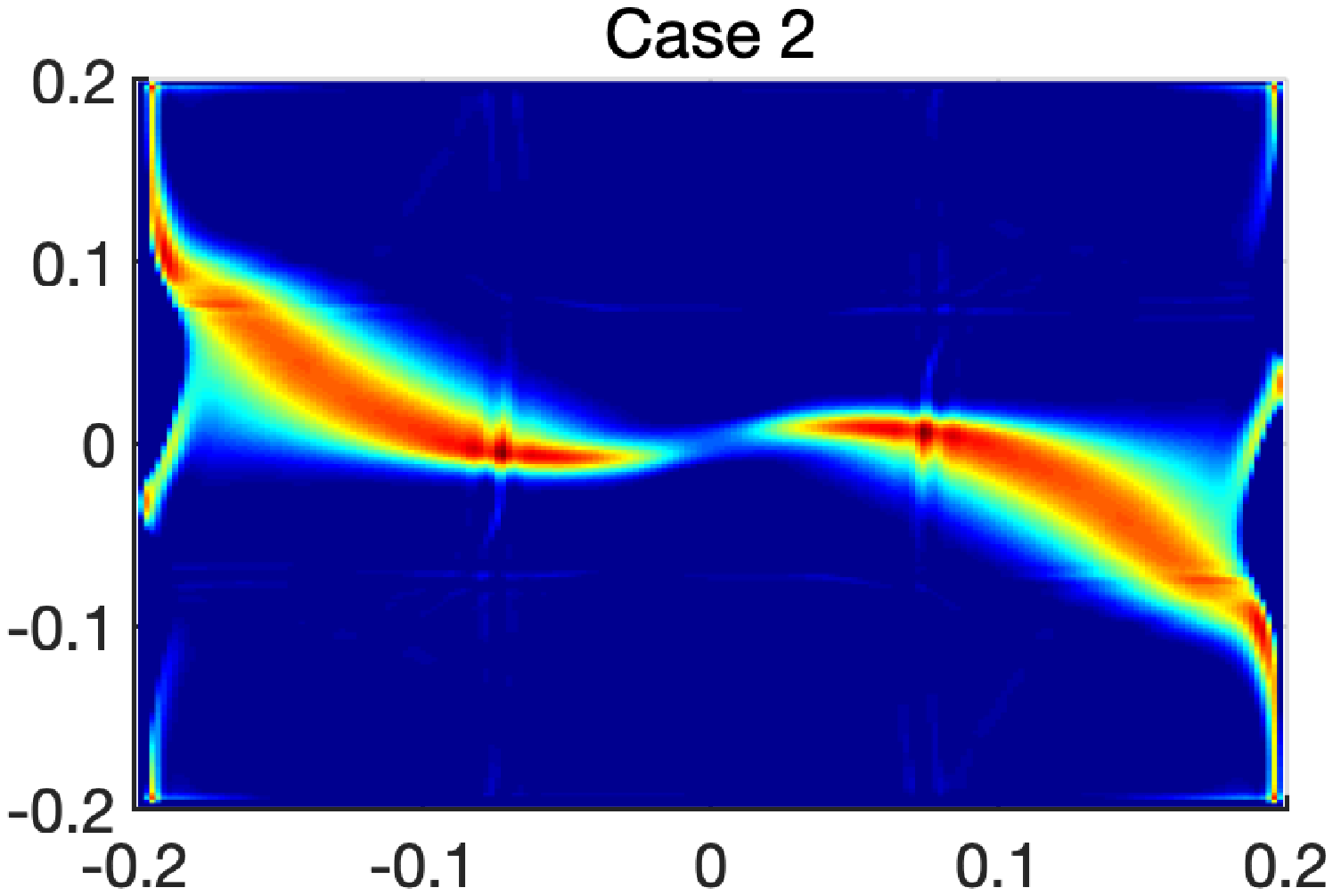}
\caption{Example 4: vorticity profiles at $T=6$. $\beta = 0.1$, $\alpha =0.1$ and $\theta _0 = \frac{\pi}{8}$. }
\label{EX51}
\end{figure}

\noindent \textit{Example 4.} $\alpha = 0.1$. 
The results at $T=6$ show that the solutions in all three cases look very similar except at the singular center. See Figure \ref{EX51}.\\

\begin{figure}[!htbp]
\centering
\includegraphics[width=.3\textwidth]{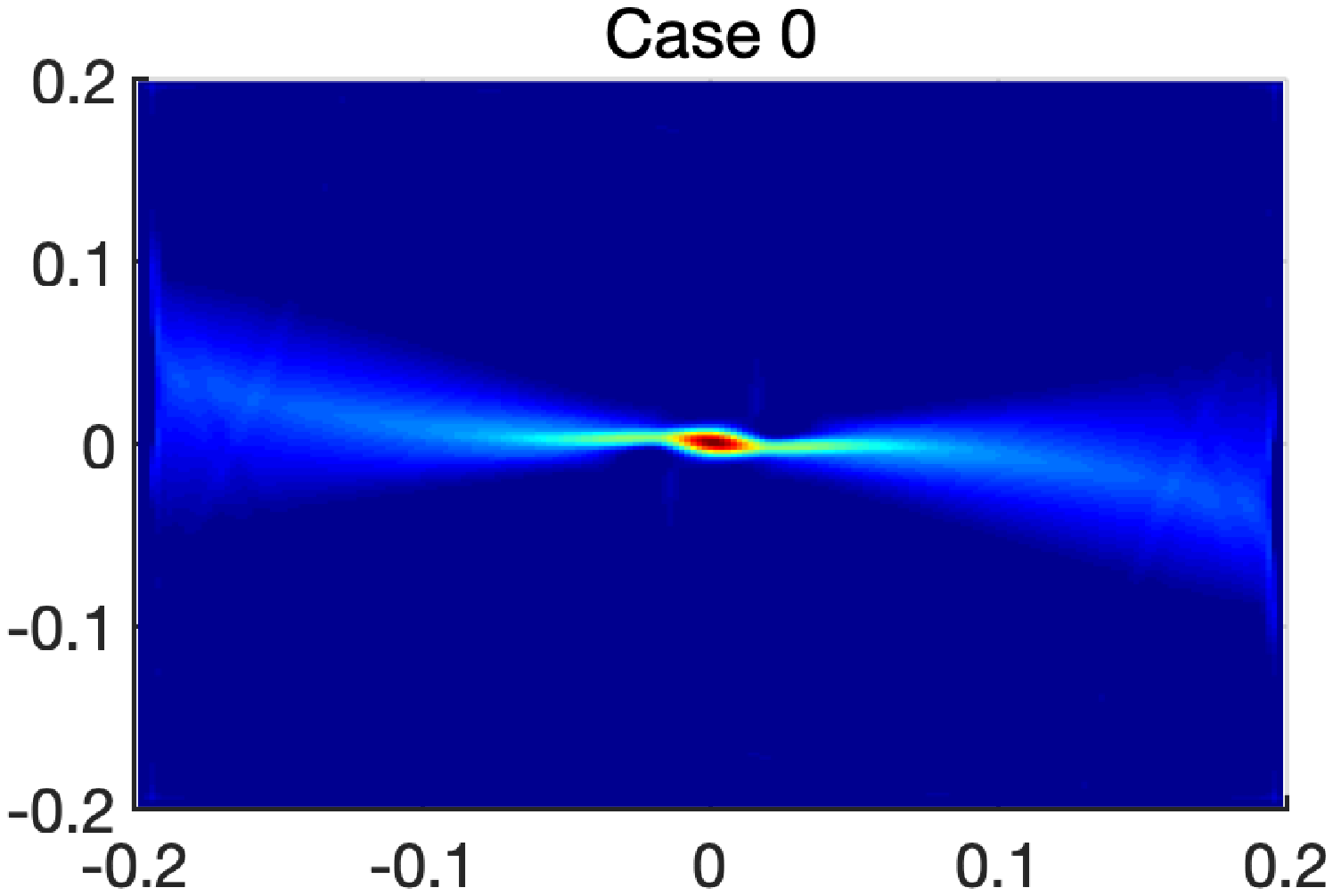}
\includegraphics[width=.3\textwidth]{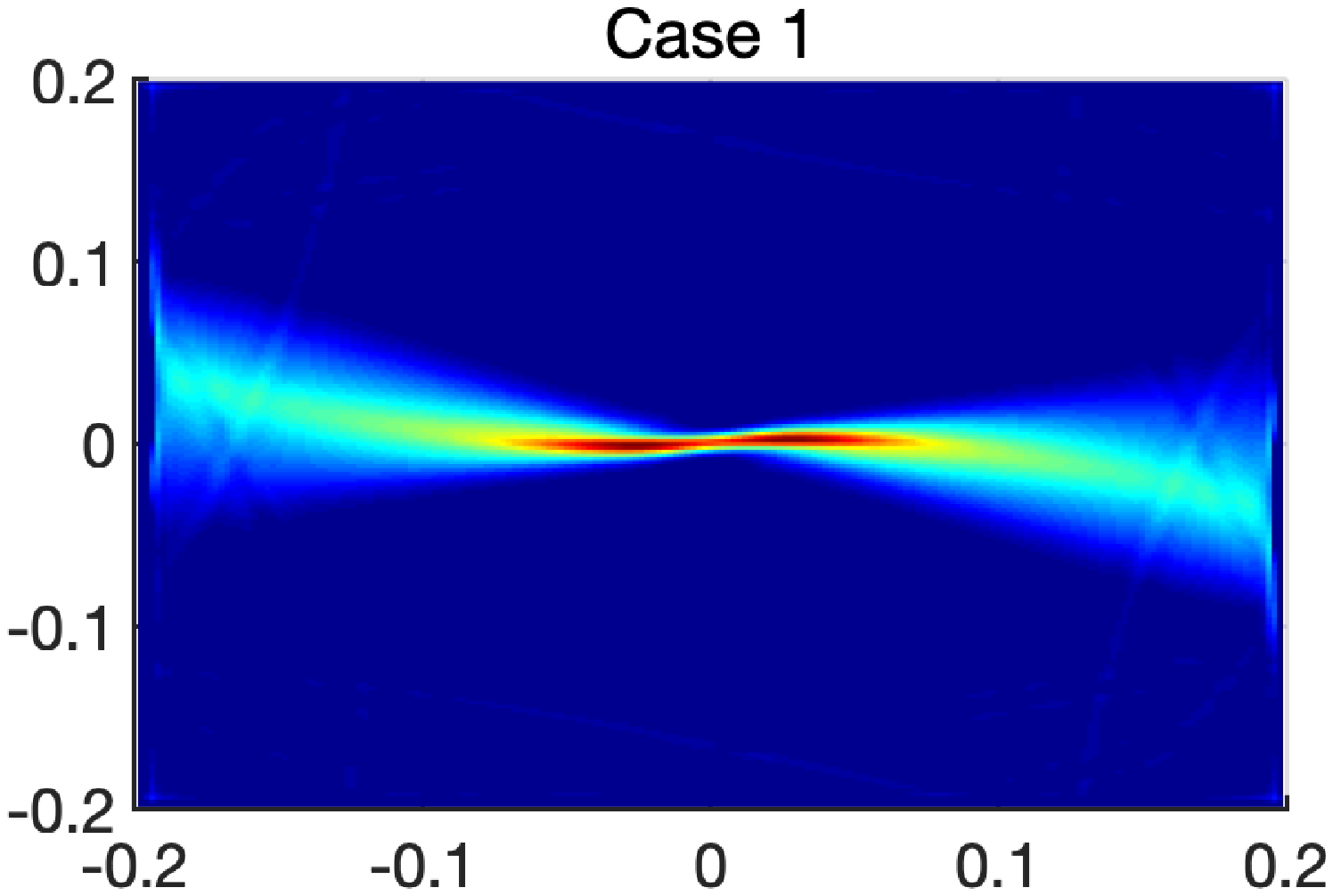}
\includegraphics[width=.3\textwidth]{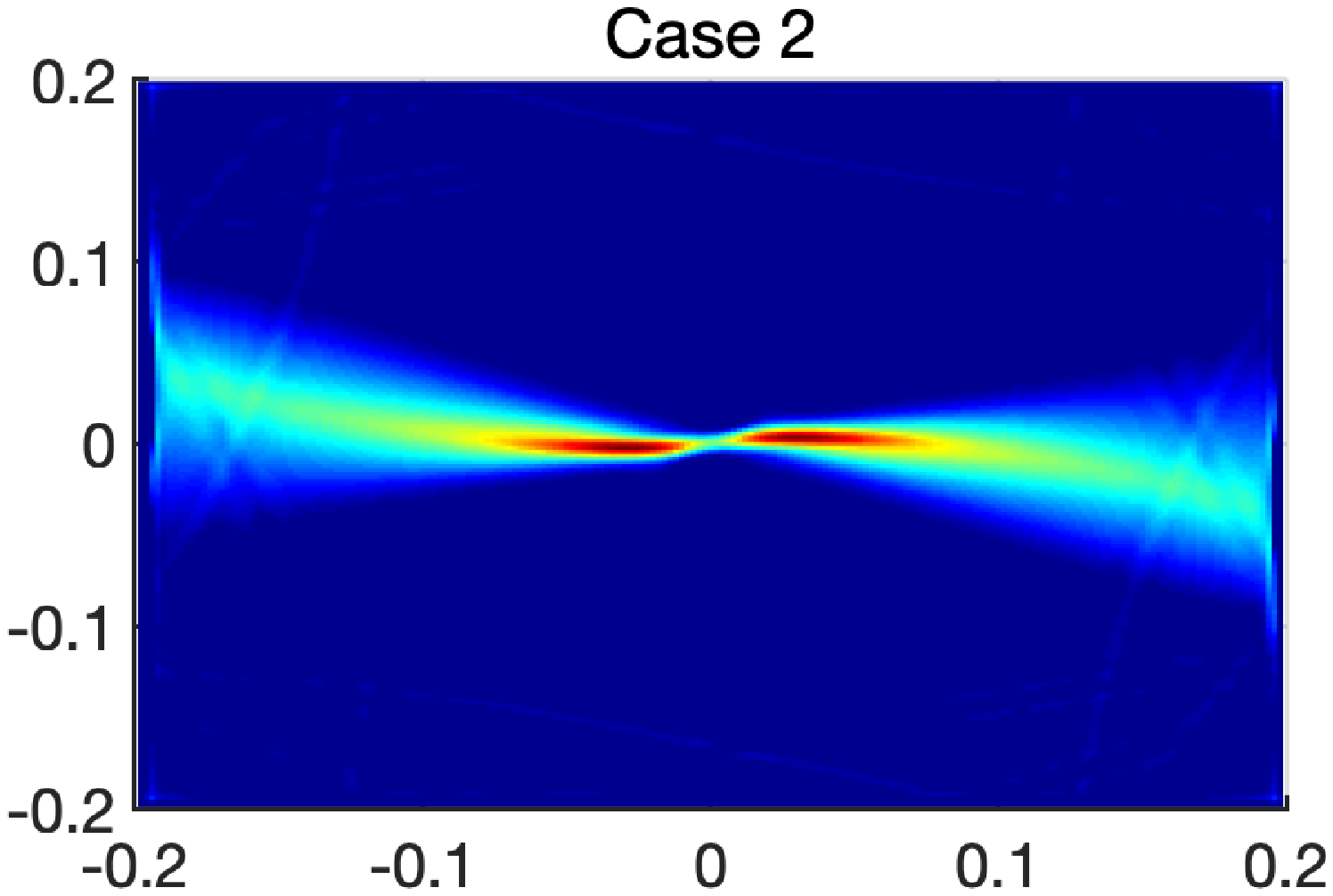}\\
\includegraphics[width=.3\textwidth]{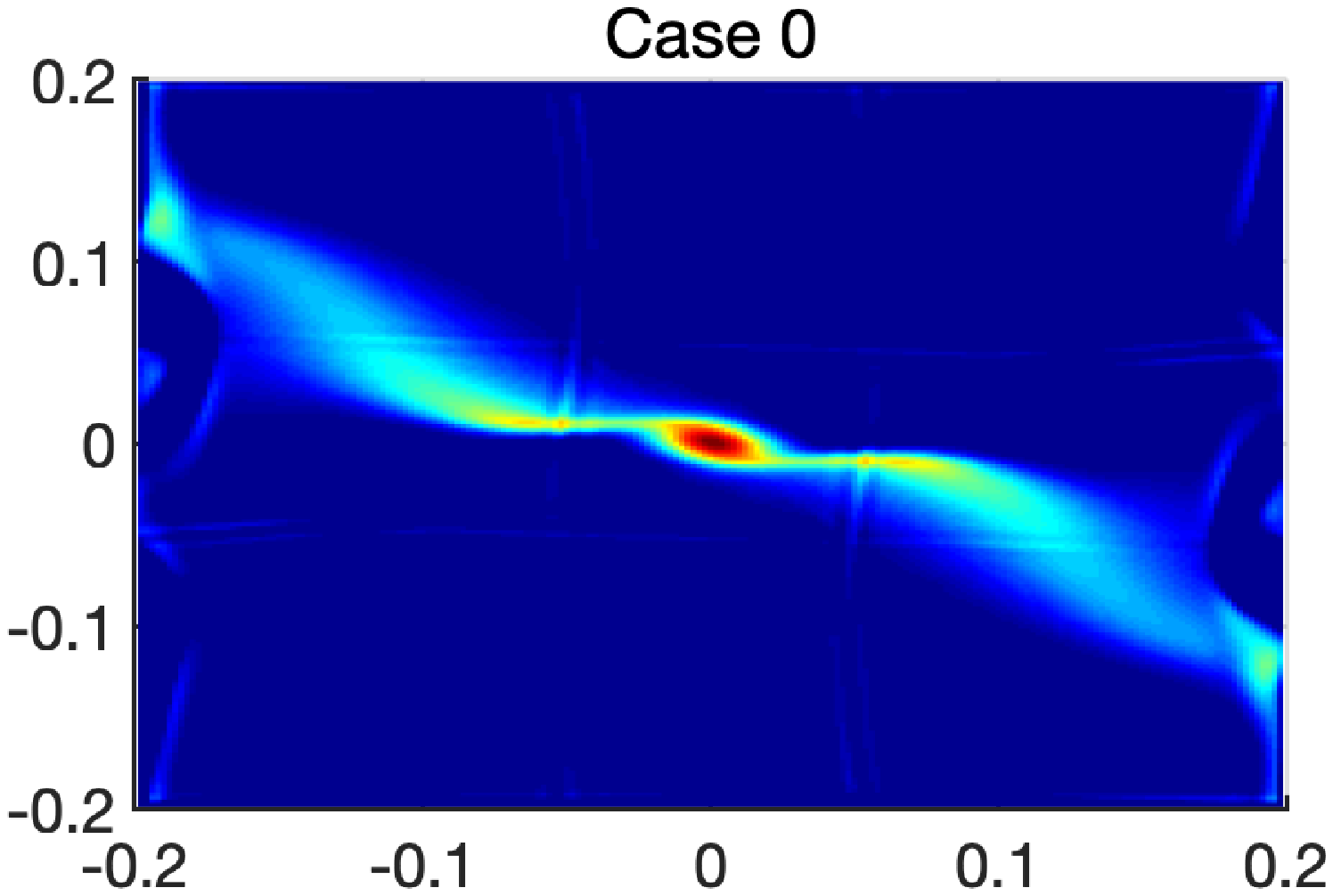}
\includegraphics[width=.3\textwidth]{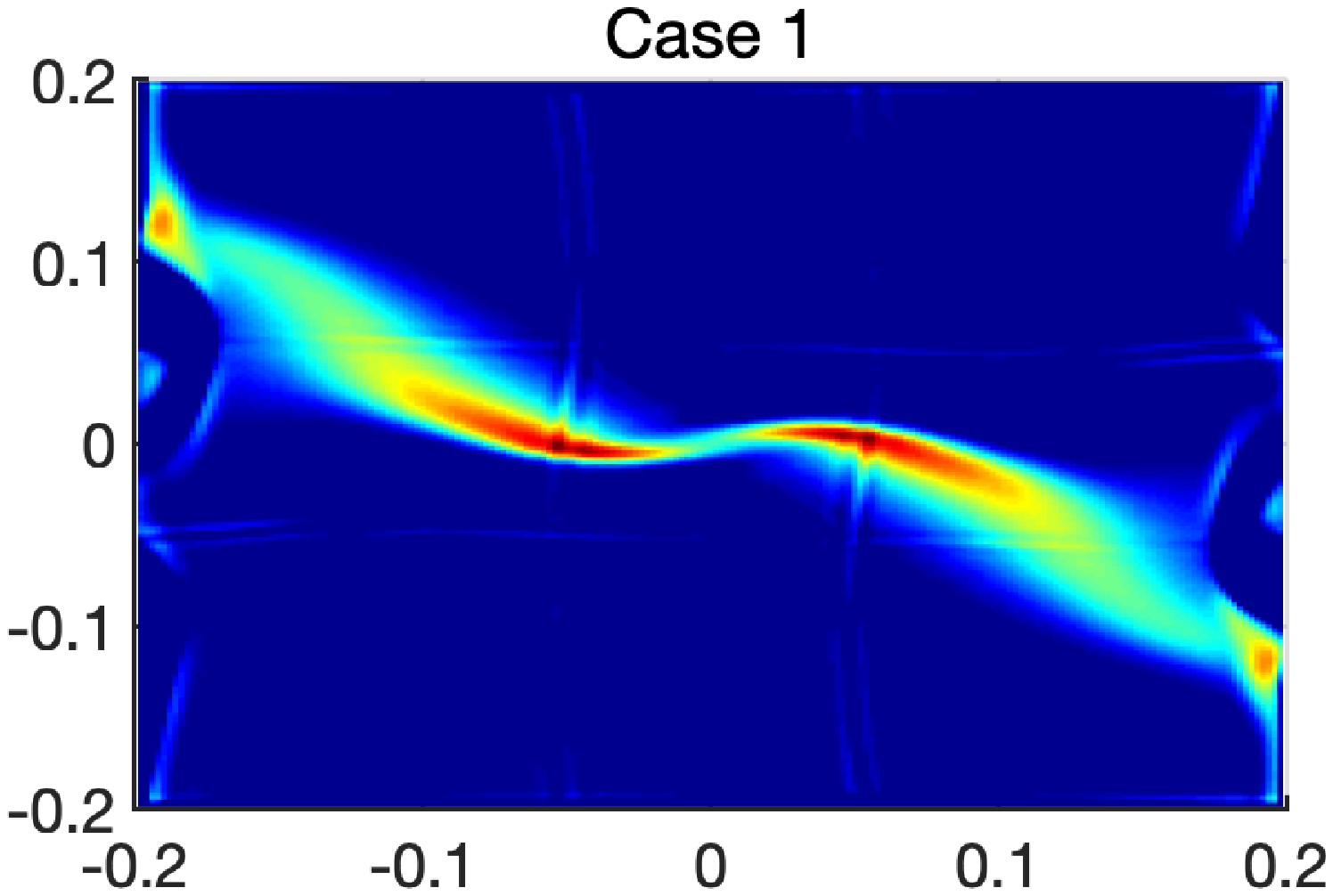}
\includegraphics[width=.3\textwidth]{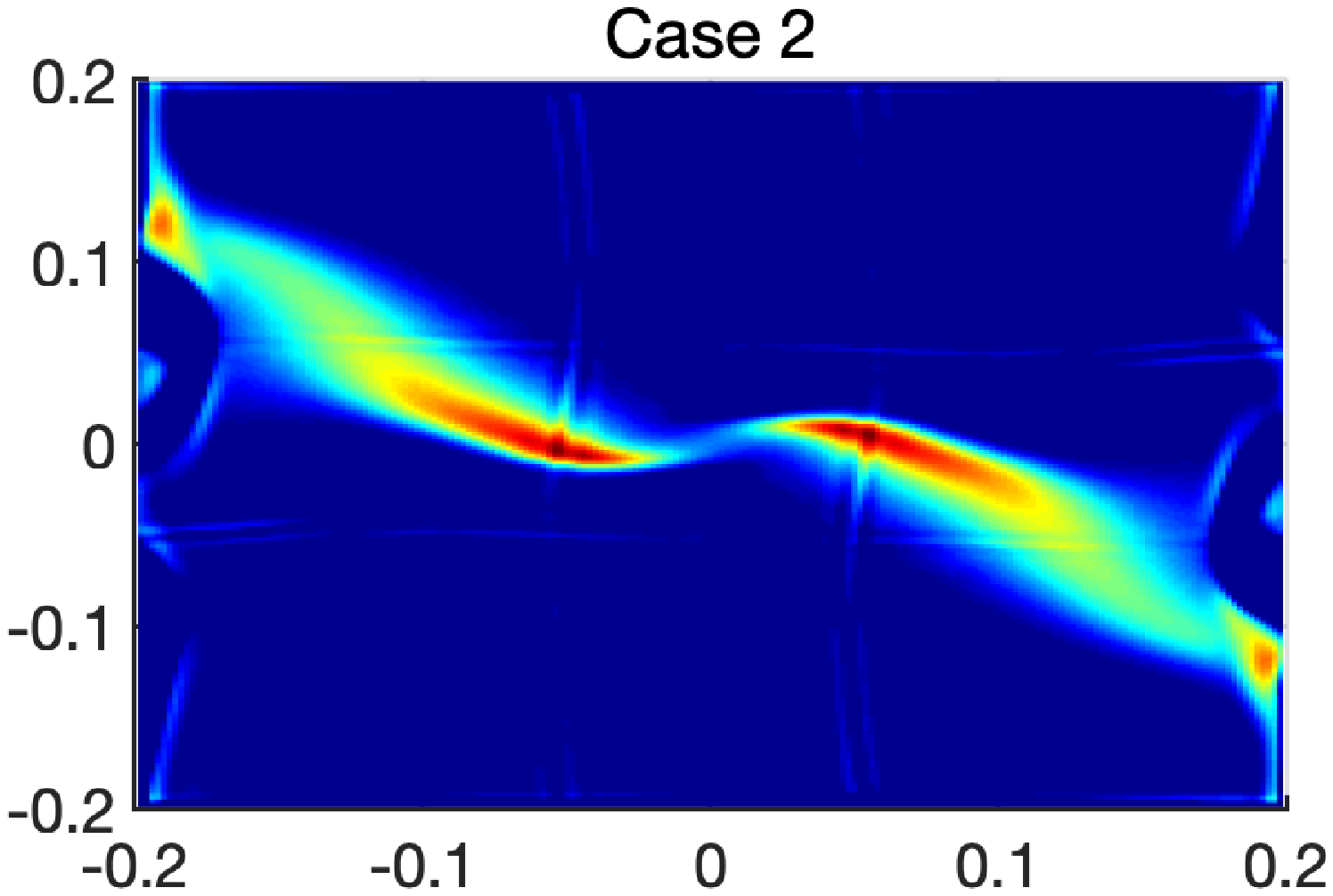}
\caption{Example 5: vorticity profiles at different times. $\beta = 0.1$, $\alpha =0.5$ and $\theta _0 = \frac{\pi}{8}$. Top: T=1; bottom: T=4.}
\label{EX6}
\end{figure}

\noindent \textit{Example 5.} $\alpha = 0.5$. It's clear that Case 0 forms a single spiral at $T=1$, while Case 1 and Case 2 sharing common features look different from Case 0. We then look at their behaviors at larger time $T=4$, and the results show that the two peaks in both Case 1 and Case 2 are more distant as time processes, which distinguishes them from Case 0.
See Figure \ref{EX6}.
\\

\begin{figure}[!htbp]
\centering
\includegraphics[width=.3\textwidth]{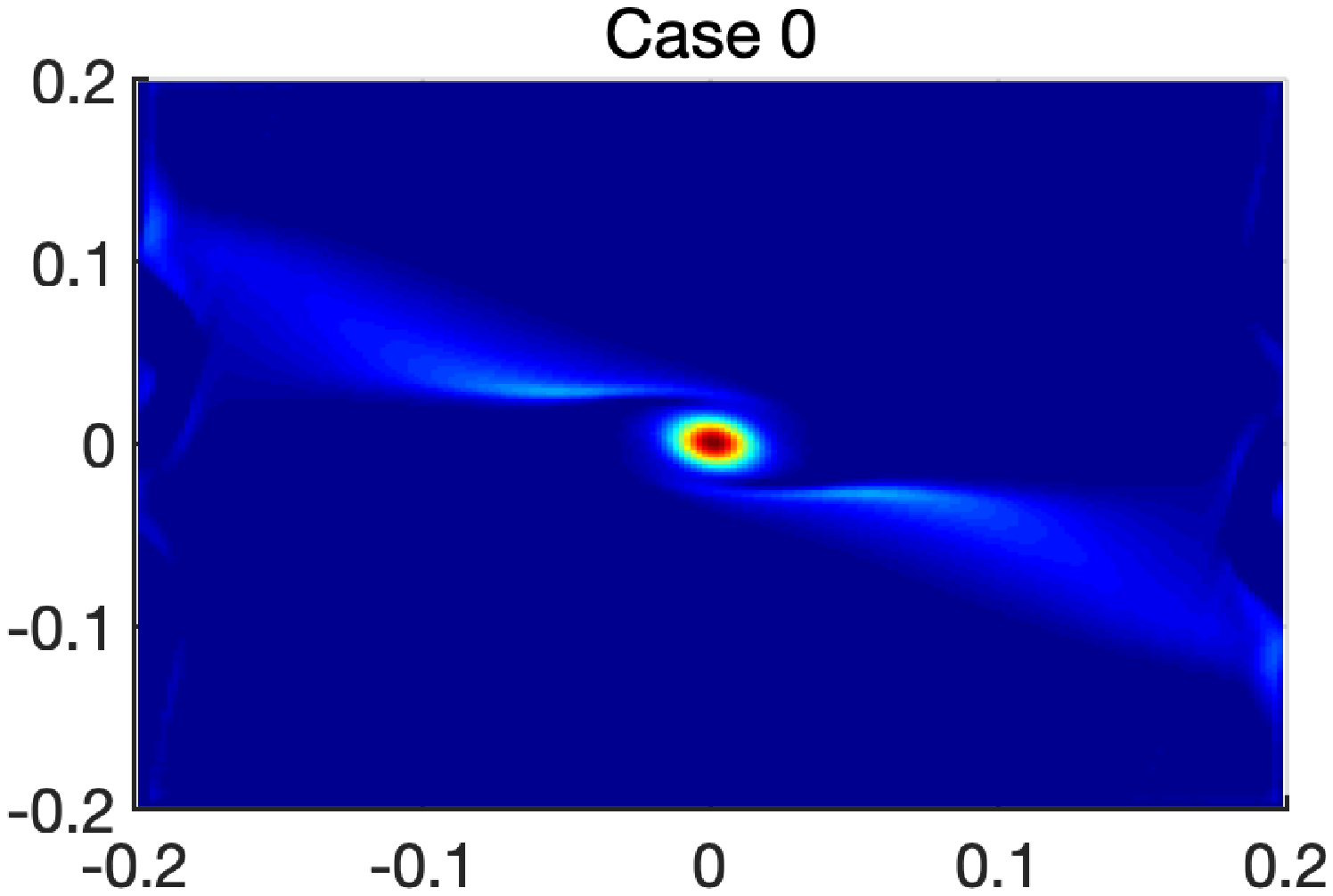}
\includegraphics[width=.3\textwidth]{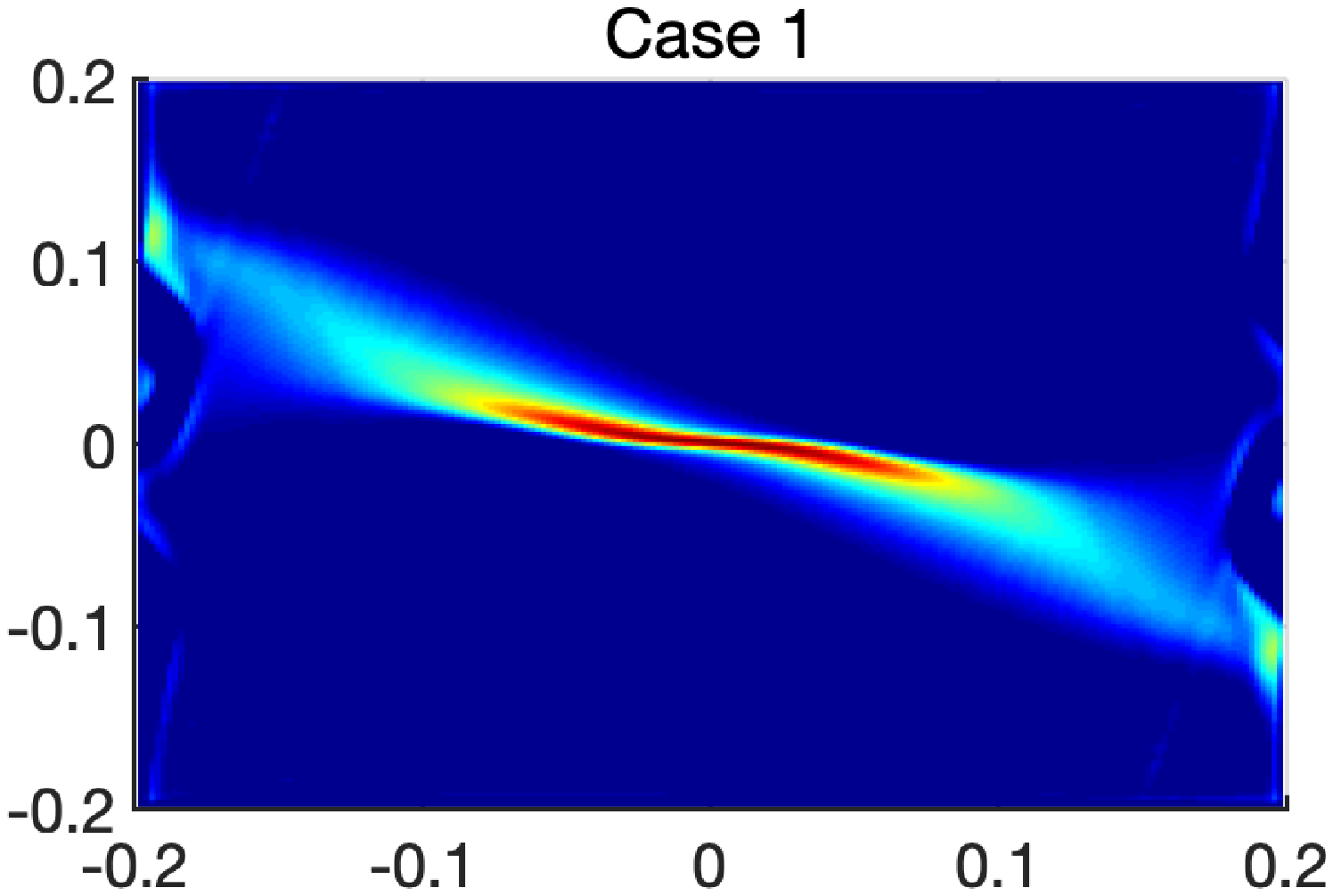}
\includegraphics[width=.3\textwidth]{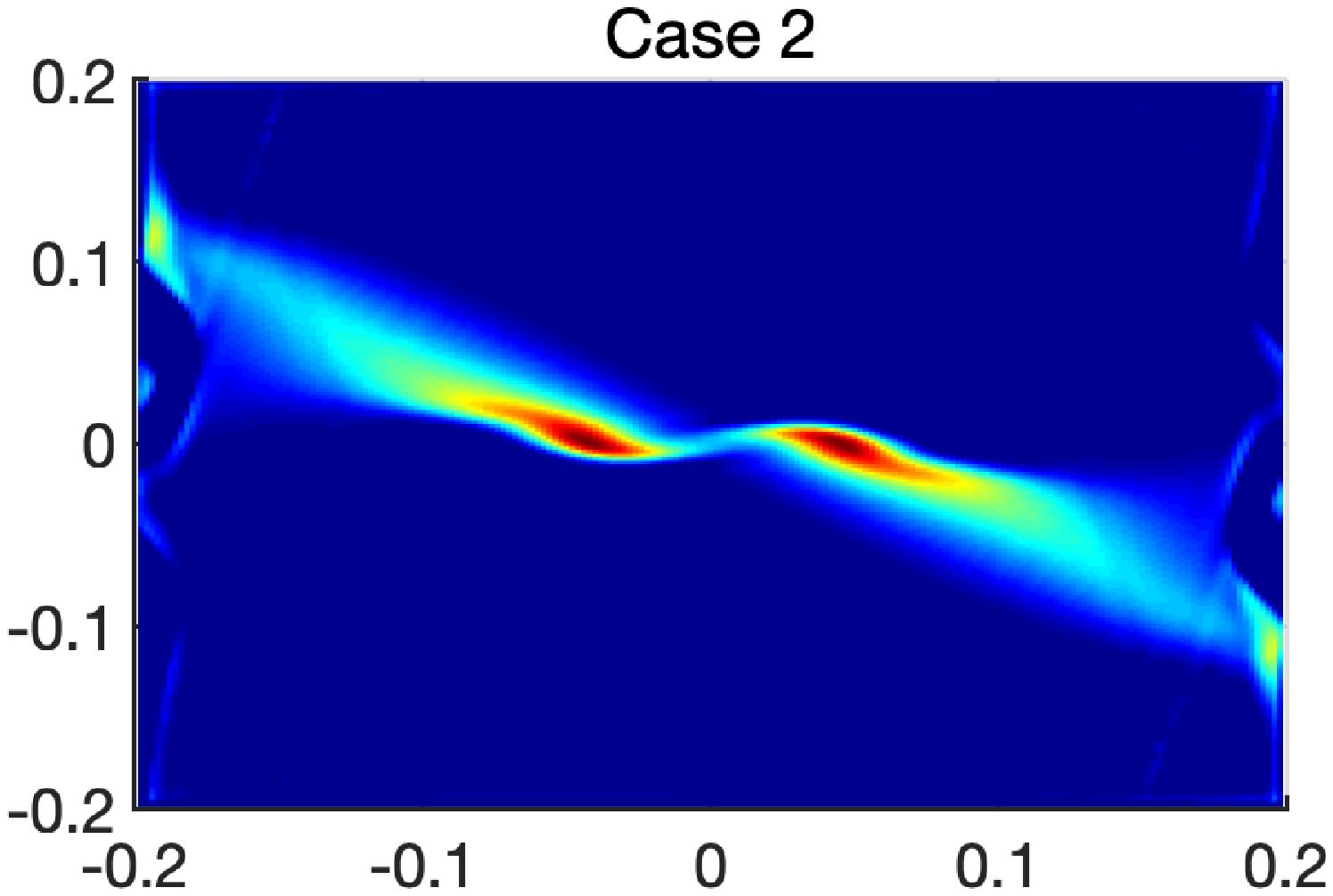}
\caption{Example 6: vorticity profiles at $T=1$.  $\beta = 0.1$, $\alpha =0.75$ and $\theta _0 = \frac{\pi}{8}$. }
\label{EX71}
\end{figure}

\begin{figure}[!htbp]
\centering
\includegraphics[width=.3\textwidth]{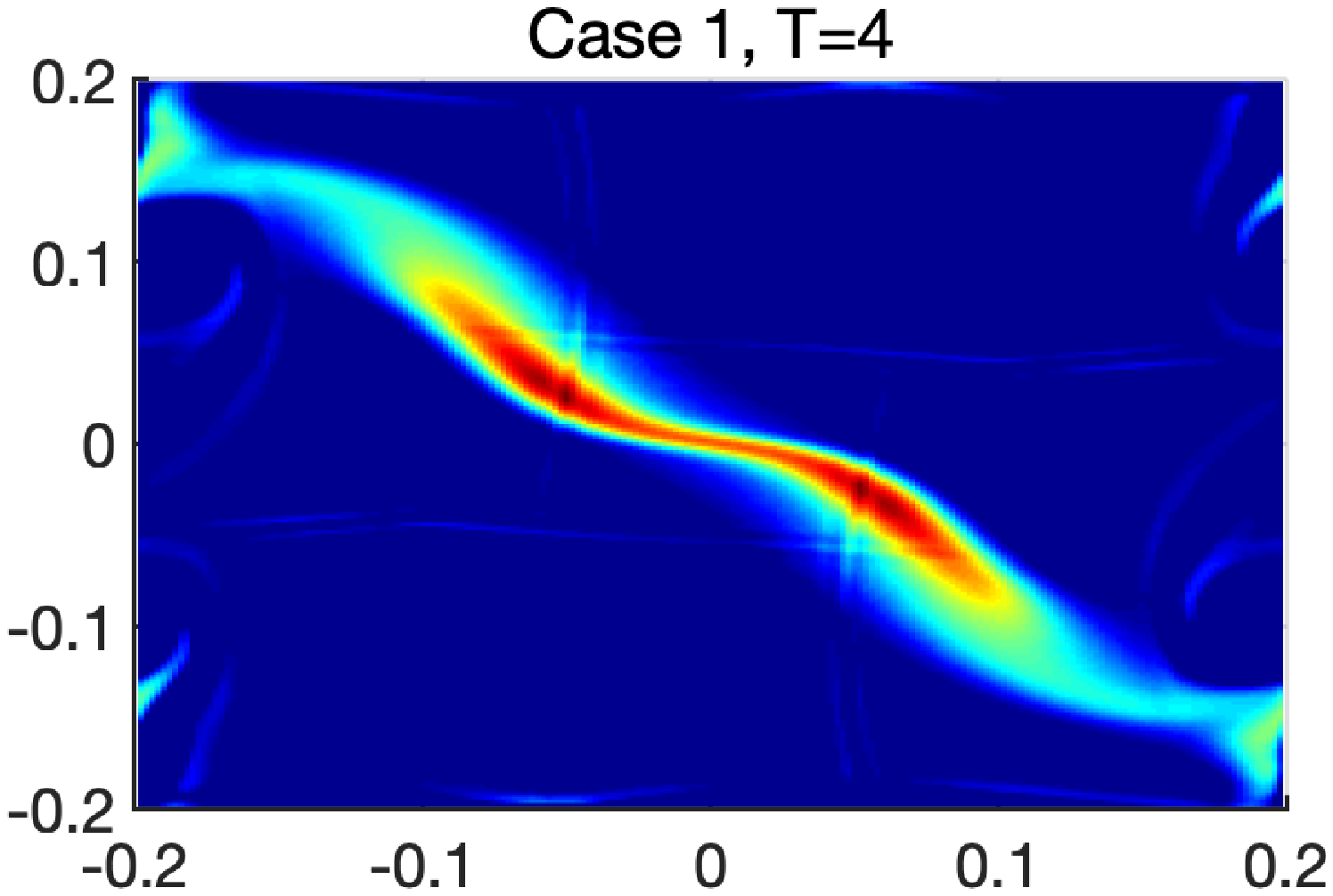}
\includegraphics[width=.3\textwidth]{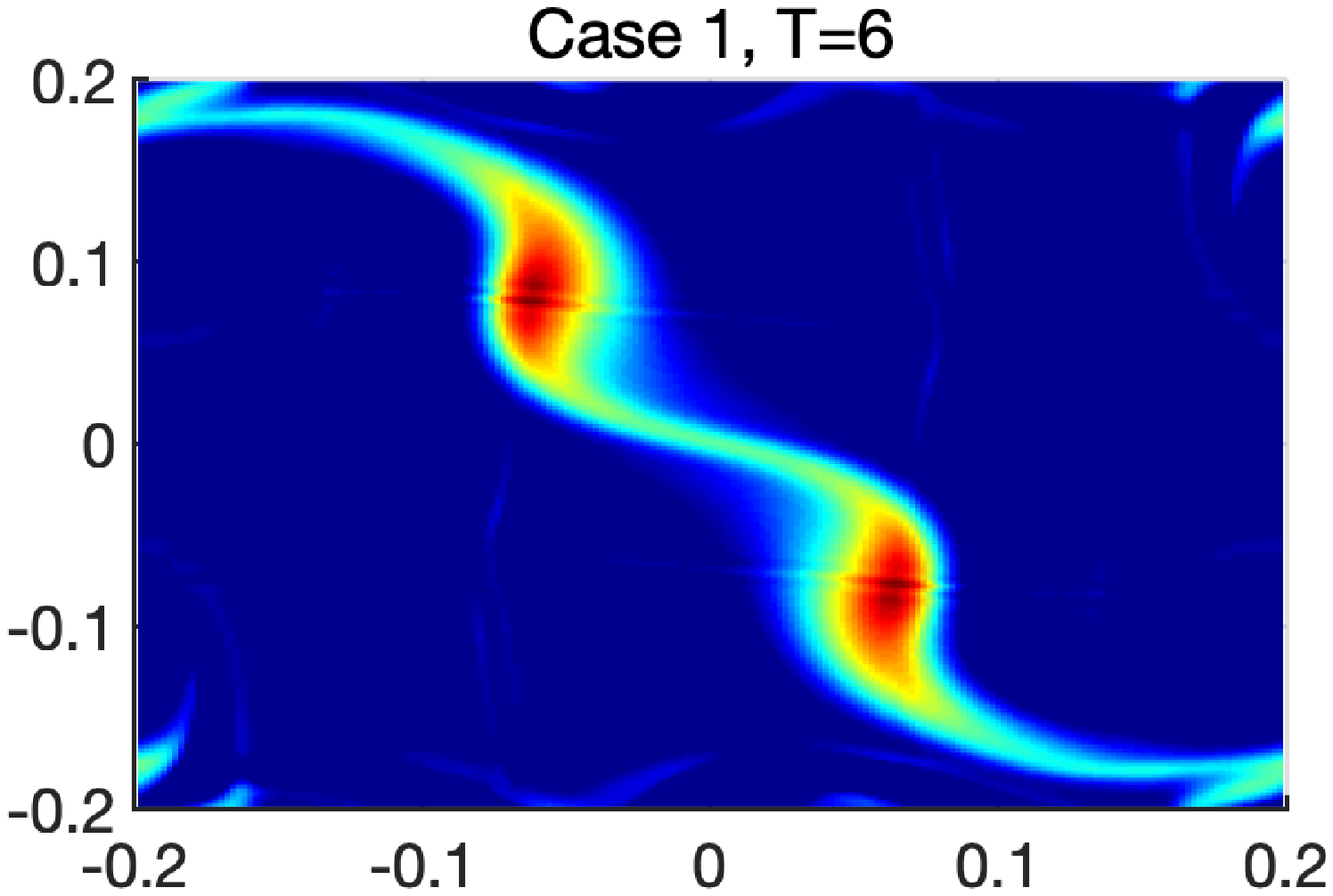}
\includegraphics[width=.3\textwidth]{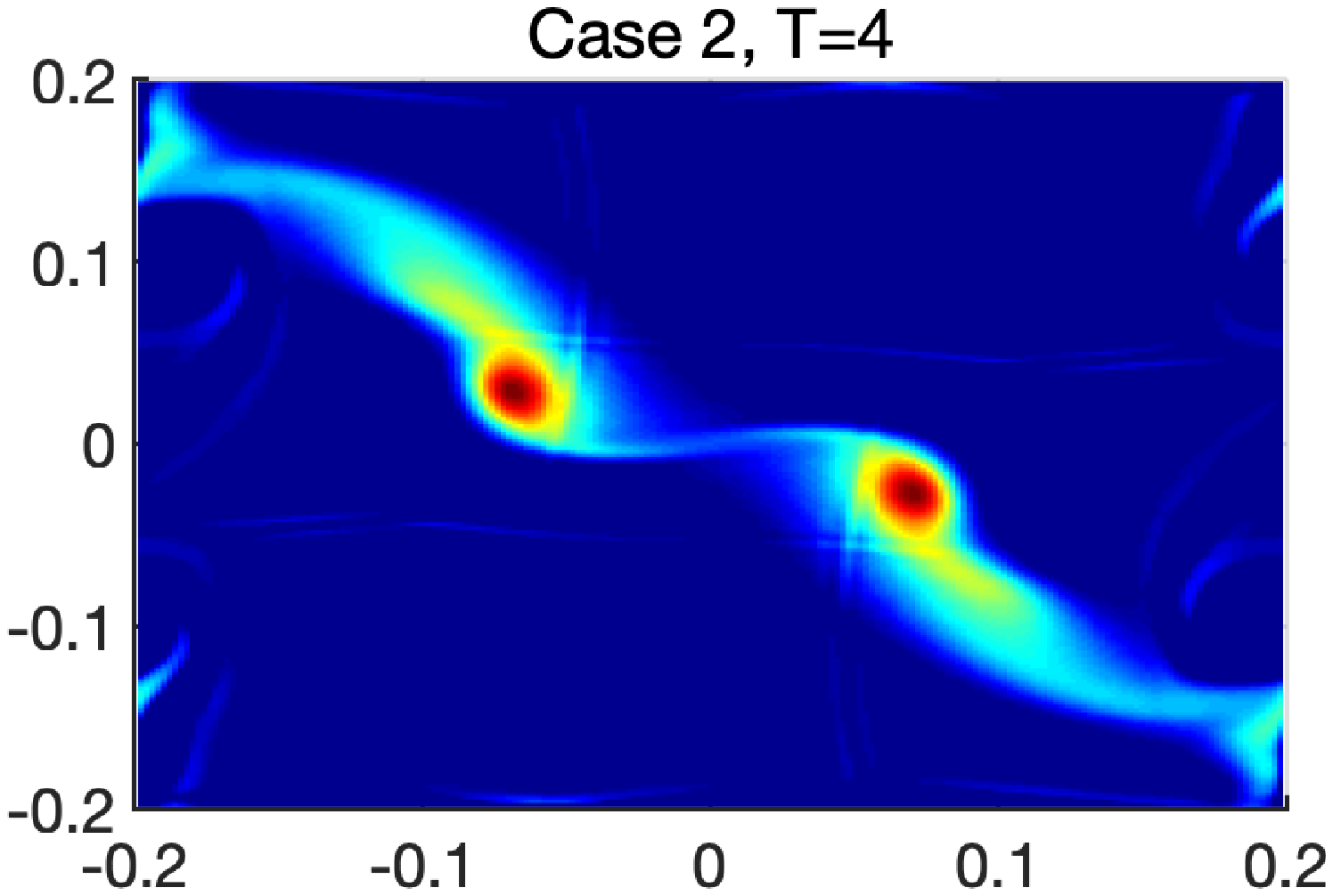}
\caption{Example 6: vorticity profiles at larger times. $\beta = 0.1$, $\alpha =0.75$ and $\theta _0 = \frac{\pi}{8}$. }
\label{EX72}
\end{figure}

\noindent \textit{Example 6.} $\alpha = 0.75$. In this example, the difference between solutions show up at $T=1$. Confirmed with a larger time test, we can see that both Case 1 and Case 2 eventually form two spirals while Case 0 results in one single spiral.  See Figure \ref{EX71} and  \ref{EX72}.\\

Note that in the section \ref{SubsecBeta}, the non-unique solutions have been observed for $\alpha =0.95$, so in the following we test more examples with larger $\alpha$ values.\\

\begin{figure}[!htbp]
\centering
\includegraphics[width=.3\textwidth]{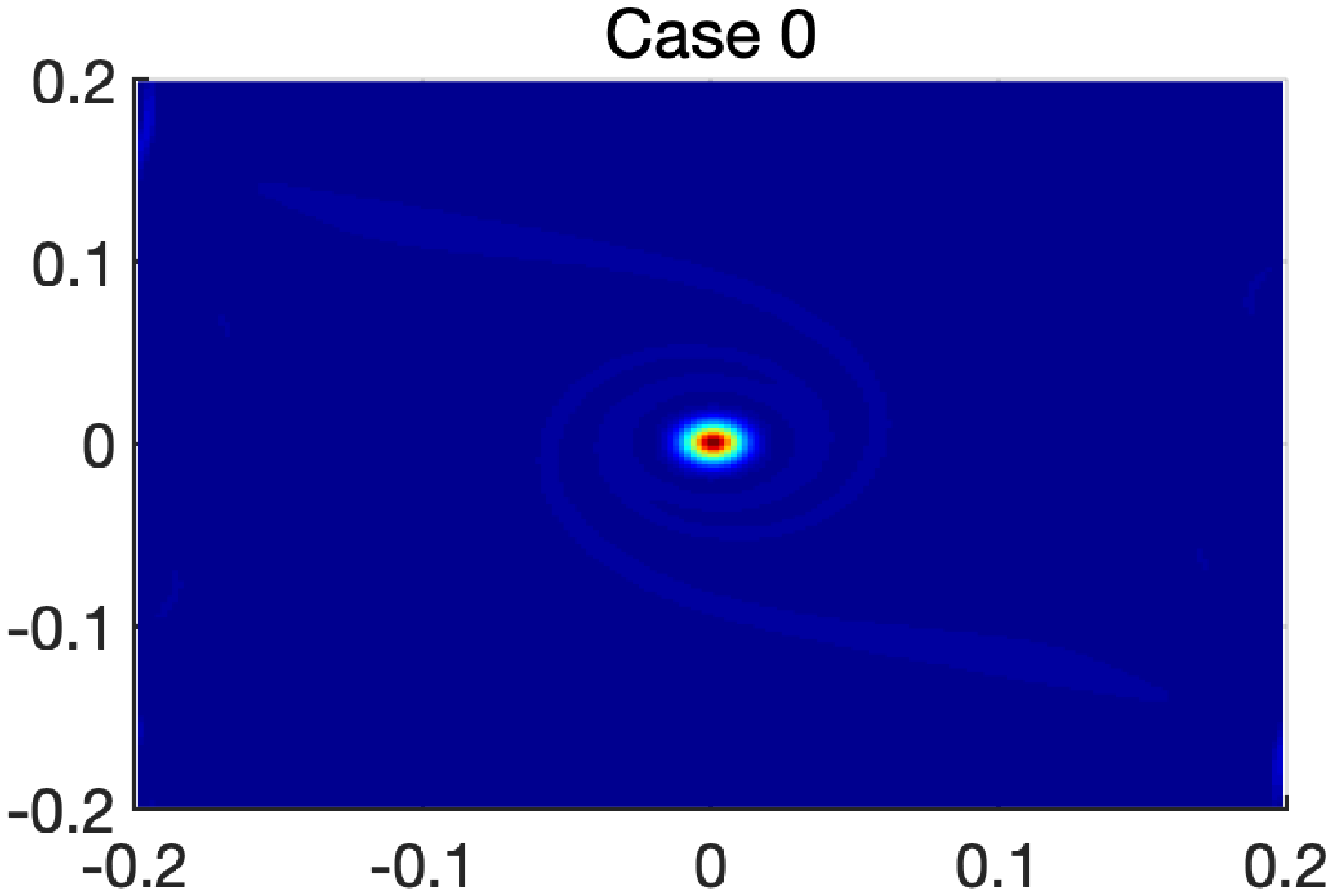} 
\includegraphics[width=.3\textwidth]{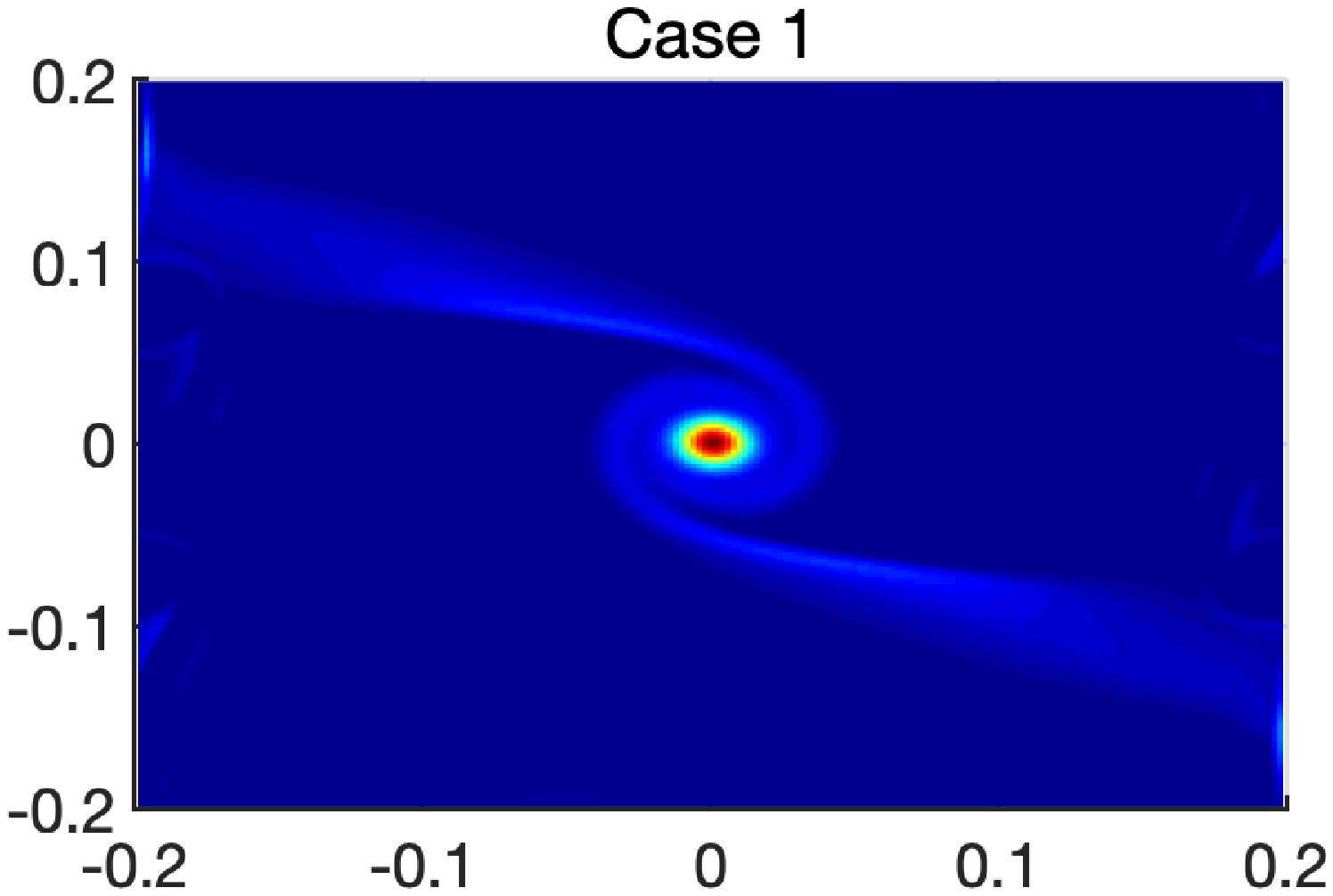}
\includegraphics[width=.3\textwidth]{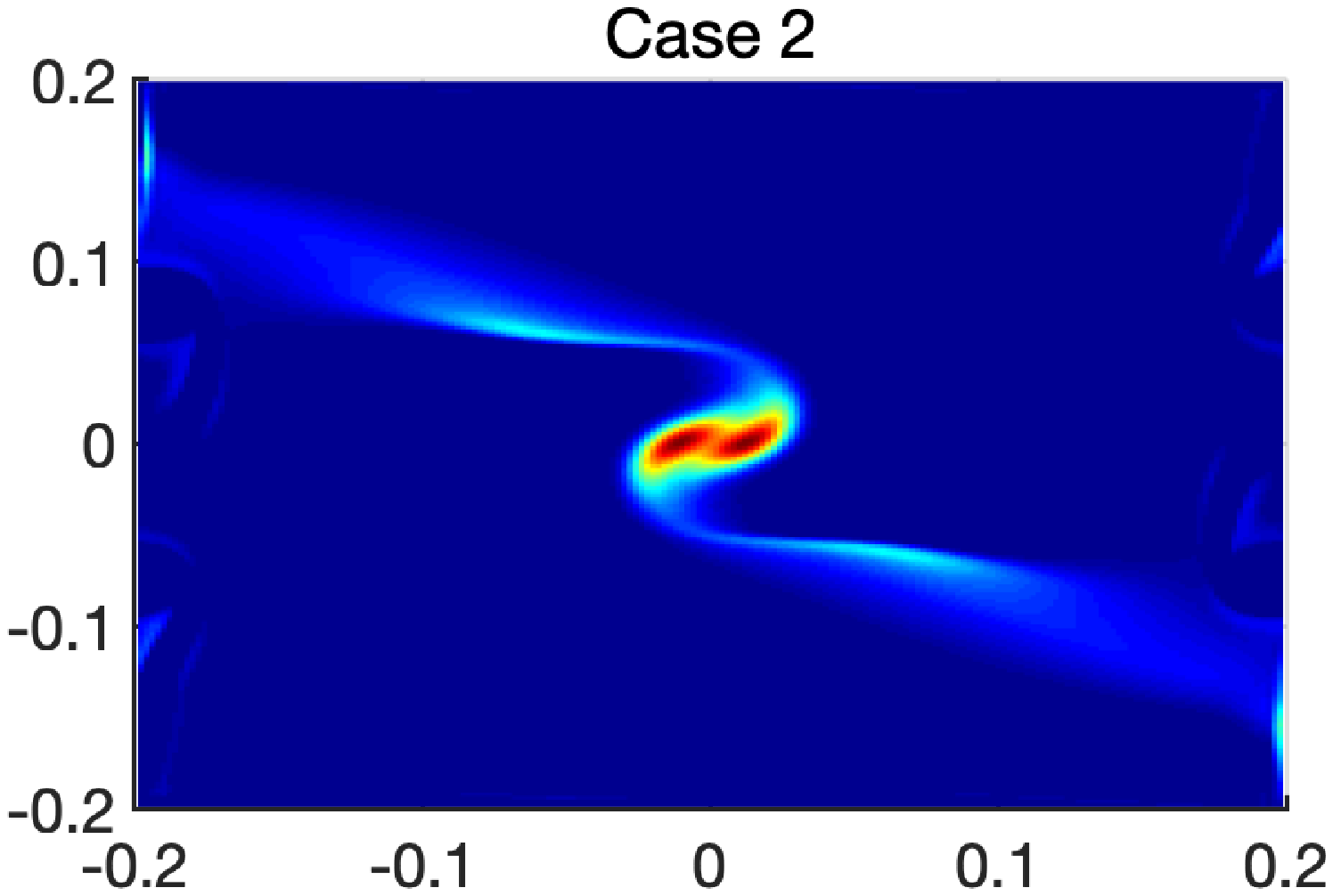}
\caption{Example 7: vorticity profiles at $T=1$. $\beta = 0.1$, $\alpha =1.2$ and $\theta _0 = \frac{\pi}{8}$.}
\label{EX9}
\end{figure}

\noindent \textit{Example 7.} $\alpha = 1.2$. All three cases lead to one single spiral. See Figure \ref{EX9}. \\

\begin{figure}[!htbp]
\centering
\includegraphics[width=.3\textwidth]{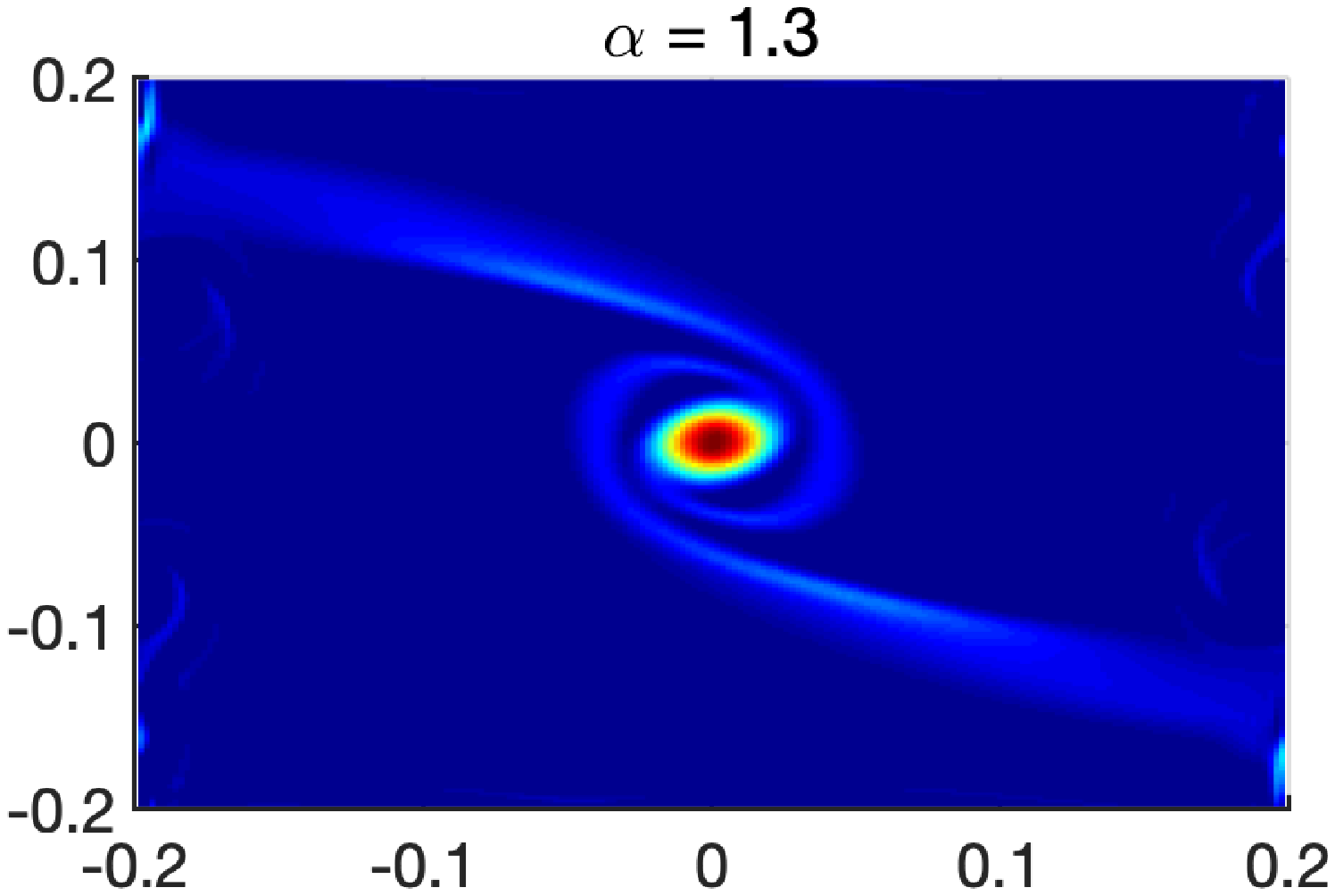}
\includegraphics[width=.3\textwidth]{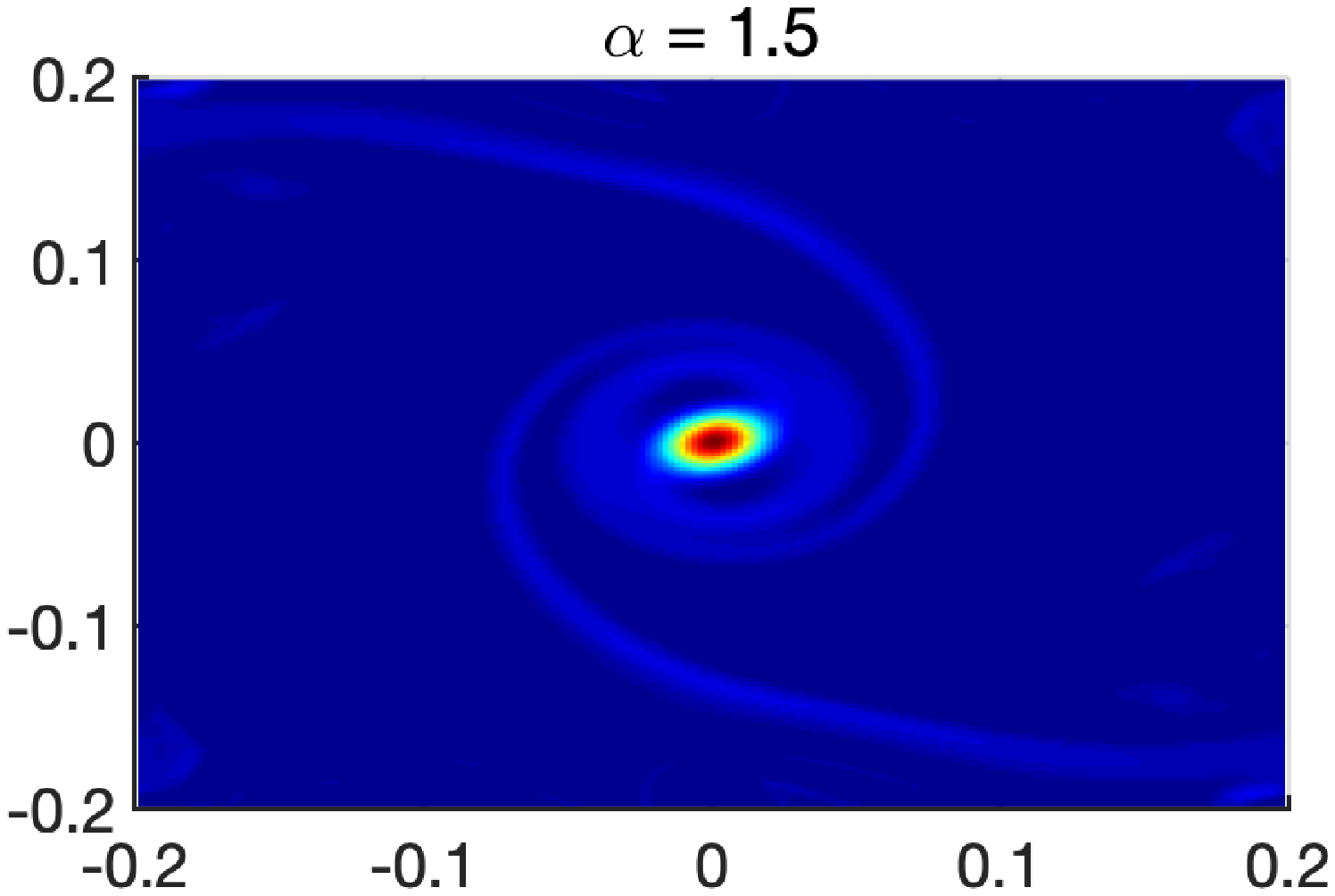}
\includegraphics[width=.3\textwidth]{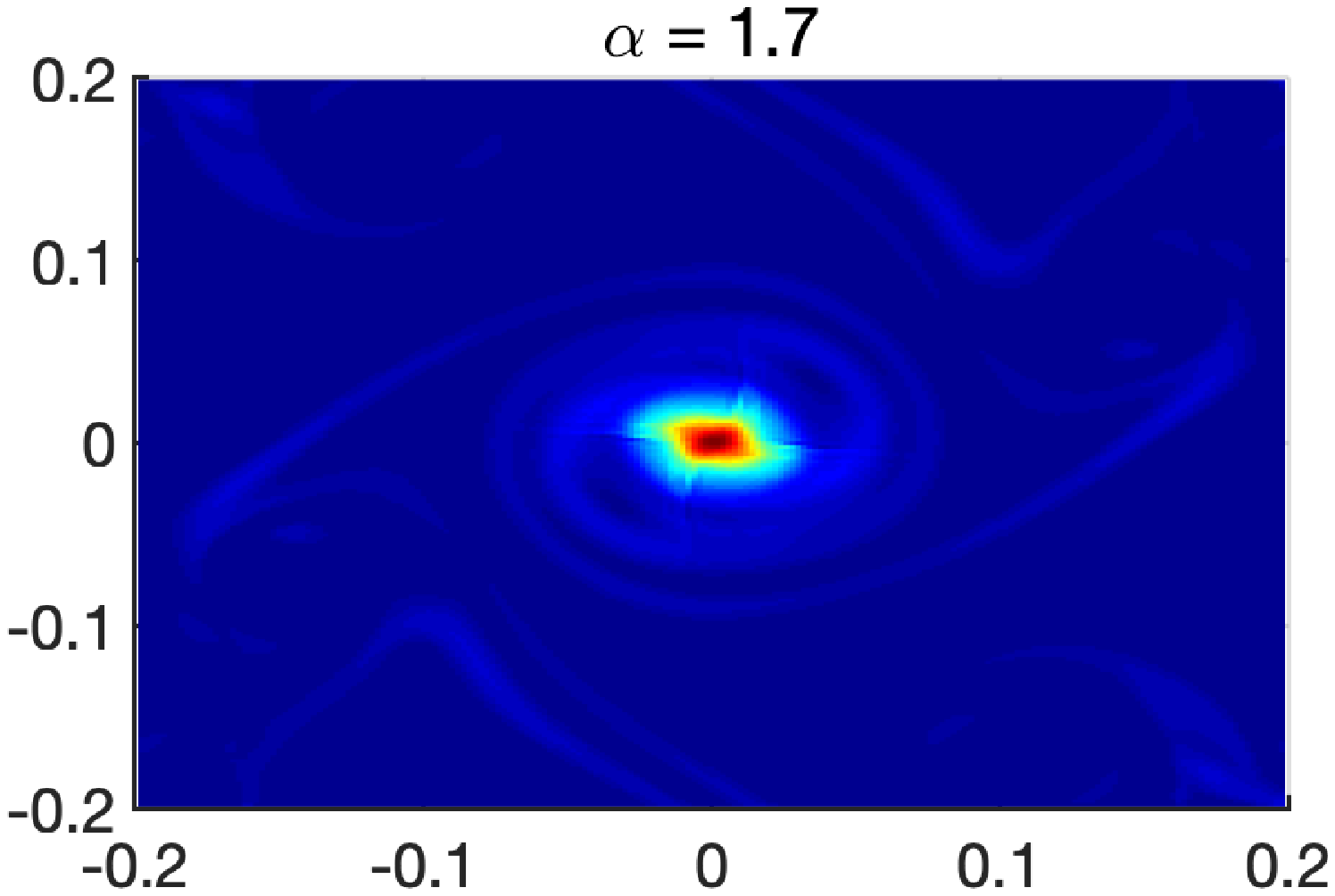}
\caption{Example 8:  Case 2 vorticity profiles at $T=1$ with some large $\alpha$ values given $\beta = 0.1$ and $\theta _0 = \frac{\pi}{8}$. }
\label{EX10}
\end{figure}

\noindent \textit{Example 8.} We test Case 2 with $\alpha =1.3$, 1.5 and 1.7 respectively. We observe that the results are all single spirals. Note that Case 0 and Case 1 always form single spirals when Case 2 does, so the uniqueness of solutions is preserved with the chosen $\alpha$ values. See Figure \ref{EX10}.\\

In summary, we conclude that with a fixed (carefully chosen) value of $\beta$ and $\theta_0$, there exist $\alpha_i\in (0,2)$, $i=1,2,3$ such that 
\begin{itemize}
\item[(1)] when $0<\alpha <\alpha _1$,  the vorticities generated from the three types of proposed initial data are similar;
\item[(2)] as $\alpha$ increases but does not exceed  $\alpha _2$, Case 0 begins to form a single spiral while the other two cases form two spirals. Non-uniqueness of the solution is indicated; 
\item[(3)] as $\alpha$ continues to increase but does not  exceed  $\alpha _3$, Case 1 begins to behave like Case 0, which results in a single spiral, while Case 2 still forms two spirals. Again, non-uniqueness of the solution is indicated;
\item[(4)] when $\alpha _3<\alpha<2$, all three cases eventually lead to vorticities as one single spiral;
\end{itemize} 

\subsubsection{Effects of $\theta _0$}

In this section, we test how the choices of $\theta _0$, the support of the function $\phi(\theta)$ in the initial vorticity, affects the solution. We fix $\beta =0.1$ and $\alpha  = 0.95$ and test only Case 0 and Case 2 initial data. \\

\begin{figure}[!htbp]
\centering
\includegraphics[width=.3\textwidth]{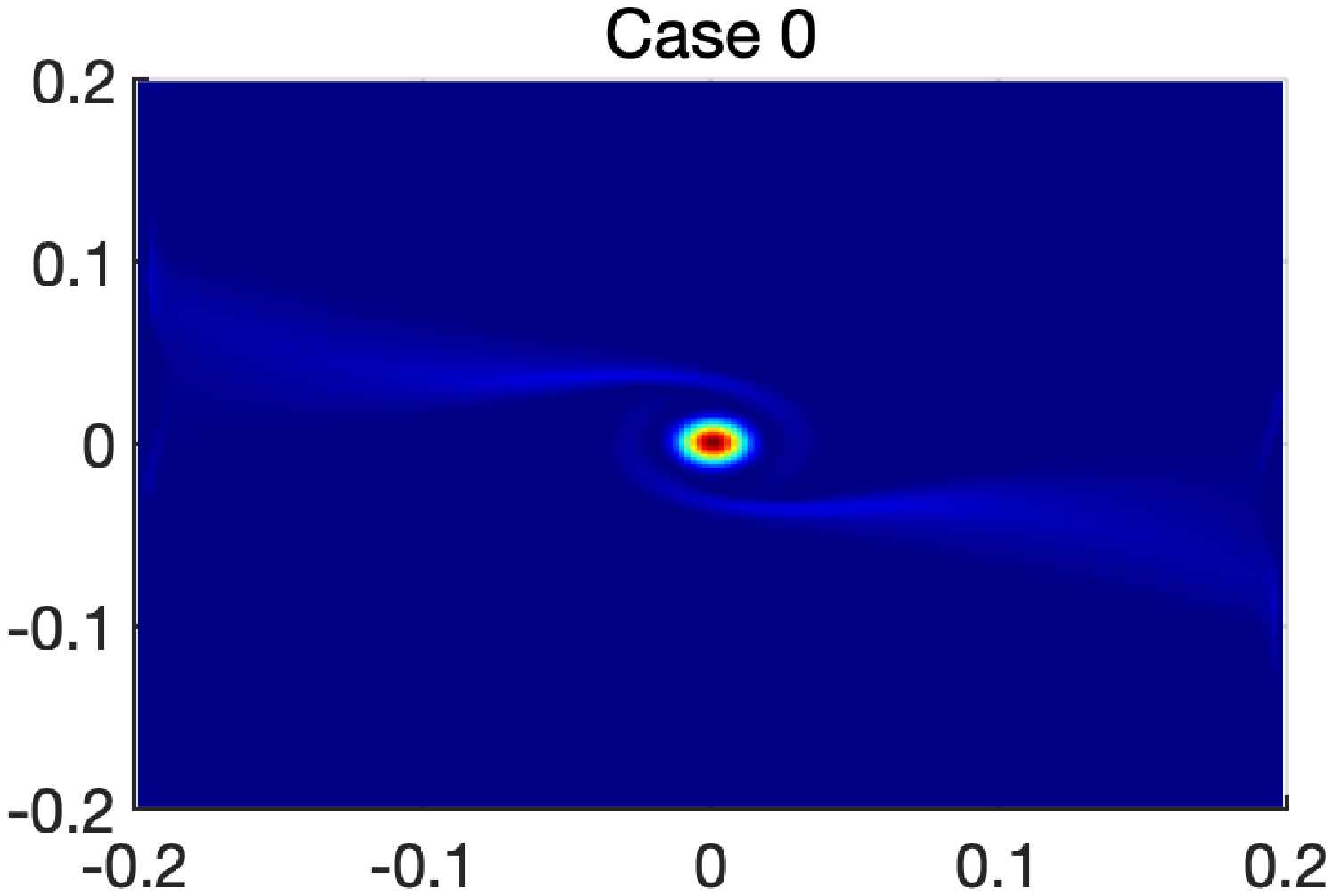}
\includegraphics[width=.3\textwidth]{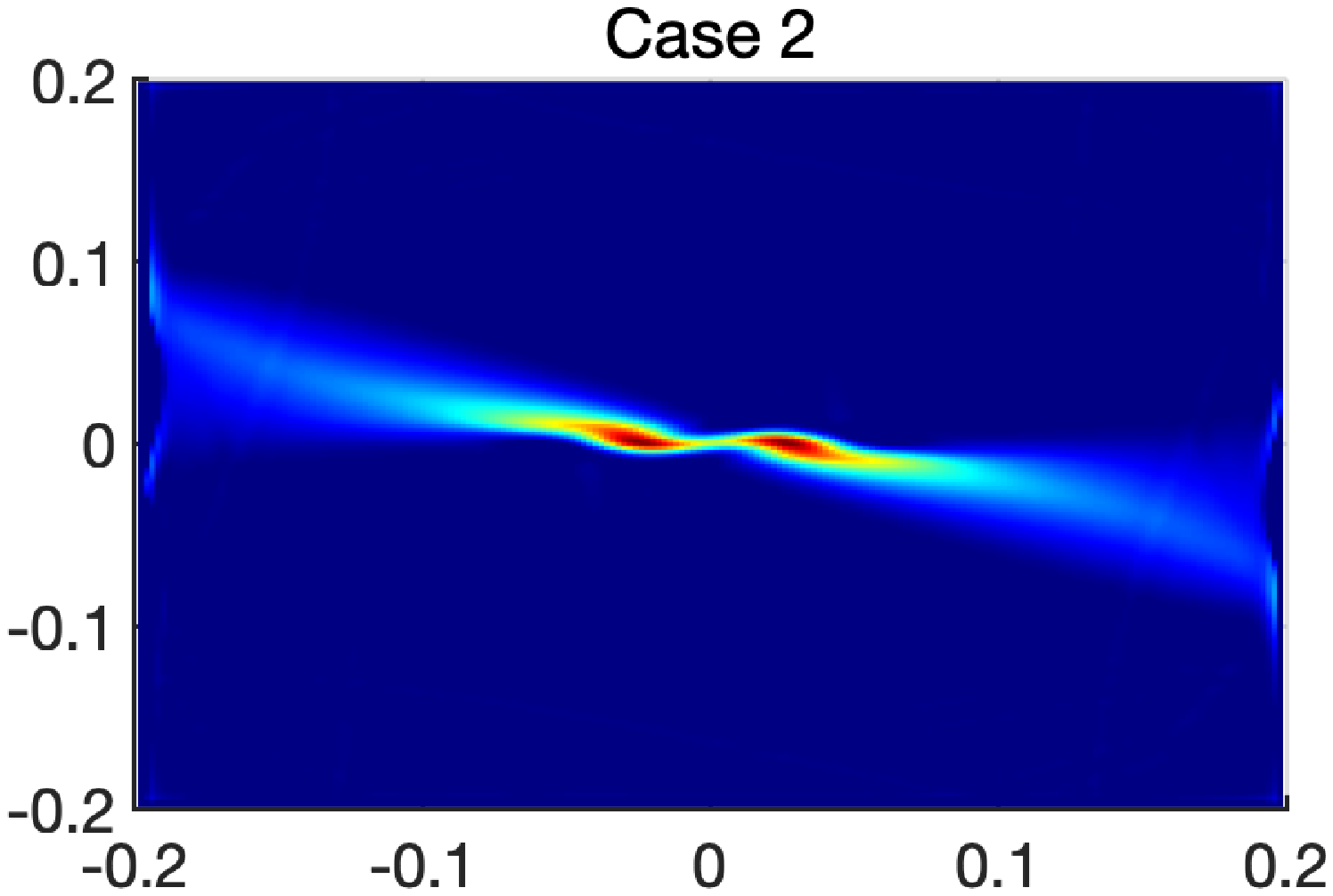}
\caption{Example 9: vorticity profiles at $T=1$. $\beta = 0.1$, $\alpha =0.95$ and $\theta _0 = \frac{\pi}{10}$.}
\label{EX11}
\end{figure}

\begin{figure}[!htbp]
\centering
\includegraphics[width=.3\textwidth]{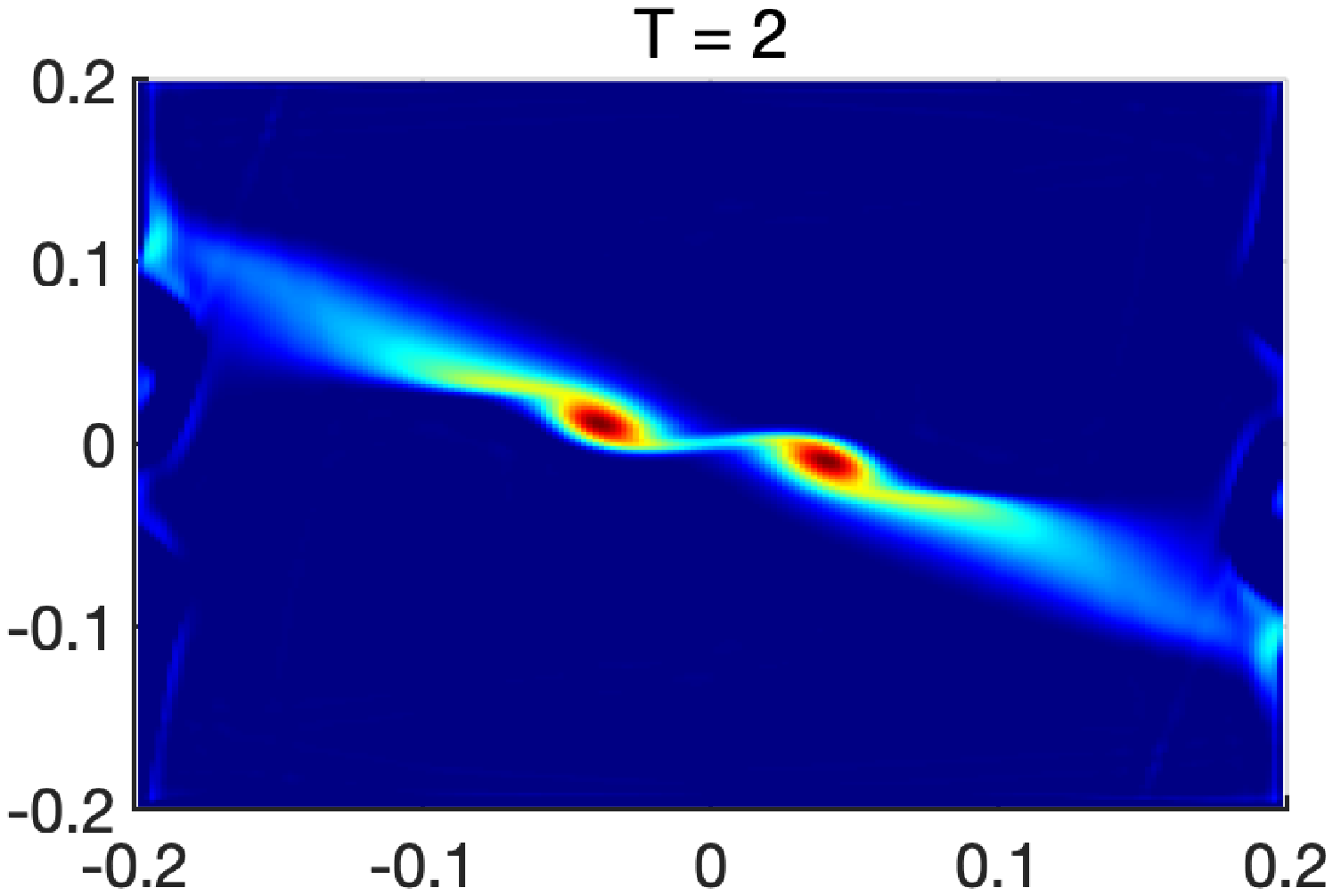}
\includegraphics[width=.3\textwidth]{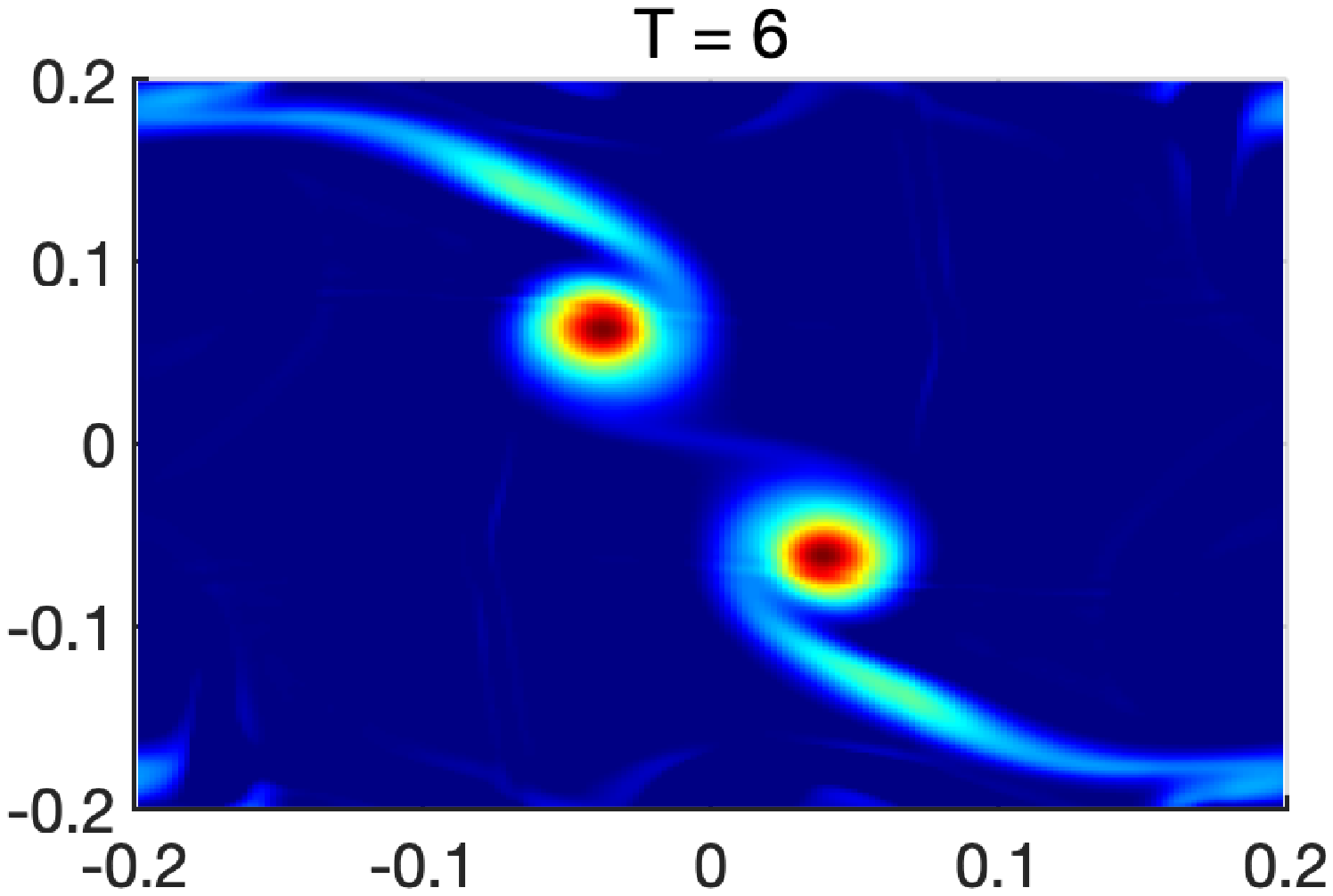}
\caption{Example 9: Case 2 vorticity profiles at larger times. $\beta = 0.1$, $\alpha =0.95$ and $\theta _0 = \frac{\pi}{10}$.}
\label{EX112}
\end{figure}

\noindent \textit{Example 9.} $\theta _0=\frac{\pi}{10}$. Non-unique vorticity profiles are observed. 
See Figure \ref{EX11} and \ref{EX112}. \\

\begin{figure}[!htbp]
\centering
\includegraphics[width=.3\textwidth]{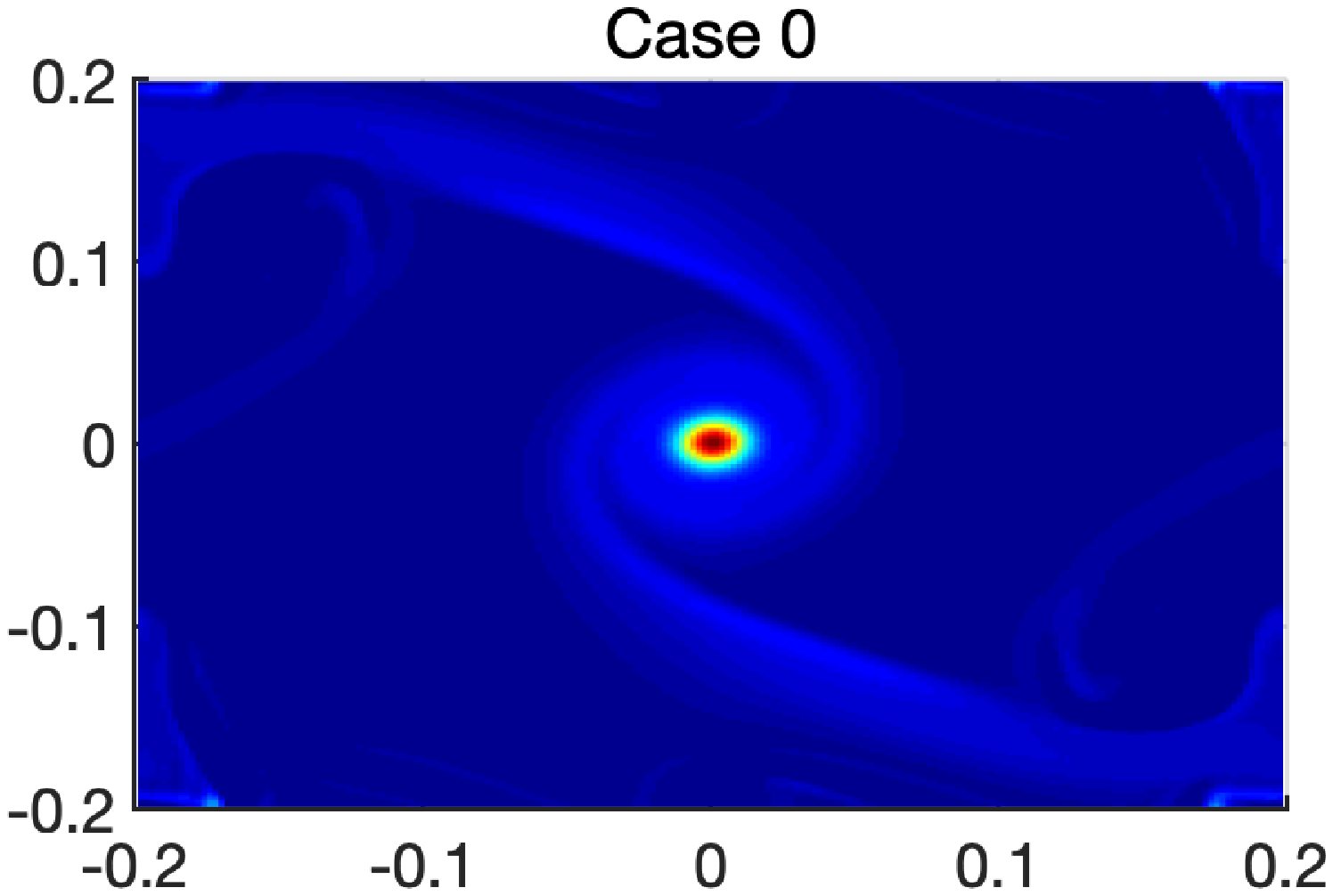} 
\includegraphics[width=.3\textwidth]{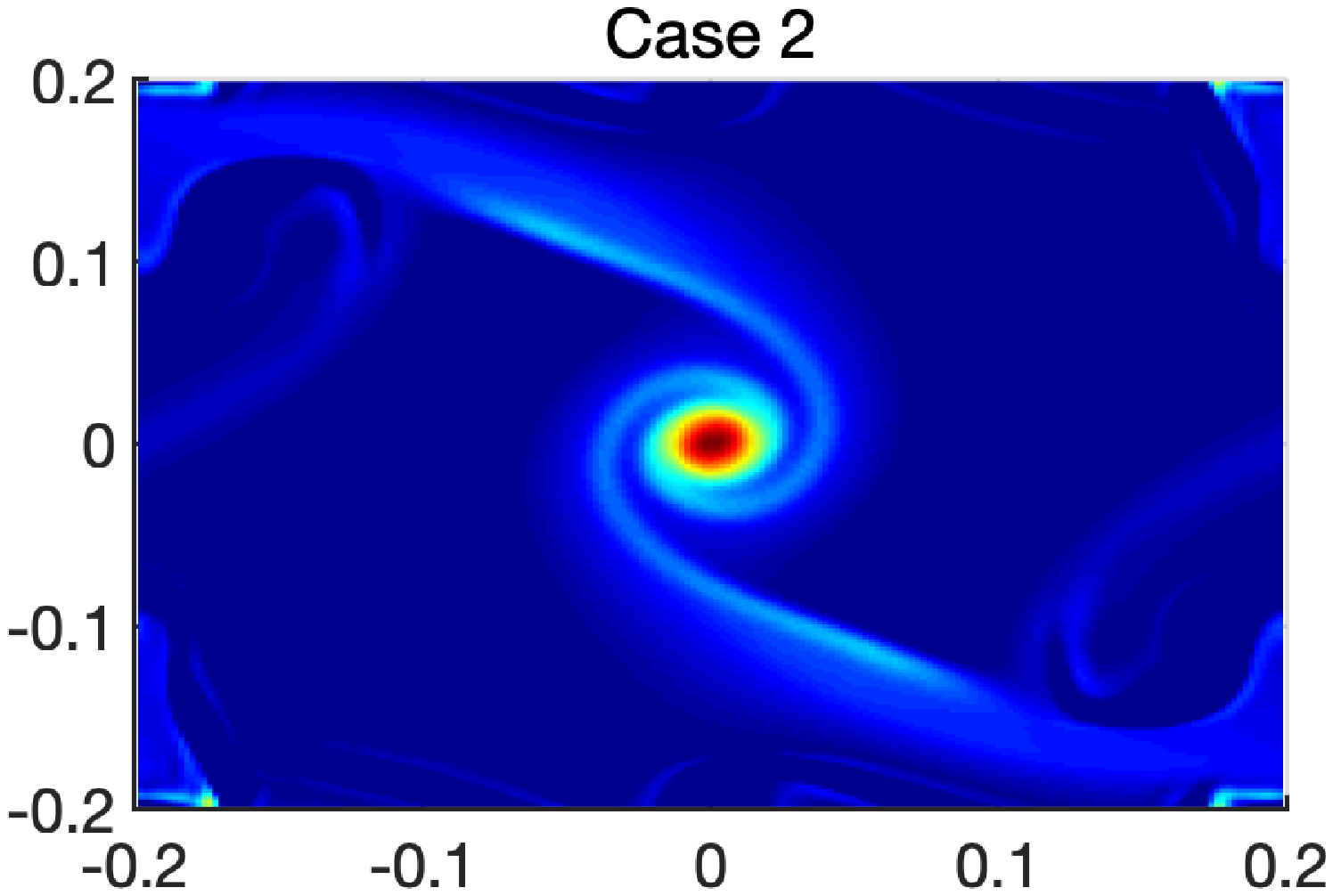}
\caption{Example 10: vorticity profiles at $T=1$. $\beta = 0.1$, $\alpha =0.95$ and $\theta _0 = \frac{\pi}{4}$.}
\label{EX12}
\end{figure}

\noindent \textit{Example 10.} $\theta _0=\frac{\pi}{4}$. Both cases result in one single spiral. 
See Figure \ref{EX12}.
\\

\begin{figure}[!htbp]
\centering
\includegraphics[width=.3\textwidth]{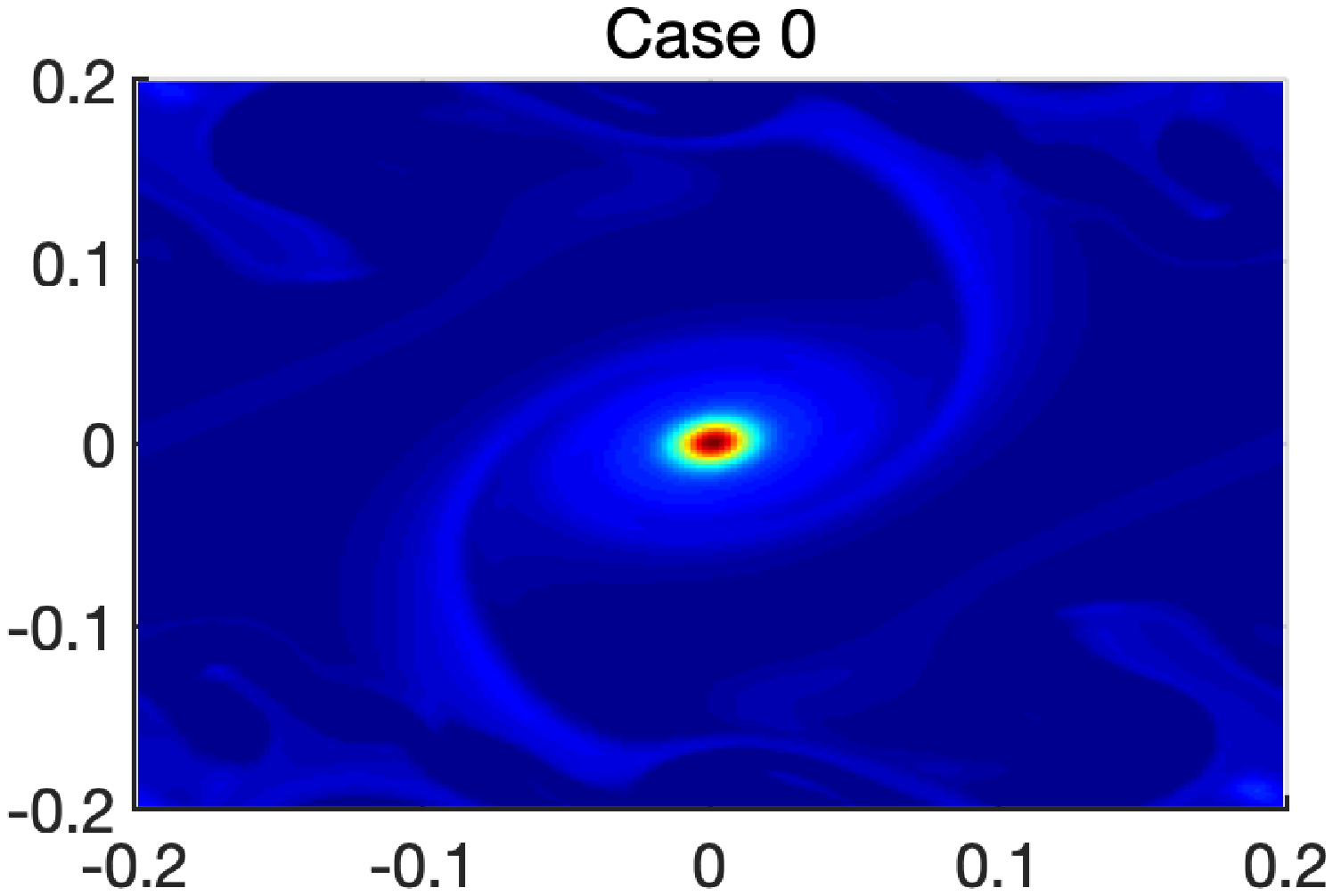} 
\includegraphics[width=.3\textwidth]{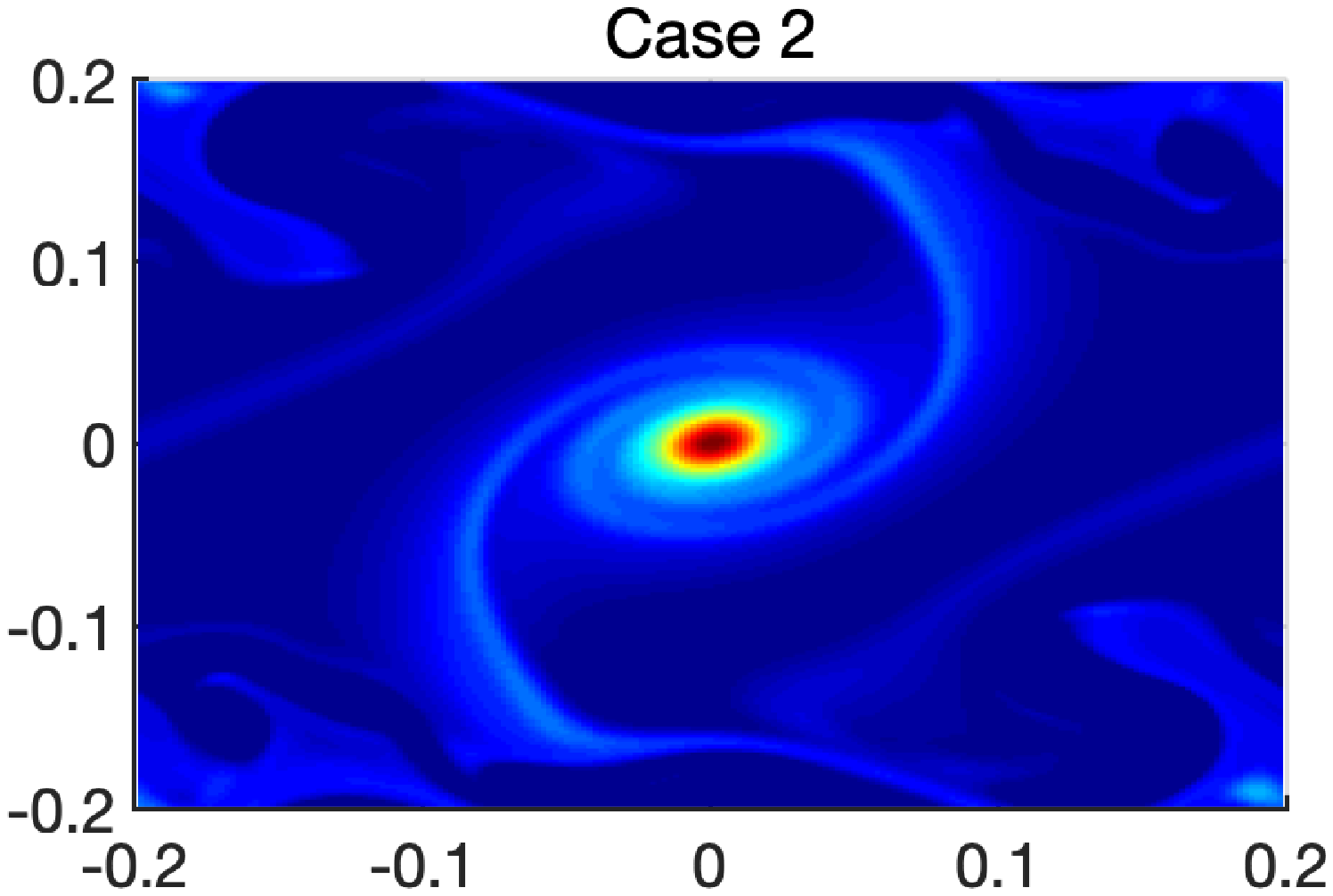}
\caption{Example 11: vorticity profiles at $T=1$. $\beta = 0.1$, $\alpha =0.95$ and $\theta _0 = \frac{\pi}{3}$.}
\label{EX13}
\end{figure}

\noindent \textit{Example 11.} $\theta _0=\frac{\pi}{3}$. Both cases result in one single spiral. 
See Figure \ref{EX13}.\\

In summary, we conclude that with a fixed (carefully chosen) value of $\beta$ and $\alpha$, there exists a positive $\theta ^*_0$ such that when $\theta _0<\theta ^*_0$, distinct vorticity profiles are observed and the non-uniqueness of (\ref{eq:2DEuler}) is indicated, while when $\theta_0  \geq \theta _0^*$ one always obtains the vorticity as a single spiral.

\section{Effects of the compressibility}
In this section, we present further numerical results to examine how the compressibility of the gas changes the solution structure. 

From a Physics view-point, the fluid should behave (asymptotically) like an incompressible one when the density is almost constant, the velocity is small and we look at large time scales.  It is known in \cite{KM82} that the rescaling of $\rho$ and $u$ (and thus $p$) via 
$$
t \to  \lambda t, \quad u \to \lambda u 
$$
will still lead to (\ref{eq:2DEuler}) with $p$ replaced by 
$$
p=\frac{A}{\lambda^2} \rho^{\gamma}.
$$
Here $\lambda$ is essentially linked to the Mach number, $M=|\bar v|(dp(\bar \rho)/d\rho)^{-1/2}$, the ratio of fluid speed to sound speed, where $\bar \rho$ is the mean density, and $\lambda = M\sqrt{A\gamma}$ upon a detailed non-dimensional scaling \cite{KM82}. Hence the moment equation indicates that $\rho$ should be like $\bar \rho +O(M^2)$ for $M$ small. For $\bar \rho=1$, one may pass to the limit $M\to 0$ to obtain  
\begin{align*}
& u_x +v_y = 0,\\
& u_t + (\rho u^2 + P)_x + (\rho uv)_y =0,\\
&v_t + (\rho uv)_x + (\rho v^2 + P)_y =0,
\end{align*}
where $P$ is the `limit' of $A(\rho^\gamma -1)/\lambda^2$.   In other words, we recover the  incompressible Euler equations, and the hydrostatic pressure appears as the limit of the “renormalized”  thermodynamical pressure. Rigorous justification of this limit can be found in \cite{KM82}.  
 
In Example 3, with parameters  $\beta = 0$, $\alpha = 0.95$, $\theta = \frac{\pi}{8}$ and $A=1$, the non-uniqueness of solutions with these parameter values have been observed, where Case 0 initial condition result in one single-spiral shape vorticity while Case 2 results in two spirals. In the following, we fix these parameters in the initial condition, and test with different values of $A$ for (\ref{eq:2DEuler}) with the pressure function $p=A\rho^{1.4}$ for Case 2 only.  
 
\begin{figure}[!htbp]
\centering
\includegraphics[width=.3\textwidth]{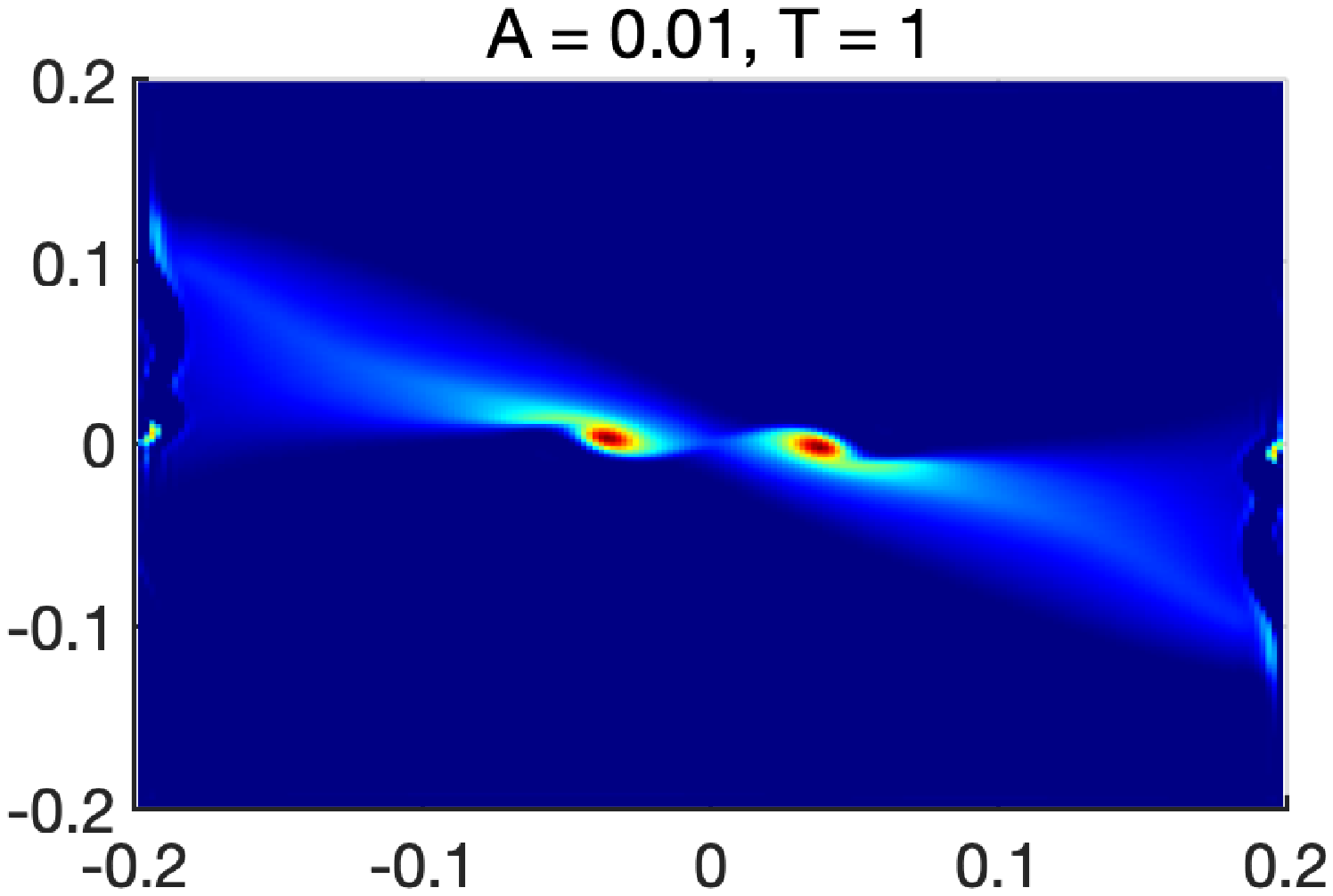}
\includegraphics[width=.3\textwidth]{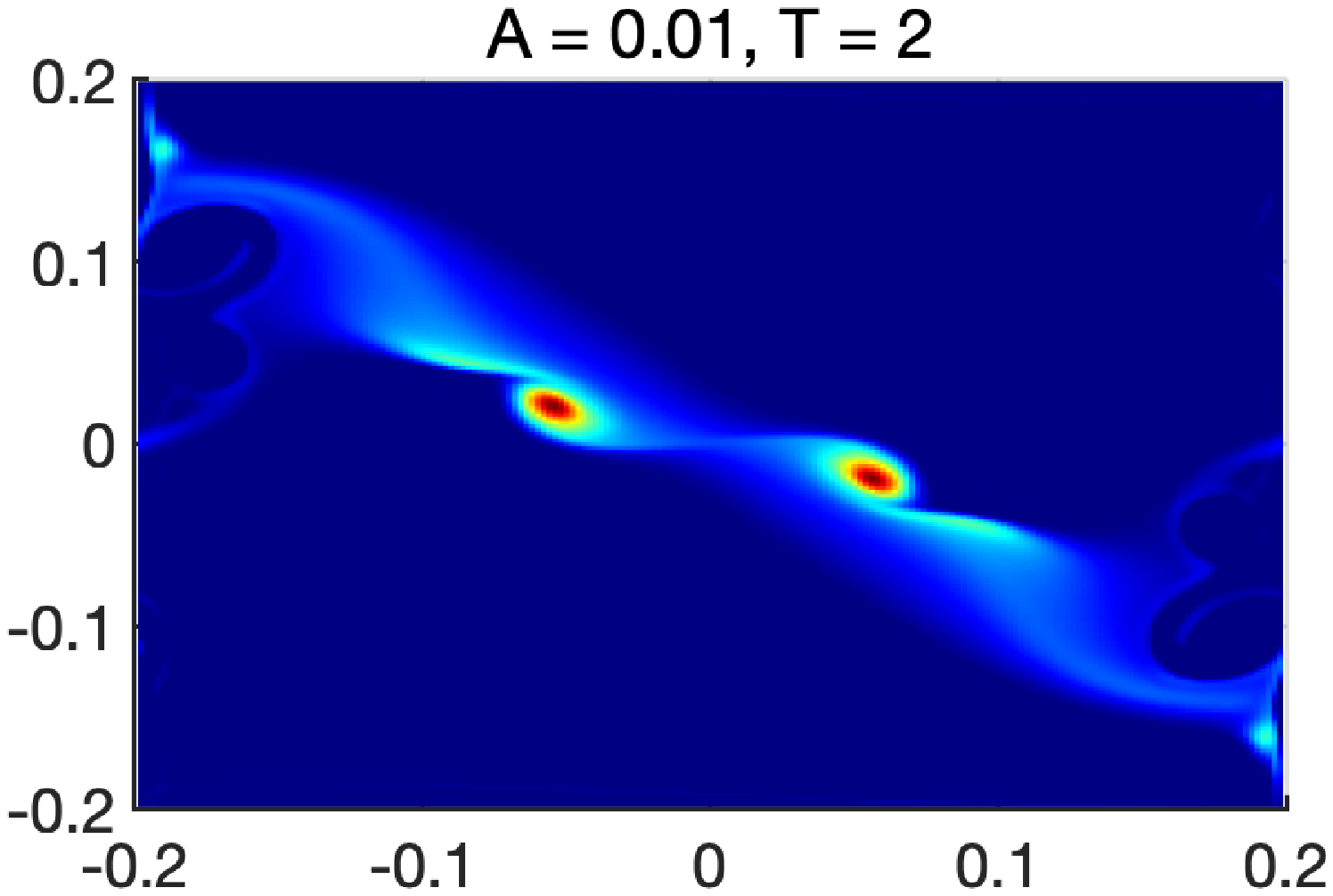}
\includegraphics[width=.3\textwidth]{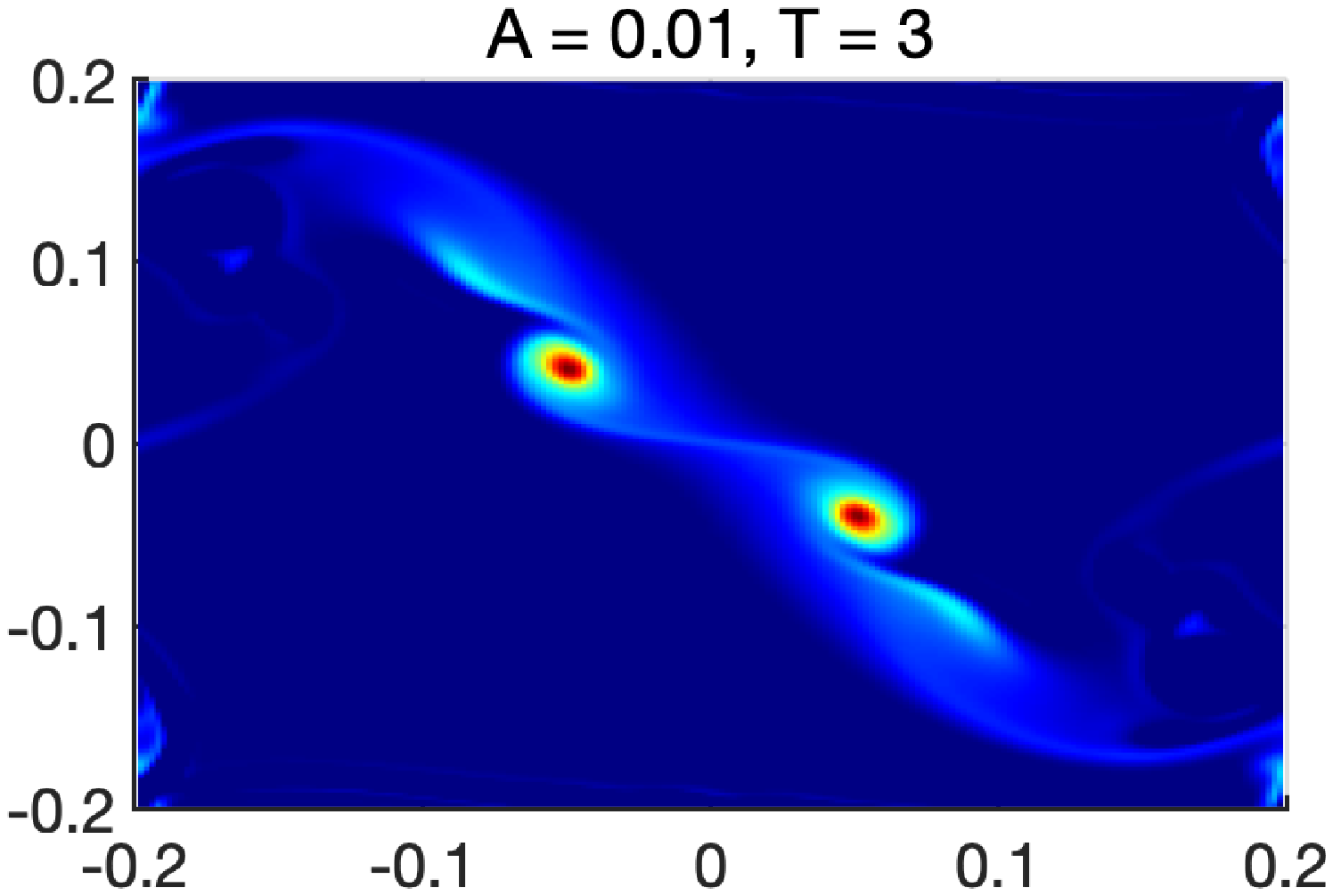}
\caption{Case 2 vorticity profiles. Highly compressible system.}
\label{A001}
\end{figure}

\begin{figure}[!htbp]
\centering
\includegraphics[width=.3\textwidth]{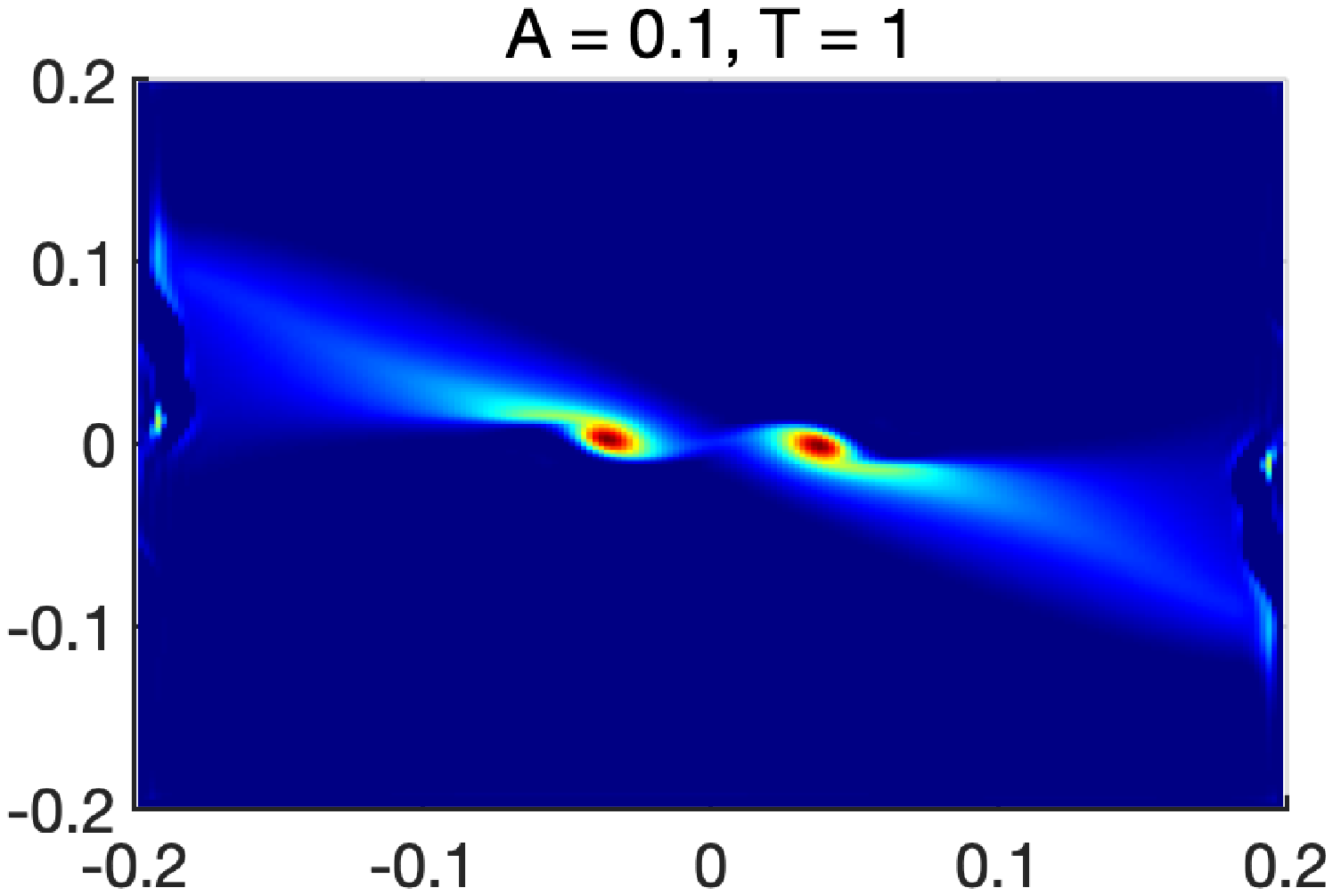}
\includegraphics[width=.3\textwidth]{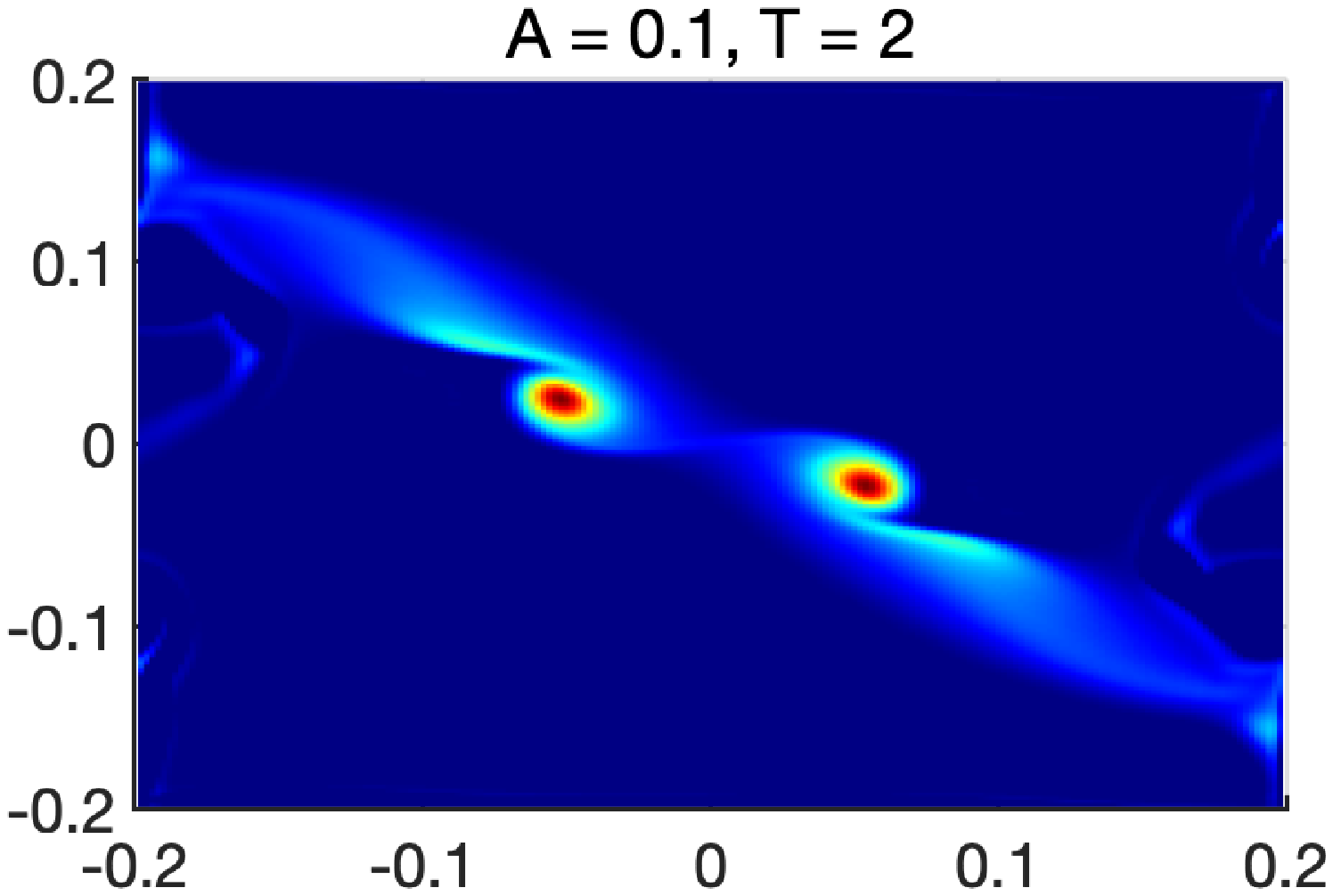}
\includegraphics[width=.3\textwidth]{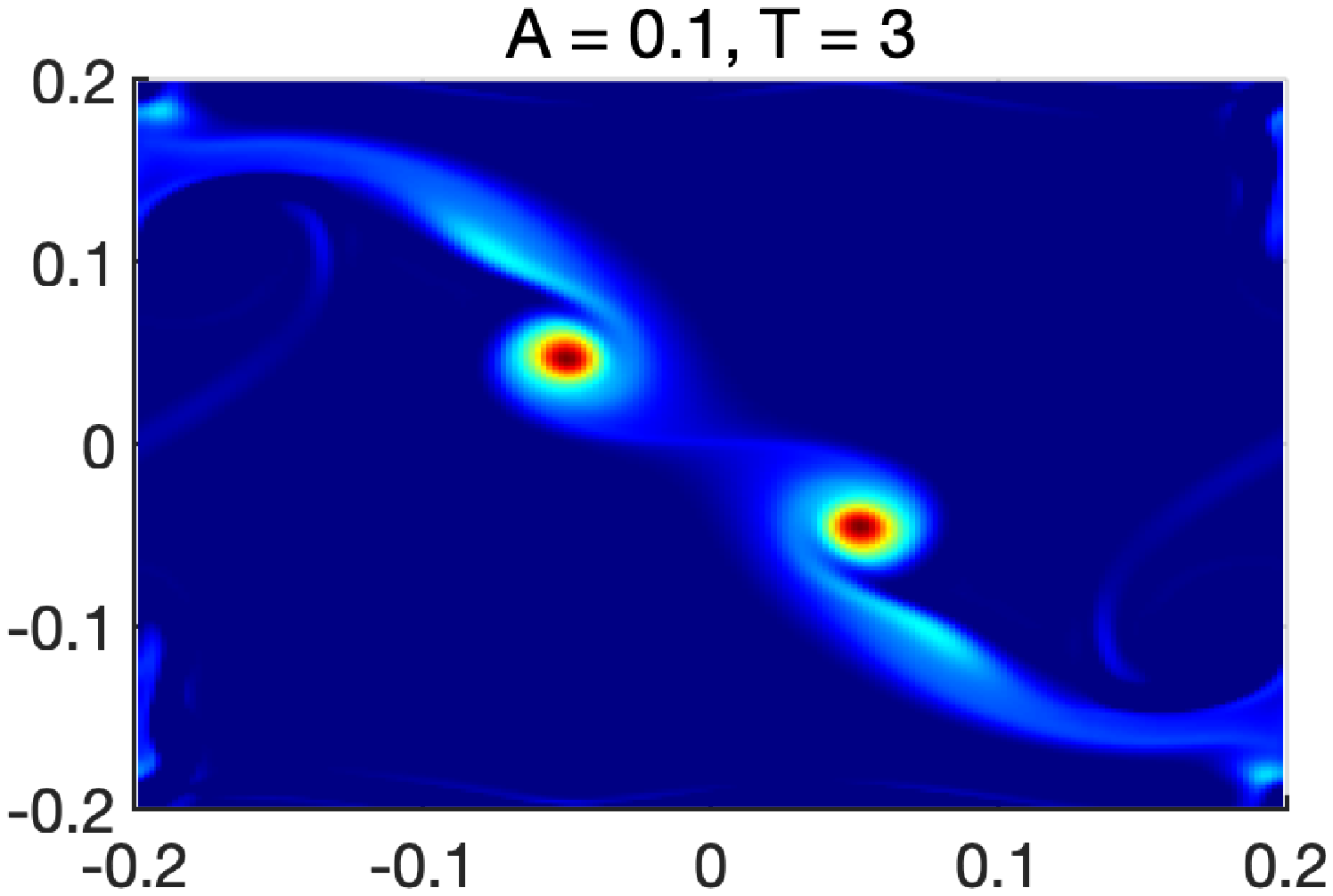}
\caption{Case 2 vorticity profiles. Very compressible system.}
\label{A01}
\end{figure}

\begin{figure}[!htbp]
\centering
\includegraphics[width=.3\textwidth]{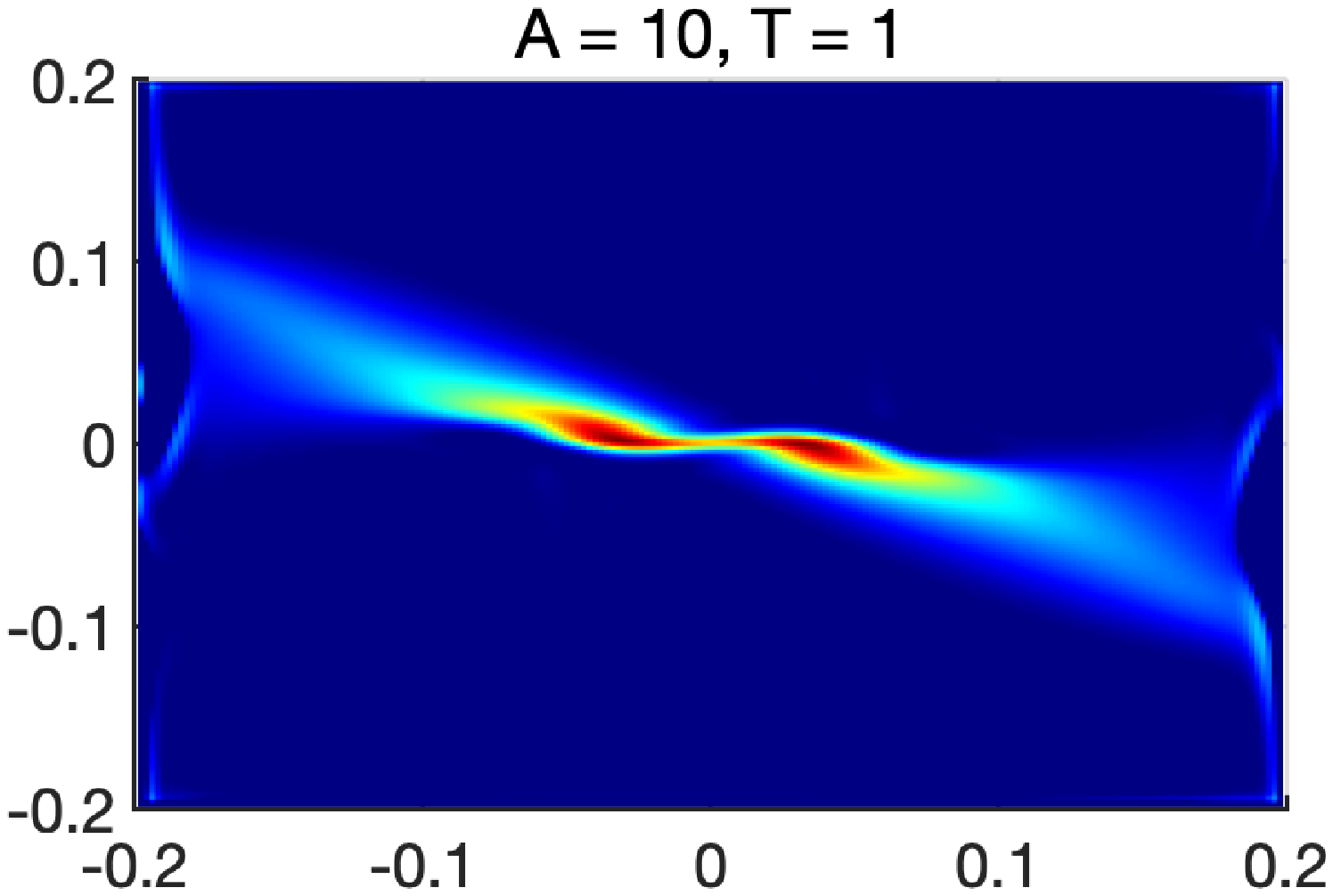}
\includegraphics[width=.3\textwidth]{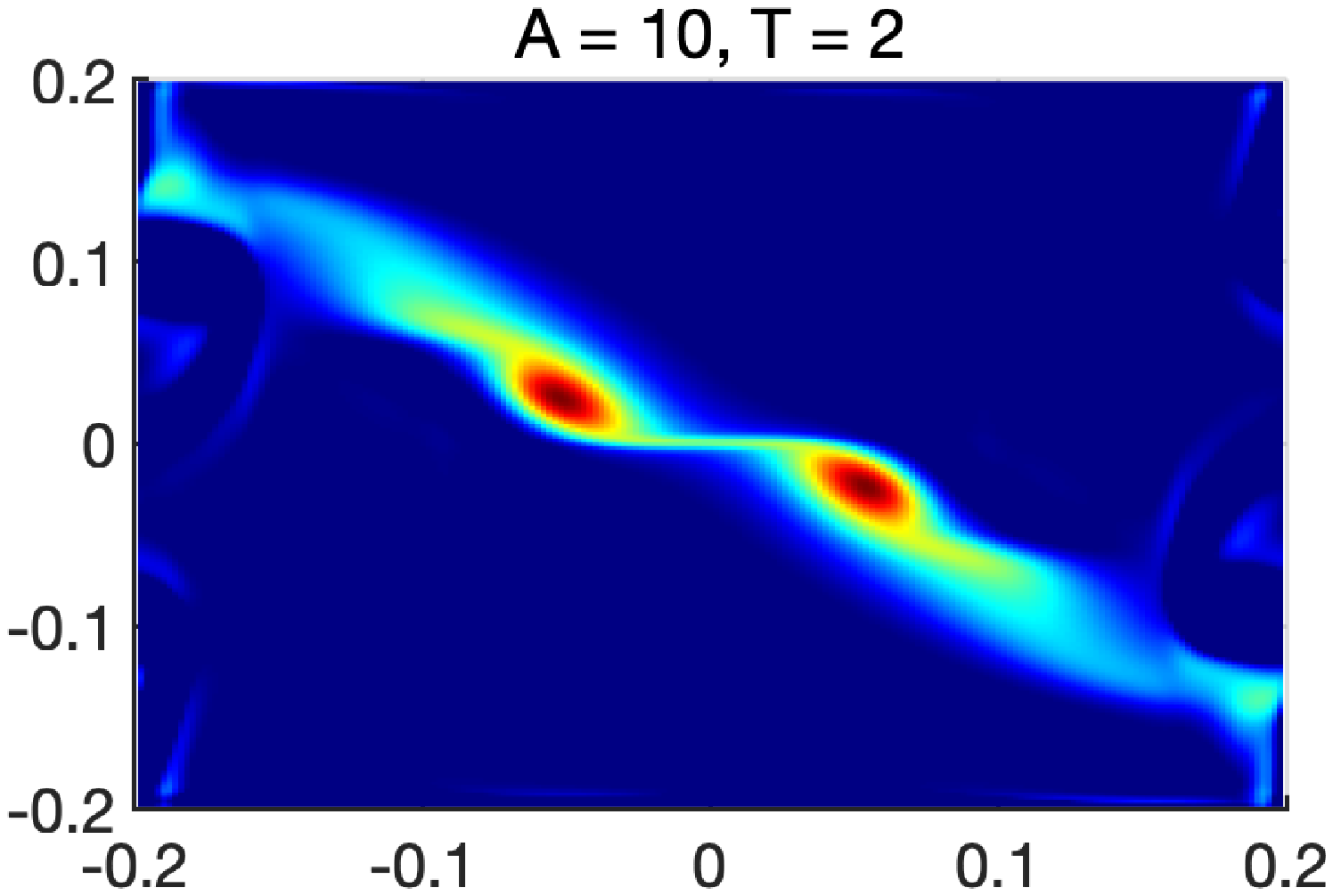}
\includegraphics[width=.3\textwidth]{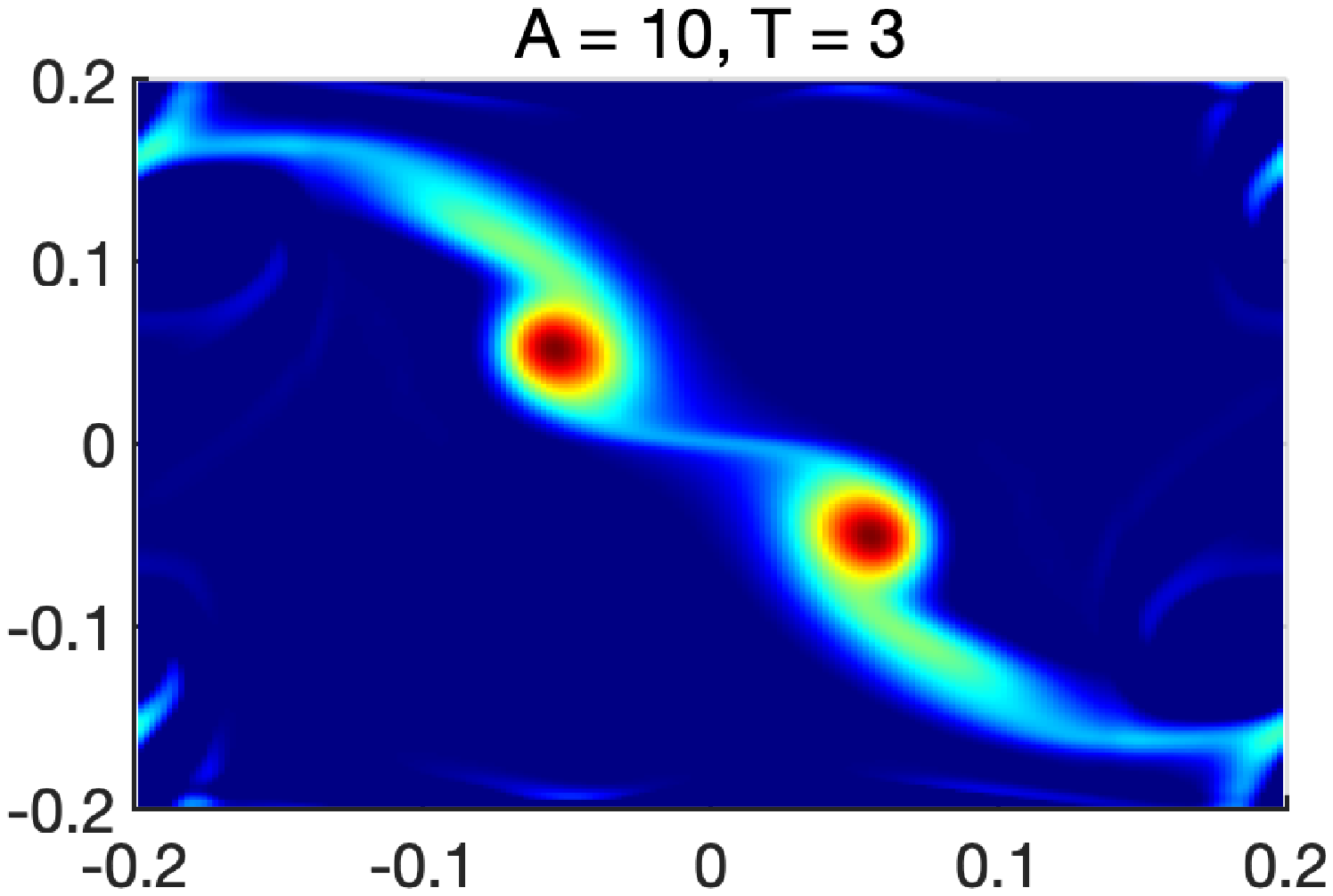}
\caption{Case 2 vorticity profiles. Moderately compressible system.}
\label{A10}
\end{figure}

\begin{figure}[!htbp]
\centering
\includegraphics[width=.3\textwidth]{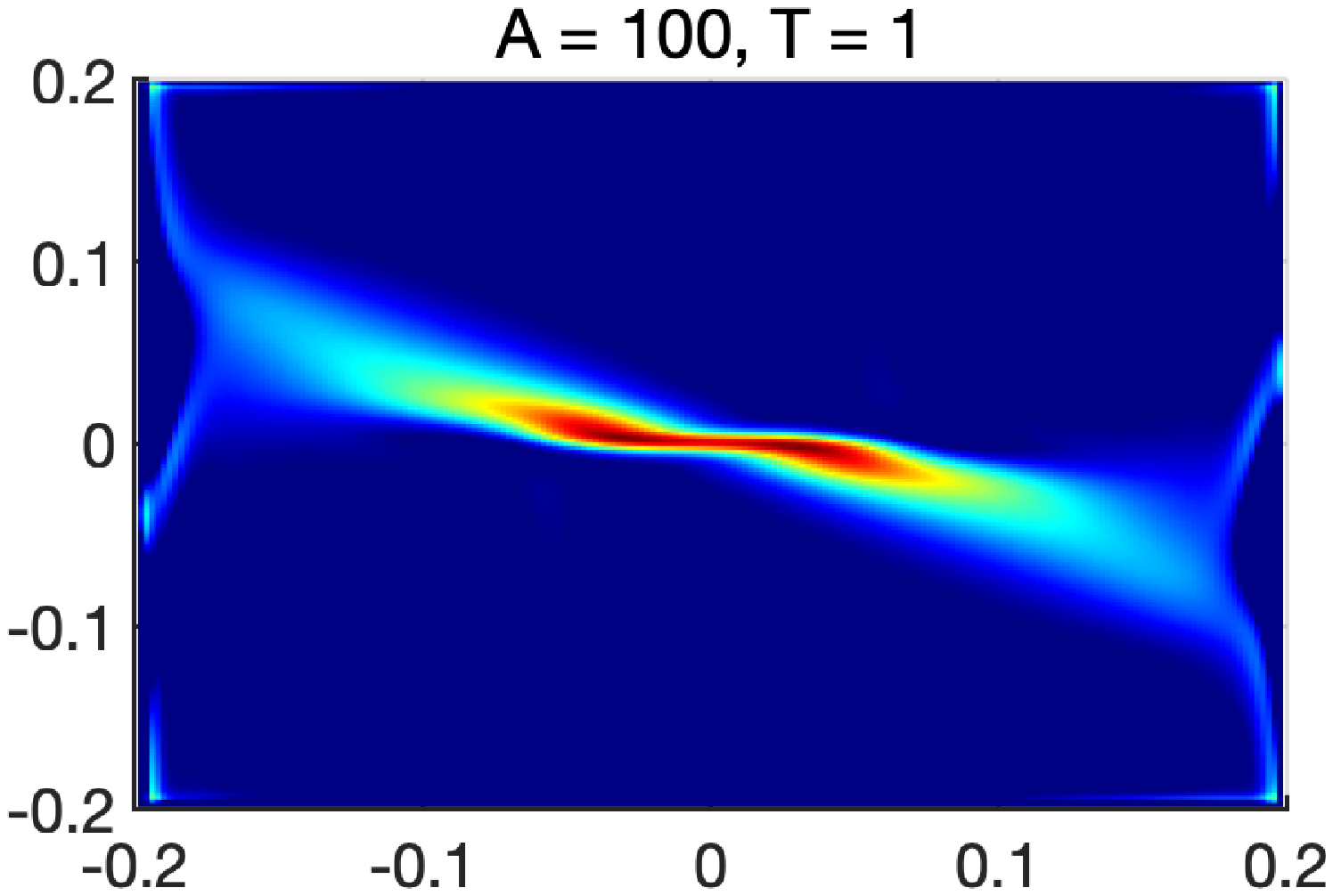}
\includegraphics[width=.3\textwidth]{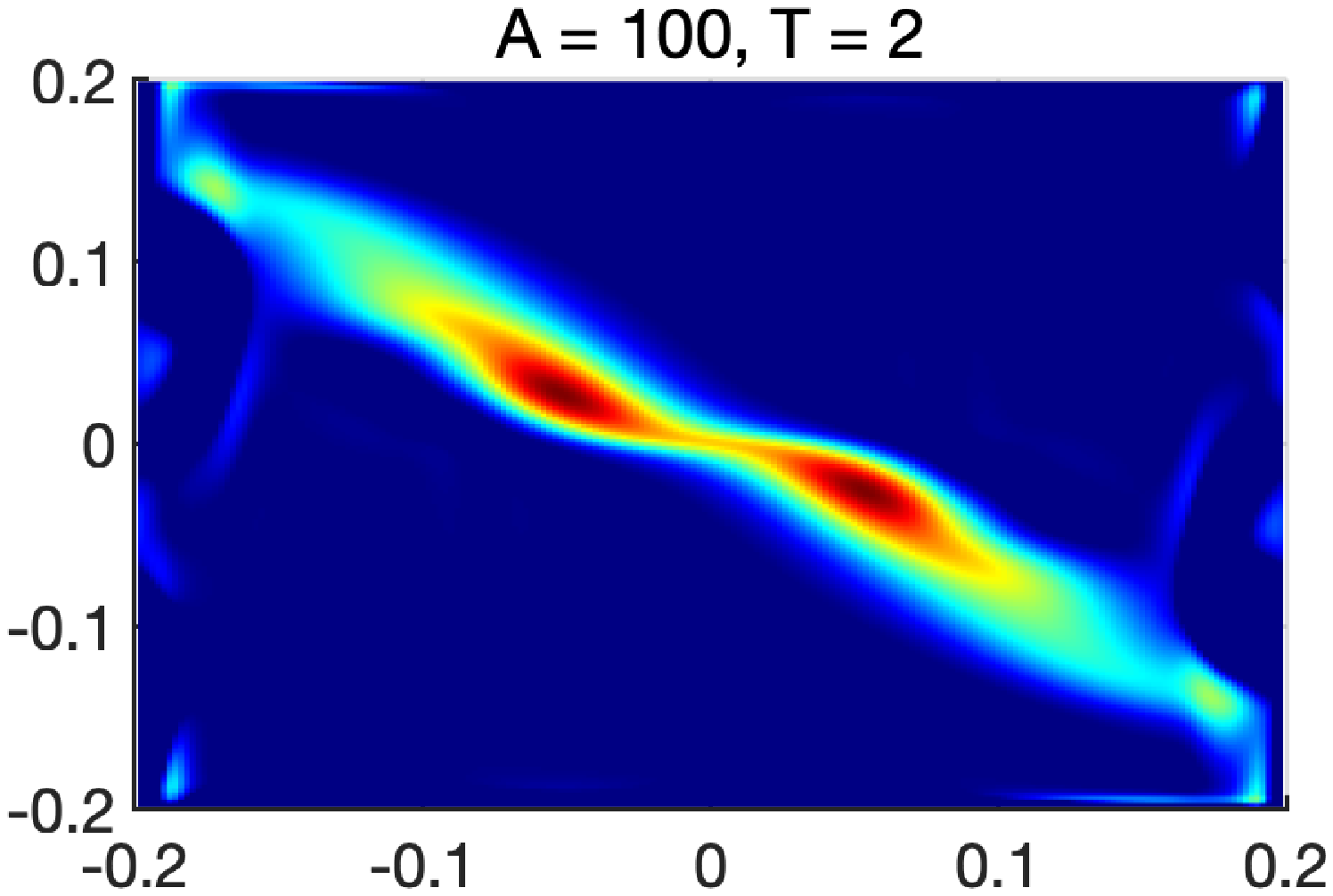}
\includegraphics[width=.3\textwidth]{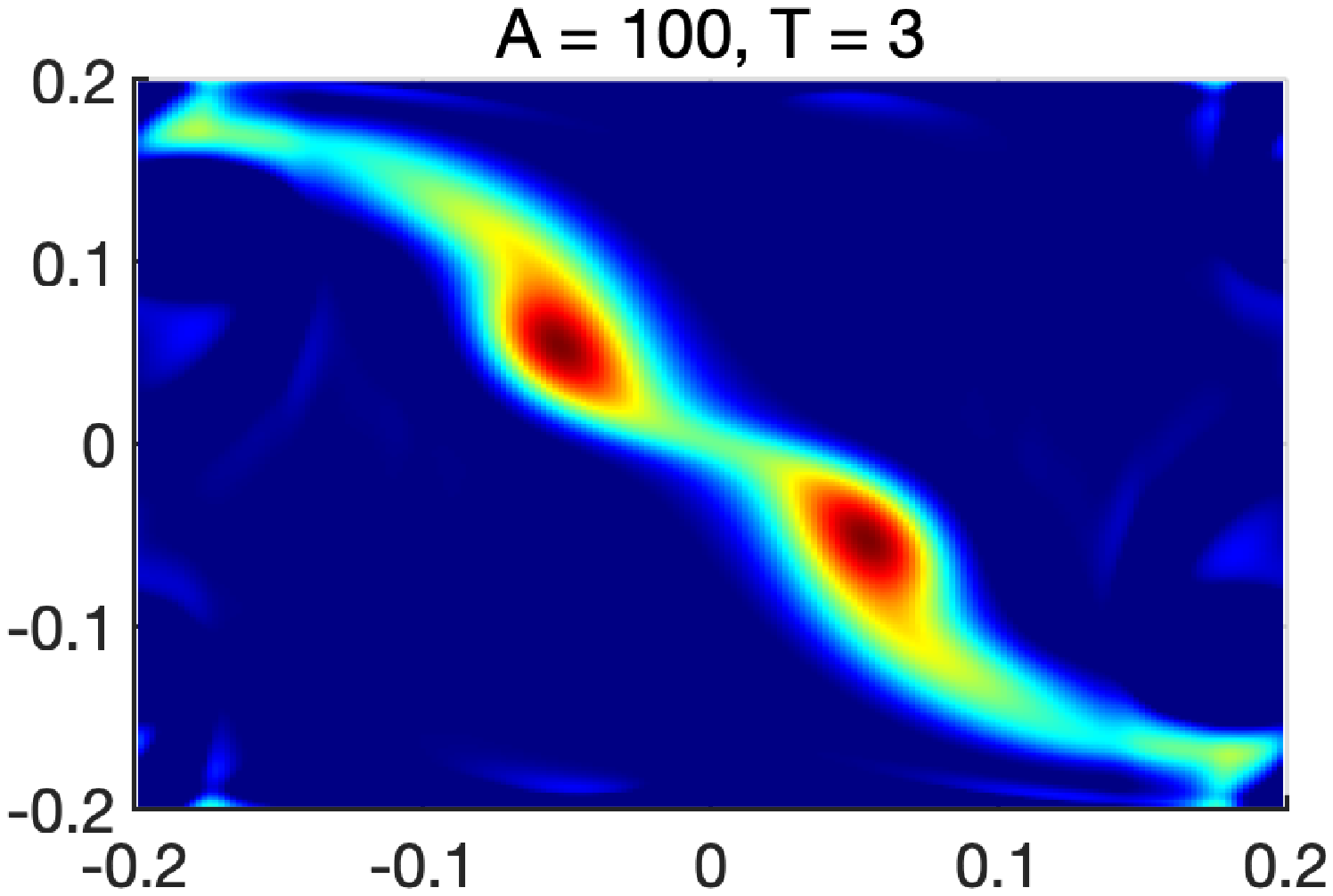}
\caption{Case 2 vorticity profiles. Mildly compressible system.}
\label{A100}
\end{figure}

\begin{figure}[!htbp]
\centering
\includegraphics[width=.3\textwidth]{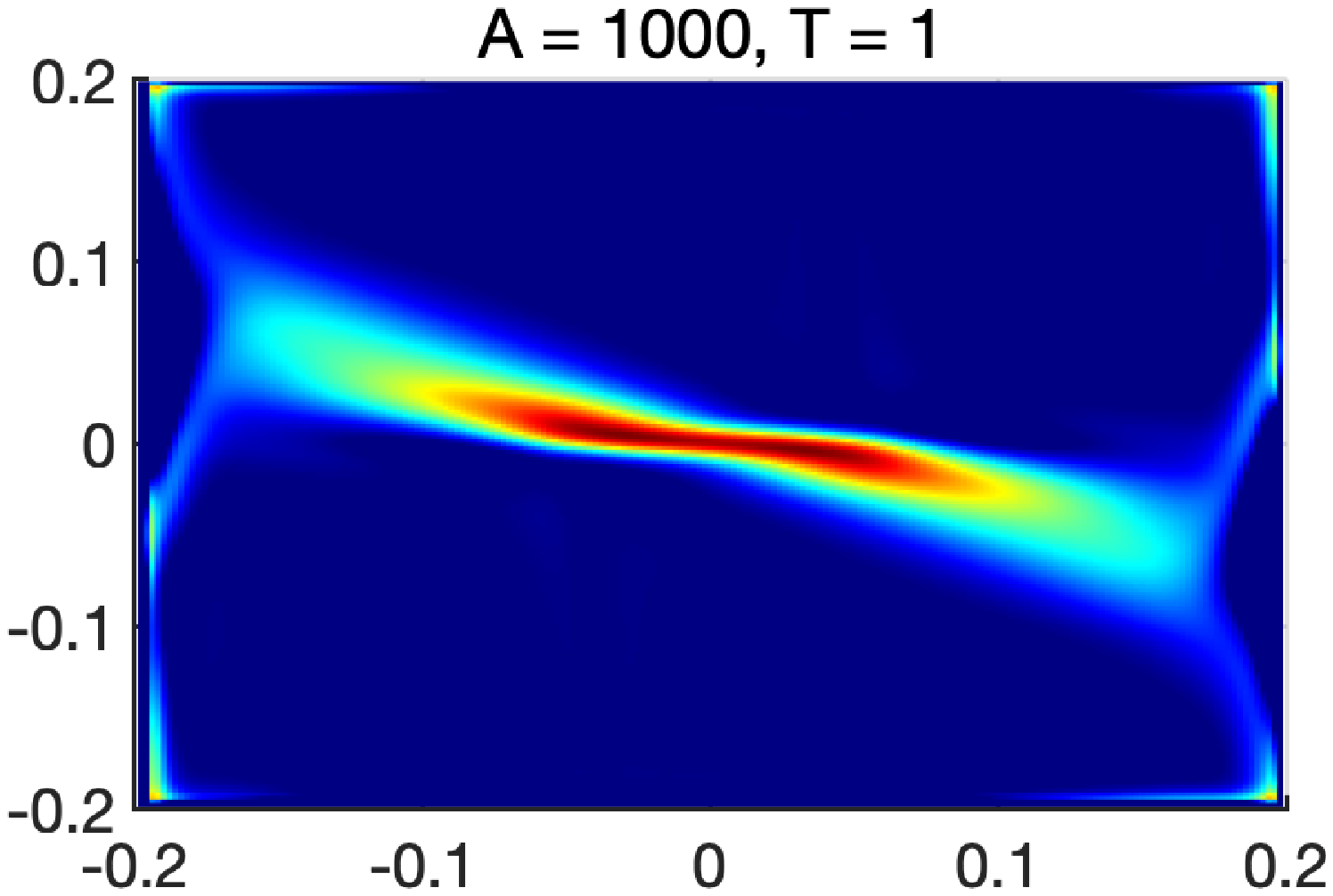}
\includegraphics[width=.3\textwidth]{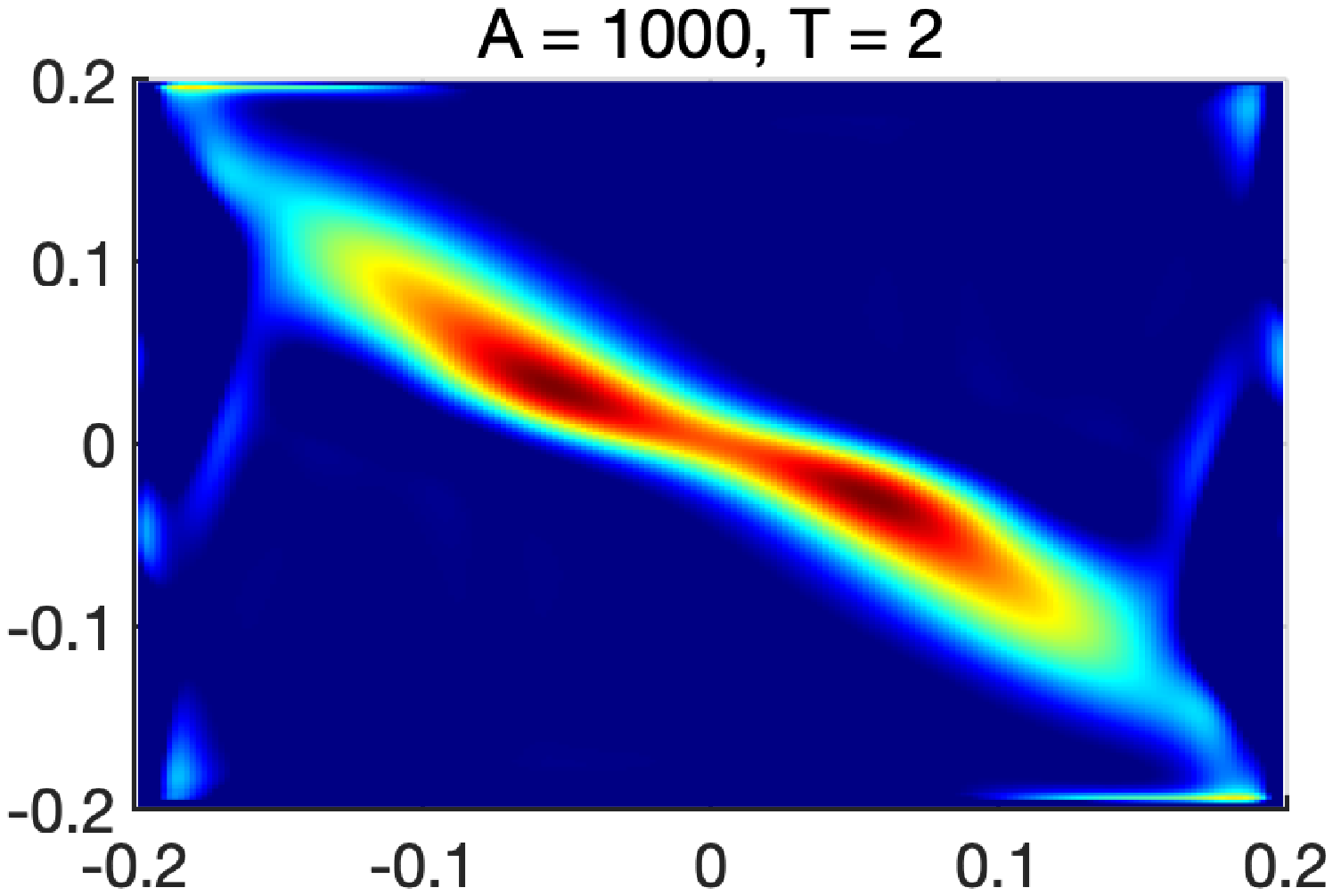}
\includegraphics[width=.3\textwidth]{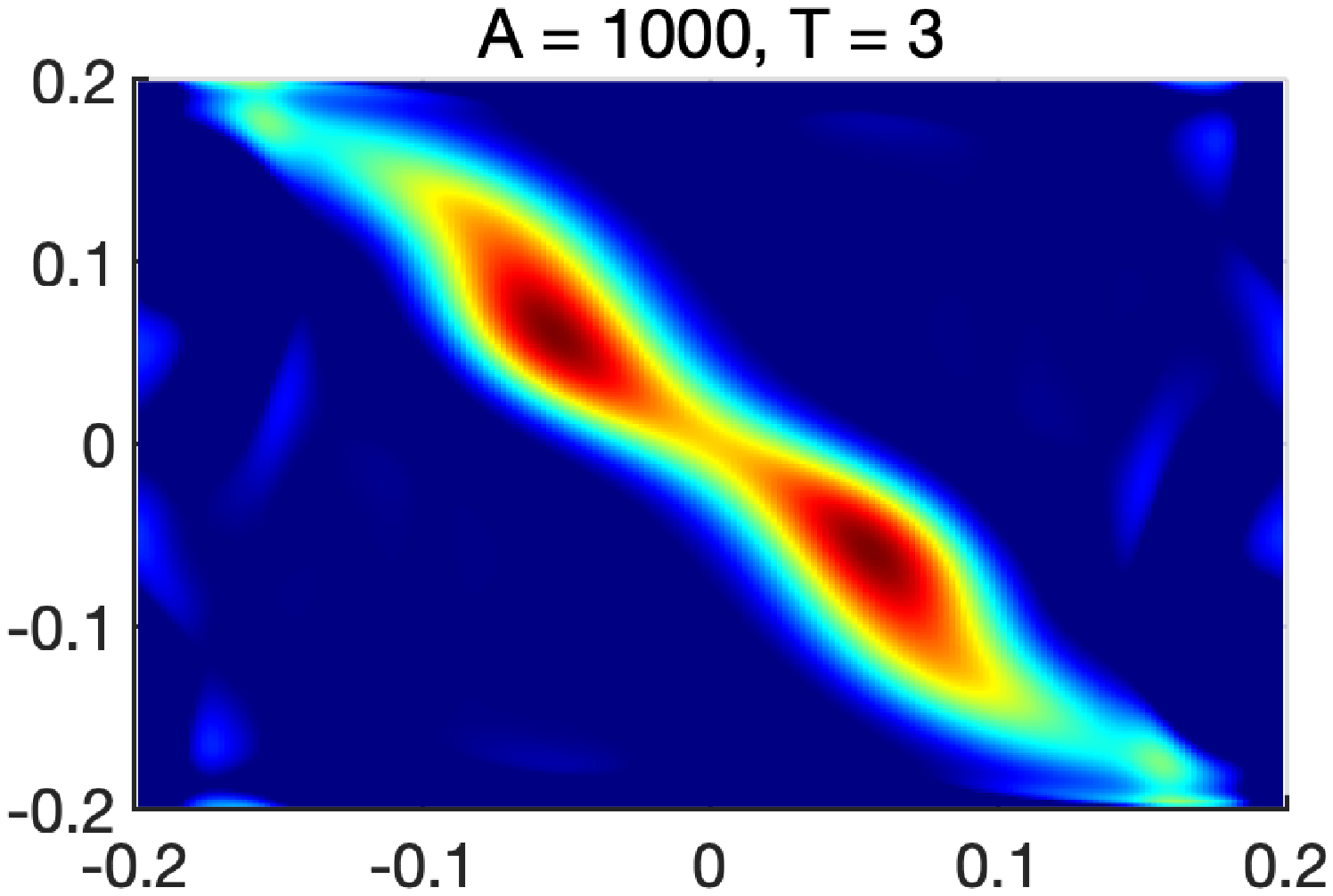}
\caption{Case 2 vorticity profiles. Slightly compressible system.}
\label{A1000}
\end{figure}

In light of the above discussion on the incompressible fluid limit,  the compressibility of the system can be enhanced by decreasing the value of $A$.  The resulting Case 2 vorticity profiles at different times  are presented from  Figure \ref{A001} to Figure \ref{A1000}. We can see that the two spirals are formed more slowly when the system is getting less compressible. 

\section{Concluding remarks}
Hyperbolic conservation laws provide the basic mathematical models for continuum physics, widely used in the scientific and engineering community. Yet, for a long time a general existence-uniqueness theorem  in several space dimensions has awaited a rigorous justification. With the numerical simulations
presented in this paper, we hope to raise the awareness that this lack of a well-posedness theory reflects a fundamental obstruction stemming from the very nature of the equations. 
At an intuitive level, when the initial vorticity is supported on two wedges and has a power singularity 
at the origin, the mechanism leading to multiple solutions can be easily understood.
This lack of uniqueness is indeed confirmed by several of our computations.
It remains a challenging open problem to rigorously validate these simulations, proving the existence of exact 
solutions having the same structure as the numerically computed ones.

\section*{Acknowledgments} 
The authors would like to appreciate the associate editor and two anonymous referees for their constructive comments that have improved the presentation of this paper. The research  of A.~Bressan was partially supported by NSF with grant  DMS-2006884, ``Singularities and error bounds for hyperbolic equations". Liu's research was partially supported by NSF with Grant DMS1812666.

\bigskip

\bibliographystyle{plain}

\end{document}